\documentclass{amsart}

\usepackage{amsmath}
\usepackage{amsthm}
\usepackage{amssymb}
\usepackage{amscd}

\newtheorem{thm}{Theorem}[section]
\newtheorem{prop}[thm]{Proposition}
\newtheorem{lemma}[thm]{Lemma}
\newtheorem{cor}[thm]{Corollary}
\newtheorem{conj}{Conjecture}

\theoremstyle{definition}
\newtheorem{defi}[thm]{Definition}
\newtheorem{ex}[thm]{Example}

\newtheorem{problem}[thm]{Problem}

\theoremstyle{remark}
\newtheorem{rmk}[thm]{Remark}
\newtheorem*{acknow}{Acknowledgment}

\numberwithin{equation}{section}

\newcommand{\ra}{\mbox{$\rightarrow$}}
\newcommand{\hra}{\mbox{$\hookrightarrow$}}
\newcommand{\proj}{\mathop{\rm Proj}\nolimits}

\newcommand{\im}{\mathop{\rm Im}\nolimits}
\newcommand{\Pic}{\mathop{\rm Pic}\nolimits}
\newcommand{\Bs}{\mathop{\rm Bs}\nolimits}
\newcommand{\Sing}{\mathop{\rm Sing}\nolimits}
\newcommand{\codim}{\mathop{\rm codim}\nolimits}
\newcommand{\Tor}{\mathop{\rm Tor}\nolimits}

\newcommand{\rank}{\mathop{\rm rank}\nolimits}
\newcommand{\Chow}{\mathop{\rm Chow}\nolimits}
\newcommand{\Hom}{\mathop{\rm Hom}\nolimits}
\newcommand{\RatCurves}{\mathop{\rm RatCurves}\nolimits}
\newcommand{\Univ}{\mathop{\rm Univ}\nolimits}

\begin{document}

\title[Classification of 
generalized polarized manifolds]{Classification 
of generalized polarized manifolds
by their nef values}

\author{Masahiro Ohno}
\address{Department of System Engineering\\
The University of Electro-Communications\\
1-5-1 Chofugaoka Chofu-shi Tokyo\\
182-8585 Japan}
\email{ohno@e-one.uec.ac.jp}



\begin{abstract}
Let $(M,\mathcal{E})$ be a generalized polarized manifold,
i.e., a pair of an $n$-dimensional smooth projective variety
and an ample vector bundle $\mathcal{E}$ of rank $r$ on $M$.
Let $\tau$ be 
the nef value of a polarized manifold $(M,\det \mathcal{E})$,
i.e., the minimum of the set of real numbers $t$
such that $K_M+t\det\mathcal{E}$ is nef;
we have $\tau r\leq n+1$ by 
Mori's theory.
In this paper we classify the pairs $(M,\mathcal{E})$
in the following two cases:
(1) $n-2\leq \tau r$ and $\tau \geq 1$;
(2) $n+1-\tau r<\tau \leq 1$.
\end{abstract}

\maketitle

\section{Introduction}
Let $M$ be an $n$-dimensional complex projective manifold
and $\mathcal{E}$ an ample vector bundle of rank $r$ on $M$.
The numerical effectiveness of the adjoint bundle
$K_M+\det \mathcal{E}$ has been studied by many authors
(\cite{yz}, \cite{p0}, \cite{fn}, \cite{p}, \cite{psw}, \cite{abw},
\cite{ma1}, and \cite{am} for the case $M$ is smooth 
and $\mathcal{E}$ is a vector bundle).
For a polarized manifold $(M,L)$
of arbitrary dimension $n$,
Ionescu~\cite{i} and Fujita~\cite{fs} first succeeded,
due to Mori-Kawamata theory,
to determine
when the adjoint bundle $K_M+rL$ is nef or not
in case $r\geq n-2$;
this corresponds, in our setting, to the case where $\mathcal{E}$
is the direct sum $L^{\oplus r}$ of $r$ copies of 
the ample line bundle $L$.
The first result for a general ample vector bundle $\mathcal{E}$ was 
made by Ye and Zhang \cite[Theorem 1]{yz} 
and Peternell \cite[Theorem]{p0} in case $r=n+1$,
and the research has proceeded to a lower rank case;
the present state of the research
is for $r=n-2$,
made by Maeda \cite{ma1}
and Andreatta and Mella \cite{am}.
By virtue of these researches,
there are lists of the pairs $(M, \mathcal{E})$
whose adjoint bundles $K_M+\det \mathcal{E}$ are, e.g., not nef,
or nef but not ample,
for every $r$ ($n+1\geq r\geq n-2$).

If we look at these lists, we find that 
the same type of the pairs $(M, \mathcal{E})$,
for example $(\mathbb{P}^n, \mathcal{O}(1)^{\oplus r})$,
appears repeatedly 
and 
the pairs in the same type differ only in the ranks of their bundles.
This indicates that 
in order to avoid this recurrence
and to obtain a more rank-independent list of the pairs $(M, \mathcal{E})$
we should classify the pairs $(M, \mathcal{E})$
by the nef values $\tau (M, \det \mathcal{E})$
of the polarized manifolds $(M, \det \mathcal{E})$;
here, by the word {\it nef value} $\tau (M, \det \mathcal{E})$,
we mean the minimum of the set of real numbers $t$
such that $K_M+t\det \mathcal{E}$ is nef.
In this paper, 
we study $(M, \mathcal{E})$
from this view point 
and give a more rank-free classification of the pairs $(M, \mathcal{E})$.
The precise results are as follows;
set $\tau =\tau (M,\det\mathcal{E})$.
Then Mori's theorem~\cite[(1.4)]{mo2} first implies the following 

\begin{prop}\label{upper bound}
We have 
$\tau r\leq n+1$.
\end{prop}

Now we have the following classification.

\begin{prop}\label{n+1/r}
If 
$\tau r=n+1$,
then $(M,\mathcal{E})$ is isomorphic to 
$(\mathbb{P}^n, \mathcal{O}(1)^{\oplus r})$.
\end{prop}
In \cite[\S 3 Theorem 2]{a}, Andreatta
has given a proof of Proposition~\ref{n+1/r}.

\begin{thm}\label{Main theorem}
Suppose that $\tau \geq 1$.
If $n-2\leq \tau r$,
then $(M, \mathcal{E})$ and the value of 
$\tau r$ are one of the following:
\begin{enumerate}
\item $(\mathbb{P}^n, 
         \mathcal{O}(1)^{\oplus r})$, and $\tau r=n+1$;
\item $({\mathbb Q}, \mathcal{O}_{{\mathbb Q}}(1)^{\oplus r})$, 
          where ${\mathbb Q}$
         is a hyperquadric in $\mathbb{P}^{n+1}$,
         and $\tau r=n$;
\item $(\mathbb{P}^n, T_{\mathbb{P}^n})$, and $\tau r=n$;
\item $(\mathbb{P}^n, \mathcal{O}(1)^{\oplus (r-1)}\oplus \mathcal{O}(2))$,
      and $\tau r=(n+1)r/(r+1)\geq n-2$ (and hence $r\geq (n-2)/3$);
\item $(\mathbb{P}({\mathcal F}), 
         H({\mathcal F})\otimes \psi ^*{\mathcal G})$,
         where ${\mathcal F}$ is a vector bundle of rank $n$ on a smooth
         proper curve $C$, $\psi :\mathbb{P}({\mathcal F})\ra C$  
         the projection, and ${\mathcal G}$ an ample vector bundle 
         of rank $r$ on $C$,
         and $\tau r=n$;
\item $M$ is a Del Pezzo manifold with $\Pic M\cong \mathbb{Z}$,
         and $\mathcal{E}\cong A^{\oplus r}$ where
         $A$ is the ample generator of $\Pic M$,
         and $\tau r=n-1$;
\item $r=n-1$, the Picard number $\rho (M)$ of $M$ is one,
         and $K_M+\det \mathcal{E}=0$ (see \cite{psw}
         for a precise classification),
         and $\tau r=n-1$;
\item  There exist a hyperquadric fibration 
         $\psi :M\ra C$ of the relative Picard number one
         over a smooth curve $C$,
         a $\psi$-ample line bundle $\mathcal{O}_M(1)$ on $M$
         and an ample vector bundle ${\mathcal G}$ of rank $r$ on $C$
         such that $\mathcal{E}\cong 
         \mathcal{O}_M(1)\otimes \psi^* {\mathcal G}$
         where $\mathcal{O}_M(1)|_F\cong \mathcal{O}_{\mathbb{Q}}(1)$
         for any fiber $F\cong \mathbb{Q}$ of $\psi$,
         and $\tau r=n-1$;
\item There exists a $\mathbb{P}^{n-1}$-bundle
          $\psi :M\ra C$ over a smooth curve $C$ such that
          $\mathcal{E}|_F\cong 
         T_{{\mathbb P}^{n-1}}$ for any fiber $F$ of $\psi$,
          and $\tau r=n-1$;
\item  $M$ is isomorphic to a projective space bundle
         $\mathbb{P}(\mathcal{F})$ over a smooth proper curve $C$
         for some vector bundle $\mathcal{F}$ of rank $n$ on $C$,
         and there is an exact sequence
         \[0\ra \pi^*\mathcal{L}\otimes H(\mathcal{F})^{\otimes 2}
            \ra \mathcal{E}\ra \pi^*\mathcal{G}\otimes H(\mathcal{F})
            \ra 0\]
         for some line bundle $\mathcal{L}$ on $C$
         and some vector bundle $\mathcal{G}$ of rank $r-1\geq 0$ 
         on $C$
         where $\pi :\mathbb{P}(\mathcal{F})\ra C$ is the 
         projection, and $\tau r=nr/(r+1)\geq n-2$
         (and hence $r\geq (n-2)/2$);
\item There exists a $\mathbb{P}^{n-2}$-fibration
         $\psi :M\ra S$,
         locally trivial in the complex (or \'{e}tale) topology,
         over a smooth surface $S$ such that 
         $\mathcal{E}|_{F}\cong 
         \mathcal{O}_{\mathbb{P}^{n-2}}(1)^{\oplus r}$
         for every fiber $F$ of $\psi$,
         and $\tau r=n-1$;
\item $M$ is the blowing-up $\psi :M\ra M'$
         of a projective manifold $M'$ at finite points,
         and there exists an ample vector bundle 
         $\mathcal{E}'$ of rank $r$ on $M'$ such that
         $\mathcal{E}\cong \psi ^*\mathcal{E}'\otimes
         \mathcal{O}_M(-E)$ where $E$ is the exceptional divisor
         of $\psi$
         and that $(M',\mathcal{E}')$ satisfies one of the following:
     \begin{enumerate}
         \item $\tau (M',\det \mathcal{E}')r<n-1$;
         \item $M'\cong \mathbb{P}^3$, and 
               $\mathcal{E}'\cong \mathcal{O}(2)^{\oplus 2}$
               or $\mathcal{O}(2)$;
         \item $(M',\mathcal{E}'),\cong (\mathbb{P}^2,\mathcal{O}(3))$;
         \item $n=2$, $(M',\mathcal{E}')$ is of type $10)$ above,
         and $\tau r=n-1$;
     \end{enumerate}
\item $(\mathbb{P}^n, 
         \mathcal{O}(1)^{\oplus (r-1)}\oplus \mathcal{O}(3))$,
         and $\tau r=(n+1)r/(r+2)\geq n-2$
         (and hence $r\geq 2(n-2)/3$);
\item $(\mathbb{P}^n, \mathcal{O}(1)^{\oplus (r-2)}\oplus 
         \mathcal{O}(2)^{\oplus 2})$,
         and $\tau r=(n+1)r/(r+2)\geq n-2$
         (and hence $r\geq 2(n-2)/3$);
\item $(\mathbb{Q}^n, \mathcal{O}_{\mathbb{Q}}(1)^{\oplus (r-1)}
         \oplus \mathcal{O}(2))$ 
         where $\mathbb{Q}^n$
         is a hyperquadric in $\mathbb{P}^{n+1}$,
         and $\tau r=nr/(r+1)\geq n-2$
         (and hence $r\geq (n-2)/2$);
\item $(\mathbb{Q}^4, \mathbf{E}(2))$ 
         where $\mathbb{Q}^4$
         is a hyperquadric in $\mathbb{P}^{5}$
         and $\mathbf{E}$ is the spinor bundle on $\mathbb{Q}^4$,
         and $\tau r=nr/(r+1)=8/3$;
\item $M$ is a Fano manifold of Picard number one
         with $K_M+(n-2)A=0$, and $\mathcal{E}\cong A^{\oplus r}$
         where $A$ is an ample line bundle on $M$,
         and $\tau r=n-2$;
\item $r=n-2$, the Picard number 
          of $M$ is one,
         and $K_M+\det \mathcal{E}=0$,
         and $\tau r=n-2$;
\item $M$ has an elementary contraction $\psi:M\to C$ onto a curve $C$
         such that 
         a general fiber $F$ of $\psi$ is a Del Pezzo manifold,
         i.e., $K_F+(n-2)A=0$ for some ample line bundle $A$ on $F$,
         and that 
         $\mathcal{E}|_F\cong A^{\oplus r}$,
         and $\tau r=n-2$;
\item $r=n-2$ and $M$ has an elementary contraction 
         $\psi :M\ra C$ onto a smooth proper curve $C$ such that
         $K_F+\det \mathcal{E}|_F=0$
         for a general fiber $F$ of $\psi$,
         and $\tau r=n-2$;
\item $M$ has an elementary contraction 
         $\psi :M\ra S$ onto a smooth projective surface $S$ 
         such that a general fiber $F$ of $\psi$ is 
         hyperquadric
         and $\mathcal{E}|_F\cong \mathcal{O}(1)^{\oplus r}$,
         and $\tau r=n-2$;
\item $r=n-2$ and there exists a $\mathbb{P}^{n-2}$-fibration
          $\psi :M\ra S$,
         locally trivial in the \'{e}tale (or complex) topology,
         over a smooth projective surface $S$ 
         such that $\mathcal{E}|_F\cong 
         T_{\mathbb{P}^{n-2}}$ for any closed fiber $F$ of $\psi$,
         and $\tau r=n-2$;
\item $r=n-2$ and there exists a $\mathbb{P}^{n-2}$-fibration
         $\psi :M\ra S$,
         locally trivial in the \'{e}tale (or complex) topology,
         over a smooth surface $S$ such that 
         $\mathcal{E}|_{F}\cong 
         \mathcal{O}
         (1)^{\oplus (n-3)}
         \oplus \mathcal{O}(2)$
         for any closed fiber $F$ of $\psi$,
         and $\tau r=n-2$;
\item $M$ has an elementary contraction 
         $\psi :M\ra S$ onto a normal projective $3$-fold $S$ 
         with only rational Gorenstein singularities
         such that except for a finite number of fibers
         each fiber $F$ of $\psi$ is isomorphic to $\mathbb{P}^{n-3}$
         and $\mathcal{E}|_F\cong \mathcal{O}(1)^{\oplus r}$
         and $\tau r=n-2$;
\item $r=n-2$, $M$ is the blowing-up $\psi\colon M\to M'$
          of a projective manifold $M'$ at a point,
          and $\mathcal{E}$ fits into the following exact sequence
         \[0\to \psi^*\mathcal{E}'\otimes \mathcal{O}_M(-2E)\to
           \mathcal{E}\to \mathcal{O}_E(-E)^{\oplus (n-3)}\to 0,\]
           where $\mathcal{E}'$ is a vector bundle of rank $n-2$ on $M'$
           and $E$ is the exceptional divisor of $\psi$,
          and $\tau r=n-2$;
\item  $M$ has a divisorial elementary contraction 
         $\psi :M\ra M'$ onto a projective variety $M'$ 
         such that the exceptional divisor $E$ of $\psi$
         and $\mathcal{E}|_{E}$ satisfy one of the following:
      \begin{enumerate}
         \item $E\cong \mathbb{P}^{n-1}$, 
            $\mathcal{O}_E(E)=\mathcal{O}(-2)$ and 
            $\mathcal{E}|_{E}\cong \mathcal{O}(1)^{\oplus r}$
         \item $E\cong \mathbb{Q}^{n-1}$ (possibly singular hyperquadric),
            $\mathcal{O}_E(E)=\mathcal{O}(-1)$ and 
            $\mathcal{E}|_{E}\cong \mathcal{O}(1)^{\oplus r}$,
         and $\tau r=n-2$;
      \end{enumerate}
\item $M$ is the blowing-up $\psi :M\ra M'$
         of a projective manifold $M'$ along a smooth curve $B$
         and $\mathcal{E}\cong \psi ^*\mathcal{E}'\otimes
         \mathcal{O}_M(-E)$
         for some vector bundle 
         $\mathcal{E}'$ of rank $r$ on $M'$ 
         where $E$ is the exceptional divisor
         of $\psi$,
         and $\tau r=n-2$.
\end{enumerate}
\end{thm}

Now let us recall the previous results more precisely 
in connection with Theorem~\ref{Main theorem}. 
Note here that saying that $K_M+t\det \mathcal{E}$
is not nef, nef but not ample, or ample
is equivalent to saying that $t<\tau $,
$\tau =t$, or $\tau<t$ respectively.
Note also that 
if $K_M+t\det \mathcal{E}$ is nef for some positive rational number
$t$ then $K_M+t\det \mathcal{E}$ is semi-ample by 
the Basepoint-free Theorem~\cite{kmm};
thus to be semi-ample and to be nef are
equivalent for $K_M+t\det \mathcal{E}$ $(t\in \mathbb{Q}_{>0})$.

Suppose that $r=1$.
Ionescu~\cite[Theorem (1.1)]{i}
and Fujita~\cite[Theorem~1]{fs} first gave 
the bound $\tau$($=\tau r$)$\leq n+1$.
Ionescu classified the pairs 
$(M,\mathcal{E})$ in case  $n\leq \tau \leq n+1$
\cite[Theorem (1.2), (1.3), (1.4)]{i},
in case $n-1<\tau \leq n+1$ and $n\geq 2$
(and thus $\tau >1$)
\cite[Theorem (1.5), (1.6)]{i},
and in case $n-2<\tau \leq n+1$ and $n\geq 3$
modulo reduction
\cite[Theorem (1.7)]{i}.
Fujita classified the pairs 
$(M,\mathcal{E})$
in case 
$n-1<\tau \leq n+1$
and $M$ has some mild singularities
\cite[Theorems~1 and 2]{fs},
in case 
$n-2<\tau \leq n-1$ and $n\geq 3$
\cite[Theorem 3']{fs} (\cite[(11.8)]{fb}),
in case $n-3<\tau \leq n-2$ and $n\geq 4$
\cite[Theorem 4]{fs}.

If $r=n+1$, then  
Ye and Zhang \cite[Theorem~1]{yz} showed
that $\tau \leq 1$ ($r=n+1$ and $\tau \leq 1$ implies $\tau r\leq n+1$),
and classified the pairs
$(M,\mathcal{E})$ with $K_M+\det\mathcal{E}=0$ 
(and thus $\tau =1$ and $\tau r=n+1$).
Peternell \cite[Theorem]{p0}
also classified the case $r= n+1$
and $K_M+\det\mathcal{E}=0$.
Fujita~\cite[Main Theorem]{fn}
classified the pairs $(M,\mathcal{E})$ in case 
$r\geq n$ and $\tau \geq 1$ 
(this is a sub-case of the case $\tau r\geq n$).
Peternell~\cite[Theorem~2]{p} also classified
in case $r=n$ and $K_M+\det\mathcal{E}=0$
(a sub-case of the case $\tau=1$ and $n=\tau r$).
If $r=n-1$ and $\tau =1$ (a sub-case of $\tau r=n-1$),
Peternell, Szurek, and Wi\'{s}niewski~\cite[Main theorem(0.3)]{psw}
classified the pairs $(M,\mathcal{E})$
in case $r=n-1\geq 4$ and $K_M+\det\mathcal{E}=0$,
Wi\'{s}niewski~\cite{w2}
(see also \cite[Theorem(0.4)]{psw})
classified 
the pairs $(M,\mathcal{E})$
in case $r=n-1=2$ and $K_M+\det\mathcal{E}=0$,
and Andreatta, Ballico, and Wi\'{s}niewski~\cite[Theorem]{abw}
classified in case $K_M+\det\mathcal{E}\neq 0$.
Finally, if $r=n-2$ and $\tau =1$ (a sub-case of $\tau r=n-2$),
then Andreatta and Mella~\cite[Theorem 5.1.\ 2) and 3)]{am} 
classified the pairs $(M,\mathcal{E})$
except for the case $K_M+\det \mathcal{E}=0$.

Suppose that $\tau >1$,
i.e., that $K_M+\det \mathcal{E}$ is not nef.
If $r=n$ (a sub-case of $n<\tau r$ under the assumption $\tau >1$), 
then pairs $(M,\mathcal{E})$
are classified by Ye and Zhang~\cite[Theorem~2]{yz}.
If $r=n-1$ (a sub-case of $\tau r>n-1$), 
then pairs $(M,\mathcal{E})$
are classified also by Ye and Zhang~\cite[Theorem~3]{yz}.
Finally, if $r=n-2$ (a sub-case of $n-2<\tau r$),
then pairs $(M,\mathcal{E})$
are classified
by Maeda~\cite[Theorem]{ma1}
and by Andreatta and Mella~\cite[Theorem~5.1.\ 1)]{am}.

Now we give several comments related to 
Theorem~\ref{Main theorem} itself.
First, the proof of Theorem~\ref{Main theorem} basically follows 
the ideas and techniques developed by the above mentioned authors.
However we need some new considerations
to deal with some cases.
For example, I could not help applying the theorem of 
Cho, Miyaoka, and Shepherd Barron~\cite{cms},
and Kebekus~\cite{ke1} \cite{ke2},
characterizing projective spaces in terms of their length,
to deal with, in particular, the case $n-2\leq \tau r<n-1$
(see, e.g., 
Proposition~\ref{l(r1)=n+1new}, or \S~\ref{LC1=2shin}).
Moreover we make the proof of Theorem~\ref{Main theorem}
independent of the value of $\tau r$ as possible,
in particular, in case the Picard number $\rho (M)$ of $M$
is one (see, e.g., \S~\ref{Slr1=n+1}, \S~\ref{Slr1=n},
and \S~\ref{Slr1=n-1}), so that, I believe,
the results obtained in the proof are also 
applicable to the case of $\tau r$ smaller,
e.g., to the case $n-3<\tau r<n-2$.
Second, we can regard $(M',\mathcal{E}')$ in the case (12) (a) of 
Theorem~\ref{Main theorem} as the first reduction
of the generalized polarized manifold $(M,\mathcal{E})$
for all $r\leq n-1$.
(See \cite{bsbook} for the first reduction 
of a polarized variety.)
Third, we did not classify the pairs $(M,\mathcal{E})$ in the case
(18) of Theorem~\ref{Main theorem}, as in 
\cite[Theorem 5.1. 2)]{am}.
Finally, note that the case (25) of Theorem~\ref{Main theorem} 
is ruled out by mistake
in \cite[Theorem 5.1 3)]{am}. 
We will also give an example 
(Example~\ref{counter example of Andreatta Mella})
of this case.
I also remark that I could not determine
the structure of $\psi$
in the case (25) of Theorem~\ref{Main theorem}
without knowing \cite[Theorem~5.1]{ao}
of Andreatta and Occhetta.

For the case $\tau \leq 1$, we have the following
\begin{prop}\label{n+1/r+1}
If $n+1-\tau r<\tau \leq 1$,
then $(M,\mathcal{E})\cong (\mathbb{P}^n, \mathcal{O}(1)^{\oplus r})$.
\end{prop}

Note that 
in the report \cite{o} of my talk at a workshop at RIMS,
I announced Propositions~\ref{upper bound}, \ref{n+1/r},
and \ref{n+1/r+1},
the classification in case $n-1<\tau r$ and $\tau \geq 1$,
and the classification in case
$n-1=\tau r$ and  $2\leq r\leq n-1$.
(Note that $\tau r<s+1$ for some integer $s$ and $\tau \geq 1$
implies $r\leq s$.)
I also gave a proof of Proposition~\ref{upper bound}, 
and, based on a characterization
of projective space by its length
due to \cite{cms},
gave a short proof of 
Propositions~\ref{n+1/r} and \ref{n+1/r+1}.
I also proved shortly 
the classification in case $n-1<\tau r$ and $\tau \geq 1$,
and outlined proofs of the classification
in case 
$n-1=\tau r$ and  $2\leq r\leq n-1$
with several typographical errors and inaccuracies.

In this paper, we extend the classification
to the case $n-2\leq \tau r$ and $\tau \geq 1$
and give a complete proof of the classification
in a more unifying way as possible
(see, e.g., Theorem~\ref{l(r1)=n}, 
Propositions~\ref{l(r1)=n+1new}, \ref{l(r1)=n-1},
\ref{no n-1 fiber}, and \ref{lr1=n-1&n-2fiber}).
\begin{acknow}
The author was a Research Fellow of the Japan Society for 
the Promotion of Science and was supported by Research
Fellowships of the Japan Society for the Promotion of Science
for Young Scientists from April 1st to September 30th 1998.
Propositions~\ref{upper bound}, \ref{n+1/r}, and \ref{n+1/r+1},
and the classification in case $\tau \geq 1$
and $n-1\leq \tau r$ 
were obtained in this period.
The author expresses 
gratitude to Professor Shoichi Kondo for his kind encouragement.
This work was partially supported by 
Grant-in-Aid for JSPS Fellows 10-6464, The Ministry of 
Education, Science, Sports and Culture,
and it was also partially supported by 
the Ministry of Education, Science, Sports and Culture, 
Grant-in-Aid for 
Encouragement of Young Scientists (A), 1274009, 2000.
Deep appreciation also goes to Professor Hajime Kaji,
who pointed out an error of 
the preliminary version of Theorem~\ref{Main theorem},
error related to the case 24) of the theorem.
This work is also partially supported by 
the Ministry of Education, Science, Sports and Culture, 
Grant-in-Aid for Young Scientists (B), 15740006, 2003.
The author expresses 
gratitude to Professor Yoshiaki Fukuma for 
inviting the author
to talk at his conference
in August 2003 at Kochi University,
which motivates
the author to improve the proof of 
the theorem
in a unifying way as possible,
and to Professors Chikashi Miyazaki and Atsushi Noma
for inviting the author
to talk at their mini-workshop
in January 2004 at University of the Ryukyus,
which makes the author to notice a gap
in the proof of Lemma~\ref{key} in the previous draft of this paper.
Deep appreciation also goes to Professor Antonio Lanteri
for sending me some reprints of him (and of him and 
Professor Hidetoshi Maeda)
and for invaluable suggestions
about the previous draft of this paper.
\end{acknow}
\noindent
{\large{{\bf Notation and conventions}}} 
\\
In this paper, we work over 
the complex number field $\mathbb{C}$.
Basically we follow the standard notation and terminology in algebraic
geometry. 
We use the word {\it manifold} to mean a smooth variety.
For a manifold $M$, we denote by $K_M$ 
the canonical divisor of $M$.
We use 
the word {\it line} to mean a smooth rational curve of degree $1$.
We also 
use the words "locally free sheaf\/" and "vector bundle" interchangeably.
For a vector bundle $\mathcal{E}$ on a variety $M$, we denote 
by $H(\mathcal{E})$ 
the tautological line bundle $\mathcal{O}_{\mathbb{P}(\mathcal{E})}(1)$
on $\mathbb{P}(\mathcal{E})$. 
We also use the terminology in the Minimal Model Program. 
For our terminology, we fully refer to \cite{kmm} and \cite{mo2}.
We will call an \textit{elementary} contraction
the contraction morphism of an extremal \textit{ray}.
For an extremal ray $R$ of $\overline{{\rm NE}}(M)$
of a projective manifold $M$,
we denote by $l(R)$ the length of the ray $R$,
i.e., 
\[l(R)=\min\{-K_M.C\,\,|\,\, C\subseteq M\textrm{ is a rational curve, 
and }[C]\in R\},\]
and by $E(R)$ the locus of the ray $R$,
i.e., the union of all curves belonging to $R$.

\section{Preliminaries and proofs of Propositions}\label{shinprelim}
In this section, we first recall some fundamental results
of Mori-Kawamata theory in our setting.
We refer to \cite{mo1}, \cite{mo2}, and \cite{kmm}
for details and proofs.

Let $M$ be a projective manifold of dimension $n$
and $\mathcal{E}$ an ample vector bundle of rank $r$ on $M$.
Set $H=\det \mathcal{E}$;
it is an ample line bundle.
For a polarized manifold $(M,H)$,
it is very important to study the adjunction bundle
$K_M+tH$ for a various $t\in \mathbb{Q}$,
as is known by the fundamental notion 
such as adjunction formula, sectional genus, etc.
If $K_M$ is not nef, we can apply Mori-Kawamata theory:
if $K_M+tH$ is nef for some positive rational 
number $t\in \mathbb{Q}_{>0}$,
the Basepoint-free Theorem \cite{kmm} implies that 
$\Bs \lvert m(K_M+t H)\rvert=\emptyset$
for all $m\gg 0$ such that $mt\in \mathbb{N}$.

Here we define the nef value $\tau (M, H)$ 
of a polarized manifold $(M,H)$
to be the minimum of the set of real numbers $t$
such that $K_M+tH$ is nef.
Set $\tau =\tau(M,H)$.
Then Kleiman's criterion for ampleness~\cite{kleiman}
implies that $K_M+tH$ is ample if $t>\tau$,
nef but not ample if $t=\tau$,
and not nef if $t<\tau$.

Suppose that $K_M$ is not nef, i.e., $\tau>0$, as above.
Then $\tau$ is a rational number by 
the Rationality Theorem \cite{kmm}.
Hence we can apply the Basepoint-free Theorem \cite{kmm} to 
$K_M+\tau H$ to see that 
$\Bs \lvert m(K_M+\tau H)\rvert=\emptyset$
for all $m\gg 0$ such that $m\tau\in \mathbb{N}$.
Let $\Psi_{\lvert m(K_M+\tau H)\rvert}:M\to
\mathbb{P}$ be the morphism defined by 
the linear system $\lvert m(K_M+\tau H)\rvert$,
and consider the Stein factorization of 
$\Psi_{\lvert m(K_M+\tau H)\rvert}$:
$\Psi_{\lvert m(K_M+\tau H)\rvert}=v_m\circ \Psi$,
where $\Psi:M\to T$ is the morphism onto a normal variety $T$
with connected fibers, and $v_m:T\to \mathbb{P}$ is the 
finite morphism.
We see that $\Psi$ is independent of $m$ chosen.
We call $\Psi$ the nef value morphism of the polarized 
manifold $(M,H)$.
Since $K_M+\tau H$ is nef but not ample
by Kleiman's criterion for ampleness~\cite{kleiman},
$\Psi_{\lvert m(K_M+\tau H)\rvert}$ is not finite, i.e., 
$\Psi$ is not isomorphic
and contracts some curve on $M$.
Let $C$ be a curve on $M$.
Then $\Psi$ contracts $C$ if and 
only if $(K_M+\tau H).C=0$.
This is the numerical characterization of curves 
contracted by $\Psi$.

Let $Z_1(M)$ be the free abelian group generated by 
all the irreducible reduced curves on $M$,
and set $N_1(M)=\{Z_1(M)/\!\!\equiv\}\otimes \mathbb{R}$,
where $\equiv$ denotes the numerical equivalence.
Let $\overline{{\rm NE}}(M)$ be the closure
in the Euclidean topology of the cone
${\rm NE}(M)$ in $N_1(M)$
generated by all the effective curves on $M$.
Set 
\[F=\{D\in N_1(M)\,|\, (K_M+\tau H).D=0\}\cap 
\overline{{\rm NE}}(M).\]
Then the Contraction Theorem \cite{kmm} says that 
$\Psi$ is characterized by the following two properties:
(i) $\Psi$ is a morphism onto a normal projective variety
with connected fibers;
(ii) $\Psi$ contracts a curve $C$ 
if and only if the numerical equivalence class $[C]$
of $C$ lies in $F$.
Namely $\Psi$ is uniquely determined by $F$,
and we call $\Psi$ the contraction morphism of 
an extremal face $F$.
This characterization of $\Psi$
implies that if a nef divisor $A$ defines
a supporting function of $F$, i.e.,
if, for a curve $C$, $A.C=0$ if and only if $[C]\in F$,
then the morphism with connected fibers onto a normal variety
induced, by the Basepoint-free Theorem \cite{kmm},
from high multiple $mA$ of $A$
is nothing but $\Psi$
and independent of $m$,
and thus $A=\Psi^*B$ for some ample Cartier divisor $B$ on $T$.
Moreover we see that a divisor $A'$ is the pull back of some 
Cartier divisor $B'$ of $T$ if and only if 
$A'.C=0$ for any curve $C$ 
(whose numerical equivalence class is) 
belonging to $F$,
because $A'+mA$ defines a supporting function of $F$
for a sufficiently large $m$,
where $A$ is assumed, as above, to be nef and to define
a supporting function of $F$.

It follows from the Cone Theorem~\cite{mo2}
that $F$ is a polyhedral cone, i.e.,
$F$ is spanned by the minimal finite set of half lines
$R$, $R_a$,\dots, $R_d$, called extremal rays,
as $F=R+R_a+\dots +R_d$,
and that each extremal ray, e.g., $R$, is 
represented, as $R=\mathbb{R}_{\geq 0}[C]$, by a rational curve $C$
such that $-K_M.C\leq n+1$.
Such a rational curve $C$ is called an
extremal rational curve of $R$.
Note here that $F$ is thus a rational cone in 
$N_1(M)$ since each extremal ray is represented by
an effective curve.
We also recall here the notion of the length
$l(R)$ of the ray $R$: the length $l(R)$ of $R$
is defined to be 
the minimum of the set of intersection numbers $-K_M.C$ of $-K_M$
and rational curves $C$ belonging to $R$, i.e.,
\[l(R)=\min\{-K_M.C\,\,|\,\, C\subseteq M\textrm{ is a rational curve, 
and }[C]\in R\}.\]
If an extremal rational curve $C_0$ of $R$
attains the length $l(R)$ of $R$, i.e.,
$-K_M.C_0=l(R)$, 
we call $C_0$ a $\it{minimal}$ extremal rational curve.

Here we give a proof of Propositions~\ref{upper bound}.
\begin{proof}[Proof of Proposition~\ref{upper bound}]
We show that an upper bound of $\tau $ is $(n+1)/r$:
if $\tau$ is positive then $\tau\leq (n+1)/r$. 
Thus we may assume, as above,  
that $\tau$ is positive.
Then, as we have seen in the above,
there exists an extremal rational curve $C$ 
of an extremal ray $R$ such that 
$(K_M+\tau \det \mathcal{E}).C=0$.
Since $\det \mathcal{E}.C\geq r$ for a rational curve $C$
and an ample vector bundle $\mathcal{E}$ of rank $r$,
this implies that $\tau r\leq n+1$.
\end{proof}

Now choose one extremal ray, say, $R$, and fix $R$.
Since $F$ is a rational cone,
if we perturb $\tau H$ a little,
we can take a $\mathbb{Q}$-Cartier divisor $H'$
such that $K_M+H'$ is nef
and that $R$ is supported by $K_M+H'$ on $M$:
\[R=\{D\in N_1(M)\,|\, (K_M+H').D=0\}\cap \overline{{\rm NE}}(M).\]
Then, by the same argument as above, applying the 
Basepoint-free Theorem \cite{kmm} to $K_M+H'$
and taking the Stein factorization,
we obtain a morphism $\psi:M\to S$ with connected fibers
onto a normal variety $S$.
We see that a curve $C$ is contracted by $\psi$
if and only if $[C]\in R$,
and that $\psi$ is uniquely determined by 
this property with 
$\psi_*\mathcal{O}_M=\mathcal{O}_S$;
thus $\psi$ is called the contraction morphism of $R$.
Since $(K_M+\tau H).C=0$ for any curve $C$ belonging to $R$,
as we have seen in case $\Psi$ and $F$,
we infer that $K_M+\tau H$ is the pull back 
of some line bundle on $S$.
Hence the nef value morphism $\Psi$ is factored through $\psi$.
Note here that the relative Picard number $\rho (M/S)$ is one.
We also call $\psi$ an elementary contraction.

Our first strategy to the classification
is to study the structure of $\psi$:
if we know the structure of $\psi$ well,
then we also know the structure of $(M,\mathcal{E})$.
Here let us recall the locus $E(R)$ of $R$:
it is the union of all curves 
belonging to $R$,
i.e., in terms of $\psi$, it is the locus contracted by $\psi$.
The first pivotal result concerning the property of $\psi$
is the following theorem \cite[Theorem (1.1)]{wl} of
Wi\'{s}niewski,
theorem which comes to this form from 
the equality $l(R)\leq n+1$ by 
Mori \cite{mo2}
through an improvement \cite[Theorem (0.4)]{i} by Ionescu:
\begin{thm}\label{wisniewski}
Let $M$ be a projective manifold
and let $R$ be an extremal ray of $\overline{{\rm NE}}(M)$.
Denote by $\psi$ the contraction morphism of $R$,
and by $E(R)$ the locus of $R$.
For each irreducible component $F(\psi)$
of any positive dimensional fiber of $\psi$
and for a given general point $x$ of $F(\psi)$,
set 
\[l_x(F(\psi))=\min \{-K_M.C\,\, |\,\,
 x\in C\subseteq F(\psi)\textrm{ is a rational curve}\}.\]
Then we have
\begin{equation}\label{upper bound of l(R)}
l(R)
\leq l_x(F(\psi))
\leq \dim F(\psi)+1-\codim (E(R),M).
\end{equation}
\end{thm}
This theorem gives an upper bound 
(\ref{upper bound of l(R)})
of the length $l(R)$.

On the other hand, we have an lower bound of $l(R)$ as follows.
Let $C_0$ be a minimal extremal rational curve of $R$.
Then we have 
\begin{equation}\label{lower bound of l(R)}
l(R)=-K_M.C_0=\tau \det\mathcal{E}.C_0\geq \tau r,
\end{equation}
since $(K_M+\tau \det\mathcal{E}).C_0=0$.

Let us apply these inequalities to deduce Proposition~\ref{n+1/r};
suppose that $n+1\leq \tau r$.
Then the lower bound~(\ref{lower bound of l(R)})
implies that $n+1\leq l(R)$,
and thus the upper bound~(\ref{upper bound of l(R)})
implies that $l(R)=n+1$, that $\dim F(\psi)=n$, 
and that $\psi$ is of fiber type.
Moreover equality holds in 
the lower bound~(\ref{lower bound of l(R)});
we have $\det\mathcal{E}.C_0=r$.
Furthermore we have $\dim S=0$,
since $\dim F(\psi)=n$,
and thus we infer that $K_M+\tau \det\mathcal{E}=0$.
In particular, we see that $M$ is a Fano manifold.
(In \cite{o}, I used \cite[Theorem (0.4)]{i} 
to deduce that $M$ is a Fano manifold
in the proof of  this  proposition.
However this is an error. Please apply Theorem~\ref{wisniewski}
instead of \cite[Theorem (0.4)]{i}.)

Now we come to apply the following 
very strong characterization
of projective space.
For a more unifying characterization,
which deduce the following as one of its corollaries,
we refer to the paper of Cho, Miyaoka and Shepherd-Barron
\cite[Theorems 0.1 and 0.2]{cms}.
We also refer to Kebekus's papers \cite{ke1} and \cite{ke2}
for a cornerstone of the proof of the following.

\begin{thm}\label{cm}
Let $M$ be a Fano manifold of dimension $n$ over the complex numbers.
If $-K_M.C\geq n+1$ for every rational curve $C\subset M$,
then $M$ is isomorphic to $\mathbb{P}^n$.
\end{thm}

It follows from this theorem that $M\cong \mathbb{P}^n$.
Hence we see that $C_0$ is a line in $\mathbb{P}^n$.
Since $\det\mathcal{E}.C_0=r$, we have 
$\det\mathcal{E}\cong \mathcal{O}(r)$.
Since $\mathcal{E}$ is ample, this implies that
$\mathcal{E}$ is a uniform vector bundle
of type $(1,\dots,1)$.

Here we recall some of the results 
on uniform vector bundles on $\mathbb{P}^n$
needed later, results due to 
Van de Ven~\cite{van-de-ven}, Sato~\cite{sato},
Elencwajg, Hirschowitz, and Schneider~\cite{e-h-s},
Elencwajg~\cite{ele-uniform}, \cite{ele-uniform-high},
Ellia~\cite{ellia},
and Ballico~\cite{ballico}.
\begin{thm}\label{uniform}
Let $\mathcal{E}$ be a uniform vector bundle
of rank $r$ on $\mathbb{P}^n$ (over an algebraically closed field
of characteristic zero).
\begin{enumerate}
\item If the type of $\mathcal{E}$ is $(0.\dots,0)$,
      then 
      $\mathcal{E}\cong \mathcal{O}^{\oplus r}$.
\item If the type of $\mathcal{E}$ is $(0,\dots,0,1)$,
      then $\mathcal{E}$ is either $\mathcal{O}^{\oplus (r-1)}\oplus 
      \mathcal{O}(1)$ or $\mathcal{O}^{\oplus (r-n)}\oplus 
      T_{\mathbb{P}^n}(-1)$.
\item If $r\leq n+1$, then $\mathcal{E}$ is one of the following:
      $\oplus_{i=1}^{r}\mathcal{O}(a_i)$,
      $T_{\mathbb{P}^n}(a)\oplus \mathcal{O}(b)^{\oplus i}$ ($i=0,1$),
      $\Omega_{\mathbb{P}^n}(a)\oplus 
      \mathcal{O}(b)^{\oplus i}$ ($i=0,1$),
      or $S^2T_{\mathbb{P}^2}$,
      where $a_i$, $a$, and $b$ are integers.
\end{enumerate}
\end{thm}

By Theorem~\ref{uniform} (1)
(see also \cite[Theorem 3.2.1]{oss} for a proof),
we have $\mathcal{E}\cong \mathcal{O}(1)^{\oplus r}$.
This completes the proof of Proposition~\ref{n+1/r}.

Note here that if we assume that $\tau\geq 1$
then we can deduce Proposition~\ref{n+1/r}
by applying the following (1) of the theorem 
of Kobayashi and Ochiai
\cite{ko}  instead of Theorem~\ref{cm}
(see \cite{o}).
\begin{thm}\label{ko}
Let $(M,L)$ be a polarized manifold of dimension $n$
over the complex numbers.
\begin{enumerate}
\item If $K_M+(n+1)L=0$, then $(M,L)$ is isomorphic to 
      $(\mathbb{P}^n,\mathcal{O}(1))$.
\item If $K_M+nL=0$, then $(M,L)$ is isomorphic to 
      $(\mathbb{Q}^n,\mathcal{O}(1))$,
      where $\mathbb{Q}^n$ denotes
      a hyperquadric in $\mathbb{P}^{n+1}$.
\end{enumerate}
\end{thm}

Finally we give a proof of Proposition~\ref{n+1/r+1},
which is,
as Proposition~\ref{n+1/r},
an easy consequence of Theorem~\ref{cm}.
\begin{proof}[Proof of Proposition~\ref{n+1/r+1}]
Since $n+1-\tau r<\tau \leq 1$,
we have $n\leq n+1-\tau<\tau r$.
The lower bound~(\ref{lower bound of l(R)})
therefore implies that $n+1\leq l(R)$.
As in the proof of Proposition~\ref{n+1/r},
the upper bound~(\ref{upper bound of l(R)})
then implies that $l(R)=n+1$ and that $M$
is a Fano manifold of Picard number one.
Hence it follows from Theorem~\ref{cm}
that $M\cong \mathbb{P}^n$.
By $n+1-\tau r<\tau $,
we also have $n+1<\tau (r+1)$.
Therefore we have $\tau \det\mathcal{E}.C_0=-K_M.C_0<\tau (r+1)$.
Hence we have $\det\mathcal{E}.C_0<r+1$,
since $\tau$ is positive.
This implies $\det\mathcal{E}.C_0=r$, 
since $\mathcal{E}$ is ample.
Thus we have $\mathcal{E}\cong \mathcal{O}(1)^{\oplus r}$,
as in the proof of Proposition~\ref{n+1/r}.
\end{proof}

\section{Relative version of Proposition~\ref{n+1/r}}
In this section, we will give 
a relative version 
Proposition~\ref{relative}
of Proposition~\ref{n+1/r}.
Proposition~\ref{relative} will be used
in the proof of Theorem~\ref{Main theorem}
in case $\psi:M\to S$ in \S~\ref{shinprelim}
is of fiber type.
Some results needed in the proof of Proposition~\ref{relative}
are collected in \S~\ref{prefibertype1}.

As in \S~\ref{shinprelim},
let $M$ be an $n$-dimensional projective
manifold, and $\mathcal{E}$
an ample vector bundle of rank $r$ on $M$.
Suppose that $K_M$ is not nef,
and denote by $\tau$ the nef value of $(M,\det\mathcal{E})$.
Let $R$ 
be an extremal ray of $\overline{{\rm NE}}(M)$
such that $(K_M+\tau\det\mathcal{E}).R=0$,
and $\psi:M\to S$ the contraction morphism of $R$.
Note that $K_M+\tau\det\mathcal{E}$ 
is the pull back of some $\mathbb{Q}$-Cartier divisor on $S$.

By inequalities (\ref{upper bound of l(R)})
and (\ref{lower bound of l(R)}) in \S~\ref{shinprelim},
we see first that 
\[\tau r\leq \dim F(\psi)+1-\codim (E(R),M)\]
for any irreducible component $F(\psi)$
of any positive dimensional fiber of $\psi$.
Note here that,
by abuse of notation,
we may assume that $F(\psi)$
takes the smallest dimension
among all the positive dimensional fibers of $\psi$.
Hence if $\dim F(\psi)+1\leq \tau r$
for some irreducible component $F(\psi)$
of some positive dimensional fiber of $\psi$
we infer that $\codim (E(R),M)=0$,
i.e., that $\psi$ is of fiber type,
and that  
\[\tau r= \dim F+1\]
for a general fiber $F$ of $\psi$.
On the other hand, 
if $\psi$ is of fiber type,
we have $K_F+\tau \det\mathcal{E}|_F=0$
for a general fiber of $\psi$,
and thus $\tau$ is also the nef value of $(F,\det\mathcal{E}|_F)$.
This is the reason why we call the following proposition
a relative version of Proposition~\ref{n+1/r}.

\begin{prop}\label{relative}
Suppose that 
$\dim F(\psi)+1\leq \tau r$
for some irreducible component $F(\psi)$
of some positive dimensional fiber of $\psi$.
Then 
$\psi$ is of fiber type,
$\dim S= n-\tau r+1$,
and 
$S$ has only rational Gorenstein singularities.
Let $U$ denote the largest open subset of $S$
such that $\psi^{-1}(U)\to U$ is smooth.
Set $d=\tau r-1$.
Then for any fiber $F$ of $\psi$
over any closed point of $U$,
we have
$(F,\mathcal{E}|_F)\cong
(\mathbb{P}^{d},\mathcal{O}(1)^{\oplus r})$
and $\codim (S\setminus U, S)\geq 3$.
Hence we have $\codim (\Sing S, S)\geq 3$.
Moreover
\begin{enumerate}
\item if $\tau r=n$, then 
      $\psi$ is a $\mathbb{P}^{n-1}$-bundle in the Zariski topology;
\item if $\tau r=n-1$, then
      $\psi$ is a $\mathbb{P}^{n-2}$-bundle in the \'{e}tale 
      (or complex)
      topology;
\item if $\tau r=n-2$, then
      $\psi$ has at most finite number of singular fibers.
\end{enumerate}
\end{prop}
\begin{proof}
As is noted above,
the assumption $\dim F(\psi)+1\leq \tau r$
implies that $\psi$ is of fiber type
and that $\tau r=\dim F(\psi)+1\geq 2$.
It follows immediately from Proposition~\ref{smooth}
that $S$ has only rational Gorenstein singularities.
Since $\tau r=\dim F+1$
for any closed fiber $F$ of $\psi^{-1}(U)\to U$
and $\tau $ is also the nef value of $(F,\det\mathcal{E}|_F)$,
as we have seen above,
Proposition~\ref{n+1/r} implies that 
$(F,\mathcal{E}|_F)\cong
(\mathbb{P}^{d},\mathcal{O}(1)^{\oplus r})$.

Let $S_2$ be the intersection
$D_1\cap\cdots \cap D_{n-d-2}$
of general very ample divisors
$D_1, \ldots, D_{n-d-2}$ of $S$.
Then the restricted morphism
$\psi^{-1}(S_2)\ra S_2$
is the contraction morphism of an extremal ray
from a smooth variety $\psi^{-1}(S_2)$.
Indeed, if this were not elementary,
then, by applying 
Proposition~\ref{Proposition 1.3}
to a $(d+2)$-dimensional
manifold $\psi^{-1}(S_2)$, 
we would have $d+1< (d+3)/2$,
i.e., $d<1$. This is a contradiction.
Hence $\psi^{-1}(S_2)\ra S_2$ is elementary.
Now it follows from Proposition~\ref{smooth}
that $S_2$ is smooth.

Take a general very ample divisor $S_1$ on $S_2$,
and let us consider the morphism $\psi^{-1}(S_1)\ra S_1$.
This morphism is again 
an elementary contraction;
indeed, if this were not elementary, again
by applying Proposition~\ref{Proposition 1.3}
to a $(d+1)$-dimensional
manifold $\psi^{-1}(S_1)$, 
we would have $d+1=(d+2)/2$, i.e., $d=0$.
This is a contradiction.
Thus $\psi^{-1}(S_1)\ra S_1$ is elementary.

Let $U_1$ denote the largest open subset of $S_1$ such that
$\psi^{-1}(U_1)\ra U_1$ is smooth.
Let $F_1$ be any closed fiber of the morphism $\psi^{-1}(U_1)\ra U_1$.
Note here 
$(F_1, \mathcal{E}|_{F_1})\cong (\mathbb{P}^{d}, \mathcal{O}(1)^{\oplus r})$.
Therefore, by Theorem~\ref{Brauer-Severi}, $\psi^{-1}(U_1)\ra U_1$ is
a $\mathbb{P}^{d}$-bundle in the \'{e}tale topology,
and on the space $\psi^{-1}(V)=V\times_{U_1}\psi^{-1}(U_1)$
over any small \'{e}tale open set $V$ of $U_1$
exists a line bundle $H_V$ such that 
the restriction of $H_V$ to any fiber of $\psi^{-1}(V)\ra V$
is isomorphic to $\mathcal{O}(1)$.
Since $\dim U_1=1$, 
Tsen's theorem
implies that $H^2({U_1}_{\textrm{et}},\mathbb{G}_m)=0$, 
where ${U_1}_{\textrm{et}}$ denotes $U_1$ with \'{e}tale topology.
(See, e.g., \cite[III p.108]{milne}.)
Hence, by modifying the glueing if necessary,
we can glue these $H_V$ in the \'{e}tale topology.
Moreover it follows from \cite[III.4.9]{milne} that
$H^1(\psi^{-1}(U_1)_{\textrm{et}}, \mathbb{G}_m)
=H^1(\psi^{-1}(U_1)_{\textrm{Zar}}, \mathcal{O}^{\times})
=\Pic \psi^{-1}(U_1)$.
Hence there exists an algebraic line bundle $H$ on $\psi^{-1}(U_1)$
such that $H|_{F_1}=\mathcal{O}_{F_1}(1)$ for any closed fiber $F_1$.
Since $\psi^{-1}(S_1)$ is smooth and $H$ is algebraic,
we can extend $H$ to a line bundle on $\psi^{-1}(S_1)$,
which we also denote by $H$ by abuse of notation.

Let $F_2$ be an arbitrary closed fiber of $\psi^{-1}(S_1)\ra S_1$.
Then $F_2$ is irreducible and reduced;
suppose, to the contrary, that $F_2=F'+F''$.
Note here that $F'$ is a Cartier divisor 
since $\dim S_1=1$ and $\psi^{-1}(S_1)$ is smooth.
Now that the Cartier divisor $F'$ satisfies the condition $F'.f=0$
for any curve $f$ in a fiber $F_1$ disjoint from $F'$,
we infer, 
by the property of a contraction morphism
of an extremal ray,
that $F'$ must be the pull back of 
some Cartier divisor on $S_1$.
Since $\dim S_1=1$,
this implies that $F'$ itself must be a fiber of 
$\psi^{-1}(S_1)\ra S_1$, 
a contradiction.

Note that the polarized variety $(F_1,H|_{F_1})$ has Fujita's delta genus
$\Delta (F_1,H|_{F_1})=0$ and  degree $(H|_{F_1})^{d}=1$.
Hence $(F_2,H|_{F_2})$ also has the same delta genus and degree,
so that $(F_2,H|_{F_2})\cong (\mathbb{P}^{d}, \mathcal{O}(1))$.
Therefore $\psi^{-1}(S_1)\ra S_1$ is a $\mathbb{P}^{d}$-bundle
in the Zariski topology;
in particular, if $\dim S=\dim S_1=1$,
i.e., $\tau r-1=d=n-1$,
we have (1) of the proposition.

Since $\psi^{-1}(S_1)\ra S_1$ is a
${\mathbb P}^{d}$-bundle,
$\psi^{-1}(S_2)\ra S_2$
has at most finite number of singular fibers.
Furthermore, since $\psi^{-1}(S_2)\to S_2$ is elementary,
it has no divisorial fibers;
since $\dim S_2=2$, this implies that 
$\psi^{-1}(S_2)\to S_2$ is equidimensional.
Therefore we infer that 
every closed fiber $F'$ of $\psi^{-1}(S_2)\ra S_2$ is isomorphic to 
${\mathbb P}^{d}$ 
and that $\mathcal{E}|_{F'}\cong \mathcal{O}(1)^{\oplus r}$
by the same argument as in 
\cite[\S 2.2, e-mail note of T. Fujita]{abw}.
Theorem~\ref{Brauer-Severi} thus implies that
$\psi^{-1}(S_2)\ra S_2$ is a $\mathbb{P}^d$-bundle
in the \'{e}tale topology.
(It is easy to see that 
$\psi^{-1}(S_2)\ra S_2$ is a $\mathbb{P}^d$-bundle
in the complex topology.)
In particular, if $\dim S=\dim S_2=2$,
i.e., $\tau r-1=d=n-2$,
we have (2) of the proposition.

Since $\psi^{-1}(S_2)\ra S_2$ is a $\mathbb{P}^d$-bundle
in the \'{e}tale (or complex) topology,
we have $S_2\subset U$,
and thus $\codim (S\setminus U, S)\geq 3$.
In particular, if $\dim S=3$,
i.e., $\tau r-1=d=n-3$,
we have (3) of the proposition.

Finally, since $\codim (S\setminus U, S)\geq 3$,
we have $\codim (\Sing S, S)\geq 3$.
\end{proof}

\section{Preliminaries for the fiber type case}\label{prefibertype1}
The following proposition is a slight improvement of 
\cite[Proposition~1.4]{abw}.
\begin{prop}\label{Ssmooth}
Let $\psi :M\ra S$ be a proper morphism of normal varieties 
with connected fibers, and suppose that $M$ is smooth
and that every divisor on $M$ is either dominating $S$
or the pull back of some Cartier divisor on $S$.
If $S$ has at worst quotient singularities then it is smooth.
\end{prop}
\begin{proof}
Assume, to the contrary, that $S$ has a singular point $s\in S$.
Since the question is local on $S$,
replacing $S$ with a small enough neighborhood around $s$
in the complex topology
if necessary,
we may assume that there exists a finite Galois cover
$p: T\ra S$ with Galois group $G$ from a smooth variety $T$.
Let $t\in T$ be a point such that $p(t)=s$.
We may moreover assume that $G$ is the stabilizer of $t$,
that $G$ is a subgroup of $GL(\dim S, \mathbb{C})$,
and, by \cite{ch}, that $G$ is small, i.e., 
$G$ contains no reflections.
Thus $p$ is \'{e}tale outside $\Sing (S)$.
Let $Z$ be the normalization of the fiber product of $M$ and $T$
over $S$. Then $G$ acts naturally on $Z$ over $M$.
Note that $Z\ra M$ is \'{e}tale in codimension one
since so is $p$
and the image of any divisor on $M$ via $\psi$
has, by assumption, codimension $\leq 1$ in $S$.
Let $n$ denote the dimension of $M$, and
$Y$ an $n$-dimensional irreducible component
of $Z$. Then 
$Y$ dominates $M$ and thus dominates $S$.
Hence it dominates $T$.
Since $\psi _*\mathcal{O}_M=\mathcal{O}_S$,
the irreducible component which dominates $T$ is unique.
Therefore $Y$ is the only $n$-dimensional component of $Z$.
This implies that the action of $G$ on $Z$ can be 
restricted to that on $Y$ over $M$.
Let $Y_t$ denote the fiber over $t$.
Since $G$ is the stabilizer of $t$, 
$G$ acts trivially on $Y_t$.
Hence $Y_t$ is contained in the ramification
locus of $Y\ra M$ because $\deg (Y\ra M)
=|G|\geq 2$.
On the other hand,
since $Y$ is normal, $M$ is smooth,
and $Y\ra M$ is \'{e}tale in codimension one,
it follows from the purity of branch loci that
$Y\ra M$ is \'{e}tale.
This is a contradiction.
\end{proof}

The following proposition is also a slight improvement
of \cite[Proposition~1.4.1]{abw};
our feature of the proof is an application 
of Koll\'{a}r's result~\cite{kollar1}.
\begin{prop}\label{smooth}
Let $\psi :M\ra S$ be the contraction morphism of 
an extremal ray,
i.e., an elementary contraction. 
Suppose that $M$ is smooth 
and that $\psi$ is of fiber type.
Then $S$ has only rational Gorenstein singularities.
Moreover if $\dim S=2$ then $S$ is smooth.
\end{prop}
\begin{proof}
First we see that every divisor on $M$
is either dominating $S$
or the pull back of some Cartier divisor on $S$,
since $\psi$ is an elementary contraction
of fiber type from a smooth variety.
By the some reason,
we also see that every integral Weil divisor on $S$
is a Cartier divisor;
in particular we see that $S$ is 1-Gorenstein.
Second note that $S$ has only rational singularities
by \cite[Cor.~7.4]{kollar1};
in particular $S$ is Cohen--Macaulay.
Now that $S$ is Cohen--Macaulay and 1-Gorenstein,
we infer that $S$ is Gorenstein.
Therefore $S$ has only rational 
Gorenstein singularities.
If $S$ is a surface, this means that
$S$ has only rational double points.
Now Proposition~\ref{Ssmooth} implies that $S$ is smooth.
\end{proof}

The following proposition due to \cite{abw1} 
(see also \cite[Proposition 1.3]{abw})
will be used in the proofs of
Propositions~\ref{relative}
and \ref{relative=n/r}.
For a proof, see also \cite[Proposition 2.9]{ao1}.
\begin{prop}\label{Proposition 1.3}
Let $\psi:M\to S$ be the contraction morphism
of an extremal face from a smooth projective variety $M$
of dimension $n$ onto a normal projective variety $S$
of dimension $<n$.
Suppose that every rational curve $C$
in a general fiber of $\psi$ satisfy a condition
$-K_M.C\geq (n+1)/2$.
Then $\psi$ is an elementary contraction except if
\begin{enumerate}
\item $-K_M.C=(n+2)/2$ for some rational curve $C$ on $M$,
      $S$ is a point, and $M$ is a Fano manifold of 
      pseudoindex $(n+2)/2$ and of Picard number $\rho (M)=2$.
\item $-K_M.C=(n+1)/2$ for some rational curve $C$ on $M$,
      and $\dim S\leq 1$.
\end{enumerate}
\end{prop}

Finally recall
the following theorem 
\cite[Th\'{e}or\`{e}me 8.2]{gro}
of Grothendieck in our setting.
\begin{thm}\label{Brauer-Severi}
Let $\psi:M\to S$ be a proper flat morphism of varieties
(over an algebraically closed field),
and suppose that a closed fiber $\psi^{-1}(s)$ of $\psi$
is isomorphic to $\mathbb{P}^{d}$.
Then there exist an open neighborhood $V'\subset S$ of $s$
and an \'{e}tale finite surjective morphism $V\to V'$
such that $M\times_{V'}V$ is isomorphic to $\mathbb{P}^{d}_V$
over $V$.
\end{thm}

\section{An overview of the proof of 
Theorem~\ref{Main theorem}}\label{overview}
We give an overview of the proof of Theorem~\ref{Main theorem}
in this section.
As in \S~\ref{shinprelim},
denote by $R$ an extremal ray of $\overline{{\rm NE}}(M)$
such that $(K_M+\tau\det\mathcal{E}).R=0$,
by $C_0$ a minimal extremal rational curve of $R$,
and by $\psi:M\to S$ the contraction morphism of $R$.

First recall the upper bound (\ref{upper bound of l(R)})
of $l(R)$ in Theorem~\ref{wisniewski} 
and the lower bound (\ref{lower bound of l(R)}) of $l(R)$:
we have 
\begin{equation}\label{bound of tau r}
\tau r
\leq \tau \det\mathcal{E}.C_0=l(R)\leq \dim F(\psi)+1-\codim (E(R),M)
\end{equation}
for any irreducible component $F(\psi)$ of 
any positive dimensional fiber of $\psi$.
By abuse of notation, we often regard $F(\psi)$
as a fiber which has the smallest dimension among all the positive
dimensional fibers of $\psi$.

Set $a=\det\mathcal{E}.C_0-r$. Since $\mathcal{E}$ is 
ample of rank $r$, $a$ is a non-negative integer.
We have 
\[\tau r\leq \tau (r+a)=\tau \det\mathcal{E}.C_0=l(R).\]
If $\tau \geq 1$,
we have, therefore,
\begin{equation}\label{bound of detE.C_0}
0\leq a\leq \tau a=l(R)-\tau r.
\end{equation}

Now we give an overview of the proof of Theorem~\ref{Main theorem}.

If $\codim (E(R),M)\geq 2$, i.e.,
$\psi$ is small,
we have $\dim F(\psi)\leq n-2$
since $F(\psi)\subseteq E(R)$.
Hence we have $\tau r\leq n-3$
by (\ref{bound of tau r}).
Since we assume that $n-3< \tau r$,
we infer that $\psi$ cannot be small.

If $\codim (E(R), M)=1$, 
the locus $E(R)$ of $R$
is a prime divisor by \cite[Proposition 5-1-6]{kmm},
and $\psi$ is called divisorial.
In case $\psi$ is divisorial,
since $F(\psi)\subseteq E(R)$,
we have 
\begin{equation}\label{divineq1}
\tau r\leq \tau \det\mathcal{E}.C_0
=l(R)\leq \dim F(\psi)\leq n-1.
\end{equation}
We divide the case into two cases: 1) $l(R)=\dim F(\psi)$;
2) $l(R)<\dim F(\psi)$.
For the case 1),
if we assume moreover that 
$l(R)=\dim F(\psi)$ for \textit{any} irreducible component of 
\textit{any} positive dimensional fiber of $\psi$
and that $\dim F(\psi)-2<\tau r$,
then we can give the classification;
it will be given in Theorem~\ref{minfiberdiv}
and Remark~\ref{relax}
in case $\dim F(\psi)-1<\tau r\leq \dim F(\psi)$,
and in Theorem~\ref{l(R)=dimF and a=1}
in case $\dim F(\psi)-2<\tau r\leq \dim F(\psi)-1$.
Note here that if $n-3<\tau r$ then
$n-2\leq \dim F(\psi)$ and thus $1\geq \dim \psi(E(R))$,
so that this additional assumption is always satisfied.
For the case 2), our assumption $n-3<\tau r$
implies that $l(R)=n-2$.
Furthermore we see that $l(R)-\tau r<1$,
and thus we have $\det\mathcal{E}.C_0=r$ by 
(\ref{bound of detE.C_0}).
This case will be treated in \S~\ref{shindiv}.
Finally note that the $\psi$ in the case (12)
of Theorem~\ref{Main theorem} is, in general,
not elementary and, in fact, the composite of 
divisorial elementary contractions of the case 1) above.
To obtain the $(M',\mathcal{E}')$ 
in the case (12) of Theorem~\ref{Main theorem},
we need to know furthermore the structure of $(M,\mathcal{E})$
with $\tau r=n-1$ and $\psi$ of fiber type.
So we will complete the classification of 
the case (12) of Theorem~\ref{Main theorem}
in \S~\ref{morebir}
after we have classified 
the case $\tau r=n-1$ and $\psi$ of fiber type.

If $\codim (E(R), M)=0$, i.e.,
$\psi$ is of fiber type,
let $F$ be a general fiber of $\psi$.
Since $(K_M+\tau \det\mathcal{E}).R=0$,
we have $K_F+\tau \det\mathcal{E}|_F=0$,
and this implies that $\tau $ is also the nef 
value of the general fiber $(F,\det\mathcal{E}|_F)$.
Set 
\[b=b(M,\mathcal{E})=\llcorner \dim M+1-\tau r\lrcorner.\]
We have $n-b<\tau r\leq n+1-b$, and we see, by (\ref{bound of tau r}),
that $b\geq 0$.
Moreover we have 
\[b(F,\mathcal{E}|_F)=b(M,\mathcal{E})-\dim S.\]
If $\dim S>0$, this implies that 
the classification of the general fiber $(F,\mathcal{E}|_F)$
is reduced to that of 
$(M,\mathcal{E})$ with smaller $b$ 
and satisfying the condition $K_M+\tau \det\mathcal{E}=0$.
(Note here that the Picard number $\rho (F)$ of $F$
is not necessarily one
and that $F$ might admit a birational contraction.)
We will proceed by induction on $b$
with the classification of 
the pairs $(M,\mathcal{E})$ with
$\dim S>0$.
Note here that, if $\tau \geq 1$, we have already classified 
the case $b=0$; Proposition~\ref{n+1/r} gives
the classification in case $\tau r=n+1$.
Since $n<\tau r$ in case $b=0$,
we have $l(R)-\tau r<1$,
and thus, if $\tau \geq 1$, inequality (\ref{bound of detE.C_0})
implies that $\det\mathcal{E}.C_0=r$. 
Therefore we have 
$\tau r=\tau \det\mathcal{E}.C_0=l(R)=n+1$.
So we may assume that 
we have already classified 
the case $b=0$.
We also see that, 
once we obtain the classification 
of some fixed $b$,
we have to make the relativization
of it to get the classification with bigger $b$
and with $\dim S>0$.
Note here that we have already made the relativization of the case $b=0$,
i.e., the case $n<\tau r\leq n+1$,
to the case $\dim F<\tau r\leq \dim F+1$
by Proposition~\ref{relative}
since $\tau\geq 1$;
indeed,
if $\dim F<\tau r$,
it follows from (\ref{bound of tau r}) above
that $l(R)=\dim F+1$,
and that $l(R)-\tau r<1$.
If $\tau \geq 1$, this implies that
$\det\mathcal{E}.C_0=r$ by (\ref{bound of detE.C_0}) above.
Therefore $\tau r=\tau \det\mathcal{E}.C_0=l(R)=\dim F(\psi)+1$.
Hence 
Proposition~\ref{relative} with $\tau\geq 1$ implies that
we have already classified the case
$\dim F<\tau r\leq \dim F+1$.

On the contrary to the case $\dim S>0$, 
we will classify the case $\dim S=0$ in a 
way as independent of the value of $b$ and $\tau r$ as possible.
The setup and the strategy to deal with the case $\dim S=0$
will be given in \S~\ref{the common setup}.
We will denote sections from \S~\ref{the common setup}
to \S~\ref{-K=n-2shin}
to the case $\dim S=0$.
I believe that the most results
obtained here,
in particular, those in \S~\ref{Slr1=n+1}, \S~\ref{Slr1=n}, 
and \S~\ref{Slr1=n-1}
are also applicable
to the case of $\tau r$ smaller,
e.g., 
to the case $n-3<\tau r<n-2$.
The classification of the case $\dim S=0$
is the main part of this paper.

Except for the case (12) of Theorem~\ref{Main theorem},
after we established the classification of the case $\dim S=0$ 
and of the case $\psi$ is birational,
we will proceed inductively on $b$
with the classification of the case $\dim S>0$
as follows.
Suppose that $b=1$.
Namely suppose that $n-1<\tau r\leq n$,
that $\psi$ is of fiber type,
and that $\dim S>0$.
Then we see by (\ref{bound of tau r})
that $\dim F=n-1$,
and thus $\dim F<\tau r\leq \dim F+1$ holds.
Therefore the case $b=1$ is already finished
by Proposition~\ref{relative} with $\tau\geq 1$.
We will make the relativization of the case 
$n-1<\tau r\leq n$
to the case 
$\dim F-1<\tau r\leq \dim F$ with $\dim S>0$ in 
\S~\ref{dim F-1<tau r leq dim F}.
Note that although the possible list of 
the general fiber $(F,\mathcal{E}|_F)$ with $\dim S>0$
is obtained by the inductive procedure,
it is another problem 
whether the possible fiber $(F,\mathcal{E}|_F)$
does really occur or not as a fiber of an elementary 
contraction of fiber type.
Moreover we have the problem
to determine the singular fiber of $\psi$.
Applying the classification of 
the case $\dim F-1<\tau r$,
we will obtain the classification 
of the case $b=2$,
i.e., the case $n-2<\tau r\leq n-1$
in \S~\ref{section=n-1/r}.
We will make the relativization
of the case $n-2<\tau r\leq n-1$
to the case $\dim F-2<\tau r\leq \dim F-1$
in \S~\ref{tau r=dimF-1}.
Finally applying the classification of 
the case $\dim F-2< \tau r$
we obtain the classification
of the case $\tau r=n-2$
in \S~\ref{section=n-2/r}.

\section{The case $\psi$ is divisorial and $l(R)=\dim F(\psi)$
}\label{prelimdiv}
The following theorem of Andreatta and Occhetta
\cite[Theorem 5.1]{ao}
plays a key role in case $\psi$ is divisorial
with $l(R)=\dim F(\psi)$
in the proof of Theorem~\ref{Main theorem}.
Note here that if $\psi$ is divisorial it follows from 
(\ref{divineq1}) in \S~\ref{overview}
that 
\begin{equation}\label{bound of tau r divisorial}
\tau r\leq \tau \det\mathcal{E}.C_0=l(R)\leq \dim F(\psi),
\end{equation}
where $C_0$ is a minimal extremal rational curve of $R$.
\begin{thm}\label{Andreatta Occhetta}
Let $M$ be an $n$-dimensional complex projective manifold,
$R$ an extremal ray of $\overline{{\rm NE}}(M)$,
and $\psi:M\to S$ the contraction morphism of $R$.
Suppose that $\psi$ is divisorial,
and denote by $E$ the exceptional divisor of $\psi$.
If $l(R)=\dim F(\psi)$ 
for any irreducible component $F(\psi)$
of any positive dimensional fiber of $\psi$,
then $S$ is smooth
and $\psi$ is the blowing up along a submanifold $\psi(E)$
of $S$.
\end{thm}

Based on Theorem~\ref{Andreatta Occhetta}, we give 
the following theorem, which is a
modification  of \cite[Theorem 3.1]{am} to our case.
The feature of our proof is the application
of Theorem~\ref{Andreatta Occhetta} and 
Ishimura's theorem~\cite{ishimura},
which simplify and clarify the argument
in \cite[Theorem 3.1]{am}.
\begin{thm}\label{minfiberdiv}
Let $M$ be an $n$-dimensional complex projective manifold,
and $\mathcal{E}$ an ample vector bundle of rank $r$ on $M$.
Let $\tau$ be the nef value of the polarized manifold
$(M,\det\mathcal{E})$.
Suppose that $K_M$ is not nef, i.e., $\tau>0$;
let $R$ be an extremal ray of $\overline{{\rm NE}}(M)$
such that $(K_M+\tau\det\mathcal{E}).R=0$,
and $\psi:M\to S$ the contraction morphism of $R$.
Suppose that $\psi$ is divisorial;
let $E$ be the exceptional divisor of $\psi$.
If $\dim F(\psi)\leq \tau r$ for any irreducible component
$F(\psi)$ of any positive dimensional fiber of $\psi$,
then $S$ is smooth, $\psi$ is the blowing up
along a smooth variety $\psi (E)$,
and $\mathcal{E}\cong \psi^*\mathcal{E}'\otimes \mathcal{O}(-E)$
for some vector bundle $\mathcal{E}'$ on $S$.

Let $\pi:P=\mathbb{P}(\mathcal{E})\to M$ be the projective space 
bundle associated to $\mathcal{E}$,
and $\pi':N=\mathbb{P}(\mathcal{E}')\to S$ the one associated to 
$\mathcal{E}'$.
Let $\varphi:P\to N$ 
be the morphism induced from the relation 
$\mathcal{E}\otimes \mathcal{O}_M(E)\cong \psi^*\mathcal{E}'$.
Then $\varphi$ is the contraction morphism
of some extremal ray $R_1$ of $\overline{{\rm NE}}(P)$.
Denote by $L$ the tautological line bundle $H(\mathcal{E})$
associated to $\mathcal{E}$,
and by $L'$ the one associated to $\mathcal{E}'$.
Then $\mathcal{E}'$ is ample
in the following cases;
\begin{enumerate}
\item $L\otimes \pi^*\mathcal{O}_M(E)$ is a good supporting
      divisor of $\varphi$,
      i.e., it is nef and it defines a supporting function of $R_1$;
\item $\mathcal{E}'|_{\psi(E)}$ is ample.
      In particular, if $\psi(E)$ is a point, then $\mathcal{E}'$ is ample.
\end{enumerate}
Finally, if $\mathcal{E}'$ is ample, we have 
$\tau (S,\det\mathcal{E}')\leq \tau $.
\end{thm}
\begin{proof}
Since $\psi$ is divisorial and $\dim F(\psi)\leq \tau r$,
all the inequalities 
in (\ref{bound of tau r divisorial})
above become equalities.
Hence we can apply Theorem~\ref{Andreatta Occhetta}
to see that $S$ is nonsingular
and 
that $\psi$ is the blowing-up along a submanifold $\psi(E)$ of $S$.
Now every positive dimensional fiber of $\psi$ 
is isomorphic to a projective space of dimension $l(R)$
and a minimal extremal rational curve $C_0$ 
of $R$
is a line in this projective space,
i.e., $-E.C_0=1$.
Since $\det\mathcal{E}.C_0=r$,
we have, by Theorem~\ref{uniform} (1),
$\mathcal{E}|_{F(\psi)}\cong 
\mathcal{O}_{{\mathbb P}^{l(R)}}(1)^{\oplus r}$
for every positive dimensional fiber $F(\psi)$ of $\psi$. 
Hence we have $\mathcal{E}\otimes \mathcal{O}_M(E)
\cong \psi^*\mathcal{E}'$
for some vector bundle $\mathcal{E}'$ of rank $r$ on $S$
by Ishimura's theorem~\cite{ishimura}.

Now note that the natural morphism 
$\varphi:P\to N$ is the blowing up of $N$
along ${\pi'}^{-1}(\psi (E))$.
Hence $-K_P$ is $\varphi$-ample,
and $\varphi$ is the contraction morphism
of some extremal ray $R_1$ of $\overline{{\rm NE}}(P)$.
Since $L\otimes \pi^*\mathcal{O}_M(E)\cong \varphi^*L'$,
we see that $L'$ is ample 
if $L\otimes \pi^*\mathcal{O}_M(E)$ is a good supporting divisor
of $\varphi$. 
If $\mathcal{E}'|_{\psi(E)}$ is ample,
then $L'|_{\varphi(E_1)}$ is ample,
where $E_1$ denotes the exceptional divisor of $\varphi$.
Since $L$ is ample, it follows from 
\cite[(5.7)]{fh}
that $L'$ is ample.
Therefore we infer that $\mathcal{E}'$ is ample
if (1) or (2) holds.
Finally, since we have $\tau r=l(R)$,
we see that $K_M+\tau\det \mathcal{E}=\psi^*
(K_S+\tau\det \mathcal{E}')$.
Hence $K_S+\tau \det \mathcal{E}'$ is nef.
Therefore we have
$\tau (S,\det \mathcal{E}')\leq \tau$
if $\mathcal{E}'$ is ample.
\end{proof}

\begin{rmk}\label{relax}
If we have $\tau\geq 1$ in Theorem~\ref{minfiberdiv},
we can relax 
the assumption $\dim F(\psi)\leq \tau r$ to
the one $\dim F(\psi)-1< \tau r$
for any irreducible component
$F(\psi)$ of any positive dimensional fiber of $\psi$;
indeed, if we have $\dim F(\psi)-1< \tau r$
then we have $l(R)-\tau r<1$ by 
(\ref{bound of tau r divisorial}) above,
and if we have $\tau \geq 1$
then inequality~(\ref{bound of detE.C_0}) 
in \S~\ref{overview}
implies that $\tau r=\det\mathcal{E}.C_0=l(R)$.
Hence, again by (\ref{bound of tau r divisorial}),
we have $\tau r=l(R)=\dim F(\psi)$ since $\dim F(\psi)-1< \tau r$.
\end{rmk}

\begin{rmk}
If there is no assumption such as (1) or (2) 
in Theorem~\ref{minfiberdiv},
the vector bundle $\mathcal{E}'$ in Theorem~\ref{minfiberdiv}
is not necessarily ample as the following example shows.
This is the reason why we attach the assumption such as (1) or (2) 
in Theorem~\ref{minfiberdiv}.
These kinds of assumptions are overlooked in \cite[Theorem 3.1]{am}.
The idea of showing the ampleness of the tautological line bundle
of $\mathcal{E}'$ by applying \cite[Lemma (5.7)]{fh} of Fujita
has its origin in the proof of \cite[Lemma (5.1)]{lm}
of Lanteri and Maeda.
\end{rmk}
\begin{ex}
Set 
$\mathcal{F}=
\mathcal{O}_{\mathbb{P}^p}^{\oplus q}\oplus 
\mathcal{O}_{\mathbb{P}^p}(-1)$,
where $p$ and $q$ are positive integers with $p+q=n$.
Let $f:S=\mathbb{P}(\mathcal{F})\to \mathbb{P}^p$
be the projection,
and $C$ the section corresponding to the quotient
$\mathcal{F}\to \mathcal{O}_{\mathbb{P}^p}(-1)$. 
Set $B=(H(\mathcal{F})\otimes f^*\mathcal{O}(1))^{\otimes 2}$.
Then $B$ is spanned but not ample since 
$B\otimes \mathcal{O}_C\cong \mathcal{O}_C$.
Let $\psi:M\to S$ be the blowing up along $C$,
and $E$ the exceptional divisor of $\psi$.
Set $A=\psi^*B\otimes \mathcal{O}(-E)$.
Note here that $A$ has the following expression:
\[A=\psi^*(H(\mathcal{F})\otimes f^*\mathcal{O}(2))
\otimes \psi^*H(\mathcal{F})\otimes \mathcal{O}(-E).
\]
First we see that 
$H(\mathcal{F})\otimes f^*\mathcal{O}(2)$
is ample and spanned.
Second we claim that 
$\psi^*H(\mathcal{F})\otimes \mathcal{O}(-E)$
is spanned.
The is because,
since the image of the natural map
$\mathcal{O}_S\otimes H^0(H(\mathcal{F}))\to H(\mathcal{F})$
is $\mathcal{I}_C\otimes H(\mathcal{F})$
where $\mathcal{I}_C$ is the ideal sheaf of $C$,
we have a surjection
$\mathcal{O}_M\otimes H^0(H(\mathcal{F}))\to
\psi^*H(\mathcal{F})\otimes \mathcal{O}(-E)$.
Hence $A$ is spanned.
Moreover $A$ is strictly nef;
if $C$ is a curve contracted to a point by $\psi$
then $A.C=-E.C>0$,
and if $C$ is a curve not contracted to a point by $\psi$
then we see that $A.C\geq 
(H(\mathcal{F})\otimes f^*\mathcal{O}(2)).\psi(C)>0$
by the above expression of $A$.
Therefore $A$ is ample.
Let $\mathcal{E}'=B^{\oplus r}$
and $\mathcal{E}=\psi^*\mathcal{E}'\otimes \mathcal{O}(-E)$.
Then $\mathcal{E}=A^{\oplus r}$
and thus $\mathcal{E}$ is ample.  However $\mathcal{E}'$ is not ample.
\end{ex}

\begin{thm}\label{l(R)=dimF and a=1}
Let $M$ be an $n$-dimensional complex projective manifold,
and $\mathcal{E}$ an ample vector bundle of rank $r$ on $M$.
Let $\tau$ be the nef value of the polarized manifold
$(M,\det\mathcal{E})$.
Suppose that $\tau\geq 1$.
Let $R$ be an extremal ray of $\overline{{\rm NE}}(M)$
such that $(K_M+\tau\det\mathcal{E}).R=0$,
and $\psi:M\to S$ the contraction morphism of $R$.
Suppose that $\psi$ is divisorial;
let $E$ be the exceptional divisor of $\psi$.
If $l(R)=\dim F(\psi)$ for any irreducible component
$F(\psi)$ of any positive dimensional fiber of $\psi$,
and $\dim F(\psi)-2<\tau r\leq \dim F(\psi)-1$,
then $S$ is smooth, $\psi$ is the blowing up
along a smooth variety $\psi (E)$,
$1\leq \tau<2$, $(l(R)-2)/2<r\leq l(R)-1$, 
and $\mathcal{E}$ fits in the following exact sequence
\[0\to
\psi^*\mathcal{E}'\otimes \mathcal{O}(-2E)
\to \mathcal{E}
\to {\psi'}^*\mathcal{F}\otimes \mathcal{O}_E(-E)
\to 0,\]
where $\mathcal{E}'$ is a vector bundle of rank $r$ on $S$,
$\psi'$ the restriction of $\psi$ to $E$,
and $\mathcal{F}$ a vector bundle of rank $r-1$
on $\psi(E)$.
\end{thm}
\begin{proof}
First we see 
by Theorem~\ref{Andreatta Occhetta}
that $S$ is smooth and 
that $\psi\colon M\to S$ 
is the blowing up along a submanifold
$\psi(E)$ of $S$.
Let $F$ denote a positive dimensional fiber of $\psi$;
$F$ is isomorphic to an $l(R)$-dimensional projective space.
Let $C_0$ be a minimal extremal rational curve of $R$.
We see that $C_0$ is a line in $F$,
i.e., $C_0$ is a smooth rational curve in $F$ with
$-E.C_0=1$.

Since $\dim F-2<\tau r\leq \dim F-1$
and $l(R)=\dim F$,
we see that $\lfloor l(R)-\tau r\rfloor=1$.
Since $\tau \geq 1$, 
it follows from (\ref{bound of detE.C_0})
in \S~\ref{overview}
that $\det\mathcal{E}.C_0\leq r+1$.
On the other hand,
since $\tau r\leq \dim F-1<\dim F=
l(R)=\tau \det\mathcal{E}.C_0$,
we have $r<\det\mathcal{E}.C_0$.
Therefore we have $\det\mathcal{E}.C_0=r+1$.
Hence 
$\tau (r+1)=\tau \det\mathcal{E}.C_0=l(R)=\dim F$,
and thus $\tau r=\dim F-\tau$.
Substituting this equality to
$\dim F-2<\tau r\leq \dim F-1$,
we obtain $1\leq \tau <2$.
This inequality together with
$\dim F-2<\tau r\leq \dim F-1$ then
implies that $(\dim F-2)/2<r\leq \dim F-1$.
Since $l(R)=\dim F$,
we have $(l(R)-2)/2<r\leq l(R)-1$.

Since $C_0$ is a line in $F$ and $\mathcal{E}$ is ample,
$\det\mathcal{E}.C_0=r+1$
implies that $\mathcal{E}|_F$
is a uniform vector bundle
of type $(1,\dots,1,2)$.
Note here that since $\tau\geq 1$
we have $r\leq \tau r\leq \dim F-1$.
Hence we have 
$\mathcal{E}|_F\cong \mathcal{O}(1)^{\oplus (r-1)}\oplus \mathcal{O}(2)$
by Theorem~\ref{uniform} (2)
(or \cite[Theorem 3.2.3]{oss} for a proof).

Set $\mathcal{L}_{\psi(E)}=
{\psi'}_*(\mathcal{E}(2E)|_E)$
and $\mathcal{M}_{\psi(E)}=
{\psi'}_*(\check{\mathcal{E}}(-E)|_E)$.
Then $\mathcal{L}_{\psi(E)}$ is a line bundle on $\psi(E)$
and $\mathcal{M}_{\psi(E)}$ is a vector bundle of rank $r-1$ on $\psi(E)$.
Moreover $\mathcal{E}(2E)|_E$ fits in the following exact sequence
\[
0\to {\psi'}^*\mathcal{L}_{\psi(E)}\to \mathcal{E}(2E)|_E
\to {\psi'}^*\check{\mathcal{M}}_{\psi(E)}\otimes \mathcal{O}_E(E)
\to 0.
\]
Set $\mathcal{F}=\check{\mathcal{M}}_{\psi(E)}$.
Let $\mathcal{G}$ be the kernel of the composite
of the natural maps 
$\mathcal{E}(2E)\to \mathcal{E}(2E)|_E$
and $\mathcal{E}(2E)|_E
\to {\psi'}^*\mathcal{F}\otimes \mathcal{O}_E(E)$:
$\mathcal{G}$ fits in the following exact sequence
\[0\to \mathcal{G}\to \mathcal{E}(2E) \to
{\psi'}^*\mathcal{F}\otimes \mathcal{O}_E(E)\to 0.\]
Since we have 
$\Tor^{\mathcal{O}_{M,x}}_i(k(x),\mathcal{O}_{E,x})=0$ 
for all $i\geq 2$ and all point $x\in M$,
we see that $\mathcal{G}$ is a vector bundle of rank $r$.
As is seen in the diagram of \cite[Theorem 1.3 (ii)]{maru},
$\mathcal{G}$ fits in the following exact sequence
\[0\to \mathcal{E}(E)\to \mathcal{G}\to 
{\psi'}^*\mathcal{L}_{\psi(E)}\to 0.\]
Note that the natural map
$\mathcal{G}\to 
{\psi'}^*\mathcal{L}_{\psi(E)}$ factors as the composite of
$\mathcal{G}\to \mathcal{G}|_E$ and 
$\mathcal{G}|_E\to 
{\psi'}^*\mathcal{L}_{\psi(E)}$.
Moreover we see,
as in \cite[Theorem 1.3 (ii)]{maru},
that the kernel of the map $\mathcal{G}|_E\to 
{\psi'}^*\mathcal{L}_{\psi(E)}$
is equal to the cokernel ${\psi'}^*\mathcal{F}$
of the map
$\mathcal{G}(-E)\to \mathcal{E}(E)$.
Hence we obtain the following exact sequence
\[0\to {\psi'}^*\mathcal{F}
\to \mathcal{G}|_E\to 
{\psi'}^*\mathcal{L}_{\psi(E)}
\to 0.\]
Thus we have 
$\mathcal{G}|_F\cong \mathcal{O}^{\oplus r}$.
Therefore we have, by \cite{ishimura}, $\mathcal{G}=\psi^* \mathcal{E}'$ 
for some vector bundle $\mathcal{E}'$
of rank $r$ on $S$.
This completes the proof of the theorem.
\end{proof}

Here we give a proof of Theorem~\ref{Main theorem}
in case $\psi$ is divisorial and 
$l(R)=\dim F(\psi)$ for any irreducible component of 
any positive dimensional fiber of $\psi$.
Let $E$ be the exceptional divisor of $\psi$.
In this case, 
since $n-3<\tau r$,
as we observed in \S~\ref{overview},
we have the following two cases by (\ref{divineq1})
in \S~\ref{overview}:\\
1) $\dim F(\psi)-1<\tau r\leq \dim F(\psi)$;\\
2) $n-3=\dim F(\psi)-2<\tau r\leq \dim F(\psi)-1$.

In case 1), we apply Theorem~\ref{minfiberdiv}
together with Remark~\ref{relax}
to see that $\tau r=\dim F(\psi)$,
that $S$ is smooth,
that $\psi$ is the blowing up
along a submanifold $\psi(E)$ of $S$,
and that $\mathcal{E}\otimes \mathcal{O}_M(E)
\cong \psi^*\mathcal{E}'$
for some vector bundle $\mathcal{E}'$ of rank $r$
on $S$.
Since $n-3<\tau r$, we can divide this case 1) 
into the following two sub-cases:\\
a) $\dim F(\psi)=n-1$;\\
b) $\dim F(\psi)=n-2$.\\
If $\dim F(\psi)=n-1$, then 
we infer that $F(\psi)=E$, and thus 
$\psi$ contracts $E$ to a point.
Hence, again by Theorem~\ref{minfiberdiv}, 
we infer 
that $\mathcal{E}'$ is ample,
and that 
$\tau (S,\det\mathcal{E}')\leq \tau $.
Put $M_1=S$ and $\mathcal{E}_1=\mathcal{E}'$.
We have $\tau (M_1,\det \mathcal{E}_1)r\leq n-1$.
As we noticed in \S~\ref{overview},
in order to obtain the case (12) of Theorem~\ref{Main theorem},
we need to analyze the structure of $(M_1,\mathcal{E}_1)$
in case $\tau (M_1,\det \mathcal{E}_1)r=n-1$.
This is done in \S~\ref{morebir},
after we also classified
$(M,\mathcal{E})$
with $\tau r=n-1$ and $\psi$ of fiber type.
If $\dim F(\psi)=n-2$, 
put $M'=S$.
This is the case (27) of Theorem~\ref{Main theorem}.

In case 2), we apply Theorem~\ref{l(R)=dimF and a=1}
to obtain the case (25) of Theorem~\ref{Main theorem};
note that, as can be seen from
the proof of Theorem~\ref{l(R)=dimF and a=1},
$\tau r=n-2$ implies that $\tau =1$,
and thus we have $r=n-2$.

This completes the proof of Theorem~\ref{Main theorem}
in case $\psi$ is divisorial and 
$l(R)=\dim F(\psi)$ for any irreducible component of 
any positive dimensional fiber of $\psi$.

Finally we give an example of the case (25) of the theorem,
which is also the simplest example of Theorem~\ref{l(R)=dimF and a=1}.
\begin{ex}\label{counter example of Andreatta Mella}
Let $M'$ be an $n$-dimensional projective manifold 
and $\psi :M\ra M'$ be the blowing-up of $M'$ at a point $p$ of $M'$.
Denote by $E$ the exceptional divisor of $\psi$.
By tensoring the pull back of a sufficiently ample line bundle
$\mathcal{L}$ on $M'$, we can make 
$\psi^*\mathcal{L}\otimes \mathcal{O}(-E)$
and $\psi^*\mathcal{L}\otimes \mathcal{O}(-2E)$ be ample.
Let $\mathcal{E}$ be $\psi^*\mathcal{L}\otimes
(\mathcal{O}(-E)^{\oplus (n-3)}\oplus \mathcal{O}(-2E))$.
Then $(M, \mathcal{E})$ gives an example of the case (25) of 
Theorem~\ref{Main theorem}.
\end{ex}

\section{The case 
$\psi$ is divisorial and $l(R)<\dim F(\psi)$}\label{shindiv}
Let $R$, $C_0$, $\psi:M\to S$ be as in \S~\ref{overview}.
In this section we give a proof of Theorem~\ref{Main theorem}
in case $\psi$ is divisorial
and $n-3<\tau r\leq l(R)<\dim F(\psi)\leq n-1$
for some irreducible component $F(\psi)$
of some positive dimensional fiber of $\psi$.
In this case, we have $l(R)=n-2$ and $\dim F(\psi)=n-1$.
Hence $F(\psi)$ is the exceptional divisor $E$ of $\psi$.
Since $\lfloor l(R)-\tau r\rfloor=0$
and $\tau \geq 1$,
we have $r=\det\mathcal{E}.C_0$
by (\ref{bound of detE.C_0}) in \S~\ref{overview}
as we have seen in \S~\ref{overview}.

Here, since $l(R)=n-2$, 
we come to apply the following theorem 
of Andreatta and Occhetta \cite[Theorem 5.2]{ao}.

\begin{thm}\label{Andreatta Occhetta2}
Let $M$ be an $n$-dimensional complex projective manifold,
$R$ an extremal ray of $\overline{{\rm NE}}(M)$,
and $\psi:M\to S$ the contraction morphism of $R$.
Suppose that $\psi$ is divisorial,
and denote by $E$ the exceptional divisor of $\psi$.
If $l(R)=n-2$, then one of the following cases occur:
\begin{enumerate}
\item $S$ is smooth and $\psi$ is the blowing up
      along a smooth curve $\psi(E)$;
\item $\psi(E)$ is a point and $(E,\mathcal{O}_E(-E))\cong
      (\mathbb{P}^{n-1},\mathcal{O}(2))$;
\item $\psi(E)$ is a point and $(E,\mathcal{O}_E(-E))\cong
      (\mathbb{Q}^{n-1},\mathcal{O}(1))$,
      where $\mathbb{Q}^{n-1}$ is a possibly singular quadric.
\end{enumerate}
\end{thm}

Since $F(\psi)=E$,
we see that we are in the cases  (2) and  (3) of 
Theorem~\ref{Andreatta Occhetta2}.
Moreover we have $\mathcal{E}|_E\cong \mathcal{O}(1)^{\oplus r}$
by Theorem~\ref{uniform} (1) in the case (2)
and by \cite[Lemma 3.6.1]{w3} in the case (3),
since $\det\mathcal{E}.C_0=r$
for a line $C_0$ in $E$.
These are  the cases (a) and (b) of (26) of Theorem~\ref{Main theorem}.

\section{Setup and strategy for the case $\dim S=0$
}\label{the common setup}
As in \S~\ref{overview},
denote by $R$ an extremal ray of $\overline{{\rm NE}}(M)$
such that $(K_M+\tau\det\mathcal{E}).R=0$,
by $C_0$ a minimal extremal rational curve in $R$,
and by $\psi:M\to S$ the contraction morphism of $R$.
As we stated
in \S~\ref{overview},
we assume that $\dim S=0$ in this section,
and give the setup and the strategy to 
deal with the case $\dim S=0$.

Let $P$ be the projective space bundle $\mathbb{P}(\mathcal{E})$ over $M$,
$\pi :P\ra M$ the projection, and $L$ the tautological line
bundle $H(\mathcal{E})$.

Suppose that $-K_P$ is ample in this section.
Note here that if $\tau \geq 1$
then $-K_P$ is ample;
indeed, if $\tau \geq 1$,
we have $(K_M+\det \mathcal{E}).R\leq 0$.
Since $\dim S=0$, this implies that
$-(K_M+\det \mathcal{E})$ is nef, so that
$-\pi ^*(K_M+\det \mathcal{E})$ is nef.
Therefore $-K_P$ is ample
because $L$ is ample.

Since the Picard number $\rho (P)$ is two,
this implies that $\overline{{\rm NE}}(P)$
is spanned by two extremal rays $R_{\pi}$ and $R_1$,
where $R_{\pi}$ is the ray corresponding to $\pi :P\ra M$
and $R_1$ is the other ray.
Let $\varphi :P\ra N$ be the contraction morphism of $R_1$.

Let $F({\varphi})$ 
be any irreducible component
of any positive dimensional fiber of $\varphi$.
The following lemma due to 
Ye and Zhang \cite{yz} and Peternell \cite{p0}
is the first key observation
to the study of $\varphi$.

\begin{lemma}\label{finite of yz}
The induced morphism $\pi|_{F(\varphi)}$ is finite.
In particular, $\dim F(\varphi)\leq n$.
\end{lemma}
\begin{proof}
Since any curve in $F(\varphi)$ is contracted by $\varphi$,
it (or, strictly speaking,
its numerical equivalence class) belongs to $R_{1}$,
and does not belong to $R_{\pi}$.
Therefore it is not contracted by $\pi$.
Hence  
$\pi|_{F(\varphi)}$ is finite and 
$\dim F(\varphi)=\dim \pi(F(\varphi))\leq \dim M=n$.
\end{proof}
Combining Lemma~\ref{finite of yz} and Theorem~\ref{wisniewski} (1),
we have
\begin{equation}\label{ub of r1}
\begin{split}
l(R_1)\leq l_x(F(\varphi))&\leq\dim F({\varphi})+1
-\codim (E(R_1),P)\\
&
\leq n+1-\codim (E(R_1),P).
\end{split}
\end{equation}
This inequality gives an upper bound of $l(R_1)$.

Let $C_1\subset P$ be 
a minimal extremal rational curve of $R_1$.
Since $\psi (\pi (C_1))$ is a point, $\pi (C_1)$
belongs to $R$, and therefore 
$(K_M+\tau \det \mathcal{E}).\pi (C_1)=0$.
Hence we have 
$l(R_1)= -K_P.C_1=rL.C_1+(\tau -1)\det \mathcal{E}.\pi_*(C_1)$.

We will use the following terminology.
\begin{defi}\label{def of unsplit}
Let $P$ be a projective manifold,
and $C_1'$ a rational curve on $P$.
We will say, for simplicity, that $C_1'$ is unsplit
if every maximal family $F\to S$, 
$S\subseteq \Chow (P)$,
of rational curves on $P$
containing $C_1'$ as a closed fiber
(see \cite[p.19]{cms})
is unsplitting.
In other words, 
$C_1'$ is unsplit if and only if 
$C_1'$ cannot be effectively algebraically equivalent
(see \cite[p.121--122]{kollar} for the definition)
to a sum 
$\Sigma_{i=1}^{\delta} D_i$
of $\delta$ ($\delta\geq 2$) rational curves $D_i$,
some of which may equal.
\end{defi}

Now we note the following lemma.

\begin{lemma}\label{big r+1}
Let $C_1'$ be an unsplit rational curve on $P$
with $\pi(C_1')$ a curve.
Then 
every quotient line bundle of 
$\mathcal{E}\otimes \mathcal{O}_{\tilde{C}_1'}$
has degree at least $L.C_1'$,
where $\tilde{C}_1'\to C_1'$
is the normalization
and $\mathcal{E}\otimes \mathcal{O}_{\tilde{C}_1'}$ denotes
$(\pi^*\mathcal{E})|_{C_1'}\otimes 
\mathcal{O}_{\tilde{C}_1'}$
for simplicity.
In particular, we have $\det \mathcal{E}.\pi _*(C_1')\geq rL.C_1'$,
and equality holds if and only if 
$\mathcal{E}\otimes \mathcal{O}_{\tilde{C_1'}}
\cong \mathcal{O}(L.C_1')^{\oplus r}$.
\end{lemma}
\begin{proof}
Since $\pi(C_1')$ is a curve,
we see that 
$\tilde{C}_1'\to C_1'\to \pi(C_1')\hookrightarrow M$
is finite.
Therefore the induced morphism 
$P\times_M\tilde{C}_1'\to P$ is finite.
Note that $P\times_M\tilde{C}_1'\to \tilde{C}_1'$
has the section
corresponding to the quotient 
$\mathcal{E}\otimes \mathcal{O}_{\tilde{C}_1'}
\to L_{C_1'}\otimes \mathcal{O}_{\tilde{C}_1'}$.
Denote by $D'$ the section.
Since $C_1'$ is unsplit and 
$P\times_M\tilde{C}_1'\to P$ is finite,
we infer that $D'$
is also unsplit in $P\times_M\tilde{C}_1'$.
Now the next Lemma~\ref{minimal quotient}
implies that the quotient line bundle
$\mathcal{E}\otimes \mathcal{O}_{\tilde{C}_1'}
\to L_{C_1'}\otimes \mathcal{O}_{\tilde{C}_1'}$
has the minimal degree
among all the quotient line bundles of 
$\mathcal{E}\otimes \mathcal{O}_{\tilde{C}_1'}$.
Therefore we have 
$\det \mathcal{E}.\pi _*(C_1')
=\deg \mathcal{E}\otimes \mathcal{O}_{\tilde{C}_1'}
\geq rL.C_1'$.
\end{proof}

\begin{lemma}\label{minimal quotient}
Let $\mathcal{F}=\oplus_{i=1}^{r}\mathcal{O}(d_i)$
be a vector bundle of rank $r$ on $\mathbb{P}^1$,
where $d_1\leq d_2\leq \dots\leq d_r$ are integers.
Let $D'$ be a section
of the projection $\pi:\mathbb{P}(\mathcal{F})\to \mathbb{P}^1$
corresponding to a quotient $\mathcal{F}\to \mathcal{O}(d)$.
If $d>d_1$, then $D'$ is effectively rationally equivalent
to $D+(d-d_1)l$,
where $D$ is the section
corresponding to the minimal quotient 
$\mathcal{F}\to \mathcal{O}(d_1)$
and $l$ is a line in a fiber of the projection $\pi$.
(See \cite[Definition (4.1.4)]{kollar}
for the definition of effective rational equivalence.)
In particular, if $D'$ is unsplit,
then $d=d_1$.
\end{lemma}
\begin{proof}
Suppose that $d>d_1$.
Since $D'$ defines a section $s'$
of $H^0(\check{\mathcal{F}}(d))$
such that $(s')_0=\emptyset$,
the assumption $d>d_1$ implies that $d-d_2\geq 0$.
Take $s_1\in H^0(\mathcal{O}(d-d_1))$
and $s_2\in H^0(\mathcal{O}(d-d_2))$
such that $(s_1)_0\cap (s_2)_0=\emptyset$.
Let $s''$ be a section of $H^0(\check{\mathcal{F}}(d))$
such that $s''=(s_1,s_2,0,\dots,0)$.
Then $s''$ defines a section $D''$ of $\pi$.
Let $\check{V}$ be the linear subspace of 
$H^0(\check{\mathcal{F}}(d)\oplus \check{\mathcal{F}}(d))$
spanned by $(s',0)$ and $(0,s'')$,
and $V$ the dual of $\check{V}$.
The injection 
$\mathcal{O}_{\mathbb{P}^1}\otimes \check{V}
\to \check{\mathcal{F}}(d)\oplus \check{\mathcal{F}}(d)
=\check{\mathcal{F}}(d)\otimes \check{V}$
as vector bundles induces 
the surjection $\mathcal{F}\otimes V\to \mathcal{O}(d)\otimes V$.
Hence we have a family 
$\mathbb{P}^1\times \mathbb{P}(V)
\subset \mathbb{P}(\mathcal{F})\times \mathbb{P}(V)$
over a projective line $\mathbb{P}(V)$,
which has $D'$ as a fiber and 
$D''$ as another fiber.
Therefore $D''$ is rationally equivalent to $D'$.
Note here that we can regard both $D$ and $D''$ 
as a section on a ruled surface 
$\mathbb{P}(\mathcal{O}(d_1)\oplus \mathcal{O}(d_2))$,
and we see that 
$D\in \lvert H(\mathcal{O}(d_1-d_2)\oplus \mathcal{O})\rvert$
and that 
$D''\in \lvert H(\mathcal{O}(d-d_2)\oplus \mathcal{O}(d-d_1)\rvert$
since $D''$ corresponds to the quotient
$\mathcal{O}(d_1)\oplus \mathcal{O}(d_2)\to \mathcal{O}(d)$.
Hence  $D+(d-d_1)l$ and $D''$ are
linearly equivalent on the ruled surface,
where $l$ is the fiber of the projection
$\mathbb{P}(\mathcal{O}(d_1)\oplus \mathcal{O}(d_2))\to
\mathbb{P}^1$.
Therefore $D'$ is rationally equivalent to 
$D+(d-d_1)l$.
\end{proof}

\begin{cor}\label{C_1'to pi(C_1') birational}
Let $C_1'$ be an unsplit rational curve on $P$
with $\pi(C_1')$ a curve.
Then $C_1'\to \pi(C_1')$ is birational.
\end{cor}
\begin{proof}
Set $C_0'=\pi(C_1')$, and 
let $\tilde{C_1'}\to C_1'$
and $\tilde{C_0'} \to C_0'$
be the normalizations.
Denote by $d$ the degree of the map $C_1'\to C_0'$.
As we saw in the proof of Lemma~\ref{big r+1},
$\tilde{C_1'}\to P$ induces 
an unsplit section $D'$ of 
$P\times_M\tilde{C_1'}\to \tilde{C_1'}$,
and by Lemma~\ref{minimal quotient}
we see that $D'$ is a minimal section
of $P\times_M\tilde{C_1'}\to \tilde{C_1'}$,
i.e., a section corresponding to 
a quotient line bundle
of $\mathcal{E}\otimes \mathcal{O}_{\tilde{C_1'}}$
of minimal degree.
Note here that minimal sections of 
$P\times_M\tilde{C_1'}\to \tilde{C_1'}$
and minimal sections of  
$P\times_M\tilde{C_0'}\to \tilde{C_0'}$.
are in one-to-one correspondence
via the pull back by $\tilde{C_1'}\to \tilde{C_0'}$;
let $D_0'$ be the minimal section 
of $P\times_M\tilde{C_0'}\to \tilde{C_0'}$
corresponding to $D'$.
Here we see that 
the image of $D'$ by 
$P\times_M\tilde{C_1'}\to P$
is $C_1'$ by the definition of $D'$.
Hence the image of $D_0'$
in $P$ is $C_1'$.
Thus the image of $D'$
in $P$ as a cycle is $dC_1'$.
We see also that $D'\to C_1'$ is birational 
by the definition of $D'$.
Therefore we have $d=1$.
\end{proof}

Note here that a minimal extremal rational curve $C_1$ is unsplit;
Lemma~\ref{big r+1} then implies that
\begin{equation}\label{-K_P.C_1alpha siki}
\begin{split}
l(R_1)= -K_P.C_1&=    rL.C_1+(\tau -1)\det \mathcal{E}.\pi_*(C_1)\\
      &= rL.C_1+(\tau -1)(rL.C_1+\alpha)\\
      &= \tau rL.C_1+(\tau -1)\alpha,
\end{split}
\end{equation}
where we set $\det \mathcal{E}.\pi_*(C_1)
=rL.C_1+\alpha$ for some non-negative integer $\alpha$.
This gives an an lower bound of $l(R_1)$
in case $\tau \geq 1$ as follows:
\begin{equation}\label{lb of r1=}
\begin{split}
l(R_1)&=\tau rL.C_1+(\tau -1)\alpha\\
      &\geq \tau rL.C_1.
\end{split}
\end{equation}

Here we pose the following problem:
\begin{problem}\label{C_1 C_0 comparison}
Is $\pi(C_1)$ always a minimal extremal rational curve in $R$ ?
\end{problem}

If we answer affirmatively to
Problem~\ref{C_1 C_0 comparison},
we have the following relation between $l(R_1)$ and $l(R)$.
\begin{prop}\label{lr1-lr=a}
If $\pi(C_1)$ is a minimal extremal rational curve
in $R$, then we have 
$l(R_1)+\alpha=l(R)$,
where we set $\det\mathcal{E}.\pi(C_1)=rL.C_1+\alpha$
for some non-negative integer $\alpha$ as above.
In particular we have $l(R_1)\leq l(R)$.
\end{prop}
\begin{proof}
First we have $\pi_*(C_1)=\pi(C_1)$ by 
Corollary~\ref{C_1'to pi(C_1') birational}.
Second if $\pi(C_1)$ is a minimal extremal rational curve
in $R$, we have $-K_M.\pi(C_1)=l(R)$.
Hence we have $l(R_1)=-K_P.C_1
=rL.C_1-K_M.\pi(C_1)-\det\mathcal{E}.\pi(C_1)
=l(R_1)-\alpha$.
\end{proof}

The idea of comparison of $C_1$ of $R_1$ and $C_0$ of $R$
stems from
\cite[Comparison Lemma (3.1)]{psw}:
Peternell-Szurek-Wi\'{s}niewski
gave an affirmative answer to 
Problem~\ref{C_1 C_0 comparison}
in case $r=n-1\geq 4$ and $K_M+\det\mathcal{E}=0$
(i.e., $\tau =1$).
Very roughly speaking, 
their strategy to the affirmative answer
is as follows: 
since $n+1\geq -K_M.C_0=\det\mathcal{E}.C_0\geq r=n-1$,
if $n-1=r\geq 3$, $\mathcal{E}\otimes \mathcal{O}_{\tilde{C_0}}$
has $\mathcal{O}(1)$ as quotient,
where $\tilde{C_0}\to C_0$ is the normalization.
Corresponding to this quotient $\mathcal{O}(1)$,
there is a rational curve $C_1'\subset P$ 
dominating $C_0$
such that $L.C_1'=1$.
Here it is clear but important that
$C_1'$ is unsplit since $L.C_1'=1$.
Considering some deformation family of $C_1'$,
they show that if $n\geq 5$ there exists a rational curve
effectively algebraically equivalent to $C_1'$
and contracted by $\varphi$.
The property of the contraction morphism
then implies that $C_1'$ belongs to $R_1$.
Since $L.C_1'=1$, this implies that $C_1'$
is a minimal extremal rational curve of $R_1$,
i.e., $C_1'=C_1$. Hence $\pi(C_1)=C_0$.

In \S~\ref{section comparison},
we extend the above argument of \cite[Comparison Lemma (3.1)]{psw}
to our case.
Here we do not restrict ourselves to 
consider the deformation family
of rational curve $C_1'$ such that $\pi(C_1')=C_0$
and that $L.C_1'=1$,
but we consider the deformation family
of rational curve $C_1'$
dominating a minimal extremal rational curve of $R$ on $M$
such that $L.C_1'$ takes
the minimum value 
among all such values,
namely, among all the values
$L.D$ where $D\subset P$ is a rational curve
dominating a minimal extremal rational curve of $R$ on $M$.
This is one of the key points of this paper.
Then we show in Lemma~\ref{existence of C_1'}
that the rational curve $C_1'$ chosen in this way
is also unsplit.
This observation is easy but plays a crucial role
in Lemma~\ref{comparisonshin}.
In fact, in Lemma~\ref{comparisonshin},
we will consider some deformation family
of unsplit rational curve $C_1'$
dominating $C_0$,
and we do not assume the minimality of $L.C_1'$.
Except for the points mentioned above,
Lemma~\ref{comparisonshin} is 
nothing but a reformulation
of \cite[Comparison Lemma (3.1)]{psw}
in our setting.
Lemma~\ref{comparisonshin} gives
a sufficient condition for
$R_1$ to contain
an unsplit rational curve $C_1'$
dominating $C_0$,
and in fact gives
a sufficient condition for
the affirmative answer to Problem~\ref{C_1 C_0 comparison}.
Indeed we have the following.
\begin{lemma}\label{C_1'=C_1}
Suppose that $R_1$ contains 
an unsplit rational curve $C_1'$
dominating a minimal extremal rational curve $C_0$ of $R$.
Then $C_1'$ is a minimal extremal rational curve in $R_1$,
i.e., we may assume that $C_1'=C_1$,
and thus we have $\pi(C_1)=C_0$.
\end{lemma}
\begin{proof}
Since $C_1$ and $C_1'$ are numerically propositional,
we have $C_1=\mu C_1'$ for some positive real number
$\mu$. We have $\mu\leq 1$ by the minimality of $C_1$.
Since both $C_1$ and $C_1'$ are unsplit,
Corollary~\ref{C_1'to pi(C_1') birational}
implies that $\pi_*(C_1)=\pi(C_1)$
and that $\pi_*(C_1')=\pi(C_1')$.
Hence we have $\pi(C_1)=\mu \pi(C_1')=\mu C_0$.
Here we see that $\mu\geq 1$ by the minimality of $C_0$.
Therefore we have $\mu=1$, namely, $C_1'$ can be regarded
as a minimal extremal rational curve $C_1$.
\end{proof}

In \S~\ref{Slr1=n+1},
we will apply Lemma~\ref{comparisonshin}
to the case $l(R_1)=n+1$
to obtain Proposition~\ref{l(r1)=n+1new}.
In \S~\ref{S|compaforlr1=n},
we will apply Lemma~\ref{comparisonshin}
to the case $l(R_1)=n$.
In this case, we will give an affirmative answer 
Corollary~\ref{cor for l(R_1)=n}
to the problem~\ref{C_1 C_0 comparison}
under the additional assumptions
that $L.C_1=1$ and that $\tau\geq 1$.
I have not succeeded in removing these additional assumptions.
So we will assume that $\tau\geq 1$ in the following,
and we will divide the case according to the value of $L.C_1$.

In \S~\ref{LC1=2shin}, we will deal with the case $L.C_1\geq 2$
and $\tau r\geq n-2$. 
In \S~\ref{PrelimfordimS=0},
we recall some preliminaries to deal with the following cases.

If $L.C_1=1$, we have 
the following advantage:
we have $(K_P+l(R_1)L).C_1=0$
for an ample line bundle $L$.
Hence to investigate the property of 
$(P,L)$ via $\varphi$ 
becomes much easier, in general,
than to investigate the property of 
$(M,\mathcal{E})$ via $\psi$.
Therefore we study the structure of $\varphi$
first; this strategy and idea
stem from Ye-Zhang \cite{yz} and Peternell \cite{p0}.
If $\tau r> n-3$,
we have $l(R_1)\geq n-2$ by (\ref{lb of r1=}) above.
We will divide the case $L.C_1=1$ according to 
the value of $l(R_1)$.
We will give the classification in case $l(R_1)=n$
in \S~\ref{Slr1=n},
based on the affirmative answer 
Corollary~\ref{cor for l(R_1)=n}
to the problem~\ref{C_1 C_0 comparison}.
In \S~\ref{shin prelim4} and \S~\ref{Sl(r1)=s'},
we prove some results needed to deal with the case $l(R_1)\leq n-1$.
In \S~\ref{Scompaforlr1n-1},
we will apply Lemma~\ref{comparisonshin}
to the case $l(R_1)=n-1$
to obtain an affirmative answer 
Corollary~\ref{comparison for lr1=n-1}
to the problem~\ref{C_1 C_0 comparison}
under some conditions.
Then applying Corollary~\ref{comparison for lr1=n-1}
we will give the classification in case $l(R_1)=n-1$
in \S~\ref{Slr1=n-1}.
Finally we will deal with the case $l(R_1)=n-2$
in \S~\ref{-K=n-2shin}.

\section{Comparison Lemma}\label{section comparison}
Let $M$ be an $n$-dimensional Fano manifold
with $\Pic M\cong \mathbb{Z}$,
and $\mathcal{E}$ an ample vector bundle
of rank $r$ on $M$.
Denote by $\tau$ the nef value of the polarized manifold
$(M,\det\mathcal{E})$.
Let $R$, $\pi:P\to M$, $L$, $R_1$, $C_1$, $\varphi:P\to N$ be 
as in the common setup \S~\ref{the common setup}.
In particular, we assume that $-K_P$ is ample
in this section.

The following observation is one of 
the key points of this paper.
\begin{lemma}\label{existence of C_1'}
Let $\lambda$ be the smallest integer
among $L.C_1'$'s,
where $C_1'$ moves among rational curves
dominating minimal extremal rational curves in $R$.
By abuse of notation, denote by $C_1'$
a rational curve dominating a minimal extremal rational curve 
in $R$ and attaining the number $\lambda$: $L.C_1'=\lambda$.
Then $C_1'$ is unsplit.
\end{lemma}
\begin{proof}
Suppose, to the contrary, that 
$C_1'$ is effectively algebraically equivalent
to a sum $\Sigma_{i=1}^{\delta} D_i$
of $\delta$ ($\delta\geq 2$) rational curves $D_i$,
some of which may equal.
Then $\pi_*(C_1')$
is effectively algebraically equivalent to 
$\Sigma_{i=1}^{\delta} \pi_*(D_i)$.

We claim here that $C_1'\to \pi(C_1')$ is birational.
Let $\tilde{\pi(C_1')}\to \pi(C_1')$ be the normalization.
We see that $\mathcal{E}\otimes \mathcal{O}_{\tilde{\pi(C_1')}}$
decomposes into a direct sum of line bundles
of degree $\geq \lambda$
by the minimality of $\lambda$.
Let $d$ be the degree of the morphism $C_1'\to \pi(C_1')$,
and  $\tilde{C}_1'\to C_1'$ the normalization.
Since the induced map $\tilde{C}_1'\to \tilde{\pi(C_1')}$
has degree $d$,
we infer that $\mathcal{E}\otimes
\mathcal{O}_{\tilde{C}_1'}$
is a direct sum of line bundles
of degree $\geq d\lambda$.
Therefore every section of 
$\mathbb{P}(\mathcal{E}\otimes
\mathcal{O}_{\tilde{C}_1'})\to \tilde{C}_1'$
has degree $\geq d\lambda$ with respect to 
the tautological line bundle.
On the other hand, $C_1'\subset P$ defines
a section of $\mathbb{P}(\mathcal{E}\otimes
\mathcal{O}_{\tilde{C}_1'})\to \tilde{C}_1'$,
which we also denote by $\tilde{C}_1'$ by abuse of notation,
such that $\tilde{C}_1'$ 
has degree $\lambda$ with respect to 
$L$.
Hence we conclude that $d=1$.

Note here that $\pi(C_1')$ is unsplit,
since $\pi(C_1')$ is a minimal extremal rational curve.
Therefore the claim above implies that 
except for one rational curve, say $\pi(D_1)$,
every $\pi(D_i)$ is a point,
that $D_1\to \pi(D_1)$ is birational,
and that $\pi(D_1)$ is also a minimal extremal rational curve.
Now we see that $L.D_1<\lambda$ since $\delta>1$ and $L$ is ample.
This contradicts the minimality of $\lambda$.
\end{proof}

The following lemma plays a crucial role in the study
of the case $\dim S=0$;
it is an extension
of \cite[Comparison Lemma (3.1)]{psw}.

\begin{lemma}\label{comparisonshin}
Let $C_1'$ be an unsplit rational curve on $P$
with $\pi(C_1')$ a minimal extremal rational curve
in $R$.
Set $C_0=\pi(C_1')$.
Let $T'$ be the connected component of 
the Hom scheme $\Hom_{bir} (\mathbb{P}^1, P)$
containing the normalization 
$\mathbb{P}^1=\tilde{C}_1'\to C_1'\subset P$.
Let $\tilde{T}'$ be the normalization of $T'$
and $T$ the image of $\tilde{T}'$
via the morphism $\Hom_{bir}^n(\mathbb{P}^1,P)\to \RatCurves^n(P)$
(see \cite[I.~(6.9), II.~(2.11),
and II.~(2.15)]{kollar} for
the morphism $\Hom_{bir}^n(\mathbb{P}^1,P)\to \RatCurves^n(P)$
and the notation).
Note that $T$ is proper since $C_1'$ is unsplit.
Let $V\to T$ be the universal family,
i.e., the restriction of $\Univ^{rc}(P)\to \RatCurves^n(P)$,
which is $\mathbb{P}^1$-bundle by \cite[II.~(2.12)]{kollar}.
Let $p: V\hra P\times T\ra P$ and $q: V\hra P\times T\ra T$
be the canonical projections.
Pick a point $x\in C_0$ and fix it. 
Let $T_x$ denote
an irreducible component
of $q((\pi\circ p)^{-1}(x))$ containing $[C_1']$,
and $V_x$ denote $q^{-1}(T_x)$.
Note that $T_x$ is proper since $T$ is so.
We also denote by $p_x$ and $q_x$ the projections
$V_x\ra P$ and $V_x\ra T_x$ respectively.
We see that $q_x$ is $\mathbb{P}^1$-bundle
since $q$ is so.
Set 
\[t=\min \{\dim T_i\, |\, T_i \textrm{ is an
irreducible component of $T$ containing $[C_1']$}\}.\]
Note that Lemma~\ref{big r+1} enables us
to set $\det \mathcal{E}.\pi_*(C_1')
=rL.C_1'+\alpha'$ for some non-negative integer $\alpha'$.
Then we have the following.
\begin{enumerate}
\item 
      If $\varphi$ has an $n$-dimensional fiber,
      then $\varphi$ contracts $C_1'$,
      i.e., $C_1'$ is an unsplit rational curve belonging to $R_1$.
\item 
      If $\dim p_x(V_x)\cap F\geq 1$ 
      for some positive dimensional fiber $F$ of $\varphi$,
      then $\varphi$ contracts $C_1'$.
\item 
      If $p_x(V_x)\cap F\neq \emptyset$
      for some positive dimensional fiber $F$ of $\varphi$,
      we have $\dim p_x(V_x)\cap F\geq 
      \dim p_x(V_x)+\dim F-\dim P
      \geq \dim p_x(V_x)+l(R_1)-1+\codim (E(R_1),P)-\dim P$.
\item We have $\dim p_x(V_x)=\dim V_x$.
\item 
      We have $\dim V_x\geq t+2-n$ and 
      $t\geq \tau rL.C_1'+(\tau -1)\alpha'+\dim P-3$.
      If $\tau \geq 1$, we have therefore
      $t\geq \tau rL.C_1'+\dim P-3$.
\end{enumerate} 
\end{lemma}
\begin{proof}
(1) If $\varphi$ has $n$-dimensional fibers,
then $\varphi$ contracts some curve dominating $C_0$.
On the other hand, 
we see,
by Lemma~\ref{big r+1},
that the quotient bundle 
$L_{C_1'}\otimes \mathcal{O}_{\tilde{C}_1'}$
of $\pi^*(\mathcal{E})_{C_1'}\otimes 
\mathcal{O}_{\tilde{C}_1'}$
has the minimal degree among all quotient bundles,
since $C_1'$ is unsplit with $\pi(C_1')=C_0$ a curve.
Therefore if $\varphi$ contracts some curve dominating $C_0$,
then $\varphi$ contracts $C_1'$.

(2) Suppose that $p_x(V_x)\cap F$ contains a curve.
We will show that
some curve in $p_x(V_x)\cap F$ corresponds
to a point in $T_x$. 
If this claim holds, we may take as $C_1'$ 
the curve in $p_x(V_x)\cap F$ so that 
we infer that $C_1'$ is contracted by $\varphi$.
Suppose, to the contrary, that there is 
no curve in $p_x(V_x)\cap F$
which corresponds
to a point in $T_x$,
i.e., that no curve in $p_x^{-1}(F)$ is contracted by $q_x$.
Let $B_1'$ be an irreducible closed curve in $p_x^{-1}(F)$,
and let $B'$ be the image $q_x(B_1')$ in $T_x$.
Since no curve in $p_x^{-1}(F)$ is contracted by $q_x$, 
we see that $B'$ is a curve.
Let $B\to B'$ be the normalization.
Note that $B$ is proper since $T_x$ is so.
Set $S=V_x\times_{T_x}B$.
We see that the projection $S\to B$ is $\mathbb{P}^1$-bundle
since $q_x$ is so.
We have $N_1(S)\cong \mathbb{R}^2$.
Let $B_1$ be the image in $S$
of the normalization $\tilde{B_1'}$ of $B_1'$.
First look at the morphism $S\to P\to N$.
We see that the image of $B_1$ in $P$ is contained in $F$
so that it is contracted by $\varphi$.
Since no curve in $p_x^{-1}(F)$ is contracted by $q_x$, 
we infer that 
the image in $P$ of any fiber of $S\to B$
intersects $F$ but does not contained in $F$.
Hence the image in $P$ of any fiber of $S\to B$
does not contracted by $\varphi$.
Therefore we have 
$\overline{{\rm NE}}(S)=
\mathbb{R}_{\geq 0}[B_1]
+\mathbb{R}_{\geq 0}[\textrm{a fiber of }S\to B]$.
Next look at the morphism $S\to P\to M$.
Since $F\to \pi(F)$ is finite by Lemma~\ref{finite of yz},
the image in $P$ of $B_1$, which is contained in $F$,
does not contracted by $\pi$.
Note here that 
the image in $P$ of any fiber of $S\to B$ is
numerically equivalent to $C_1'$,
which does not contracted by $\pi$.
Hence we infer,
by the property of the contraction morphism $\pi$
of an extremal ray $R_{\pi}$,
that the image in $P$ of any fiber of $S\to B$
does not contracted by $\pi$.
This implies that the pull back of an ample line
bundle on $M$ defines a positive function
on $\overline{{\rm NE}}(S)$.
Therefore $S\to M$ is finite.
On the other hand, the image in $M$
of any fiber of $S\to B$ passes through $x$,
so that $S\to M$ is not finite over $x$.
This is a contradiction.
This completes the proof of (2).

(3) Since $P$ is smooth, we see that 
$\dim p_x(V_x)\cap F\geq 
\dim p_x(V_x)+\dim F-\dim P$ if 
$p_x(V_x)\cap F\neq \emptyset$.
We apply (\ref{ub of r1}) in \S~\ref{the common setup}
to see that $\dim p_x(V_x)+\dim F-\dim P\geq
\dim p_x(V_x)+l(R_1)-1+\codim (E(R_1),P)-\dim P$.

(4)
We will show more precisely that 
$V_x\setminus (\pi\circ p)^{-1}(x) \to
p_x(V_x)\setminus \pi^{-1}(x)$ is finite.
Suppose to the contrary that we could find a curve 
$B_1'\subset V_x\setminus (\pi\circ p)^{-1}(x)$
over a point $y\in p_x(V_x)\setminus \pi^{-1}(x)$.
Set $B'=q_x(B_1')\subseteq T_x$.
Since the image in $P$ of any fiber of $q_x:V_x\to T_x$
is numerically equivalent to $C_1'$,
we see, by the same reason as in (2), that 
the image of any fiber of $q_x$ is not contracted by $\pi$.
Hence we infer that $B'$ is a curve.
Let $B\ra B'$ be the normalization,
and set $S=V_x\times_{T_x}B$.
Denote by $B_1$ the image in $S$ of the normalization
$\tilde{B_1'}$ of $B_1'$.
Note here that $V_x\setminus (\pi\circ p)^{-1}(x) \to
p_x(V_x)\setminus \pi^{-1}(x)$ is proper since $T_x$ is so.
Therefore $B$ is also proper.
Now look at the morphism $\alpha :S\to P$.
We see that $\alpha(S)$ is a surface
and that the image by $\alpha$
of any fiber of the $\mathbb{P}^1$-bundle
$S\to B$ passes through $y$.
Moreover we infer that 
$y$ is the only fixed point of 
the family of rational curves 
defined as the images of the fibers of $S\to B$,
since $B$ is proper.
Therefore $\dim \alpha(S)\cap \pi^{-1}(x)=1$
since $y\notin \pi^{-1}(x)$. 
Now that $N_1(S)\cong \mathbb{R}^2$
and that $\alpha$ contracts $B_1$,
we see that $N_1(\alpha(S))\cong \mathbb{R}^1$.
Since $\pi$ contracts all 
curves in $\alpha(S)\cap \pi^{-1}(x)$,
this implies that $\pi$ contracts all $\alpha(S)$,
i.e., that $\pi(\alpha(S))$ is a point.
This however contradicts $\pi(y)\neq x$.

(5) We have $\dim V_x=\dim T_x+1$.
First we see by the same reason as stated in (4)
that the image of any fiber of $q$ is not contracted by $\pi$.
Hence there exists no rational curve contracted by 
$q$ in $(\pi\circ p)^{-1}(x)$.
This implies that $T_x$ has dimension equal to 
the corresponding irreducible
component of $(\pi\circ p)^{-1}(x)$.
Note here that any irreducible component of $V$
containing $q^{-1}([C_1'])$ has dimension
at least $t+1$.
Hence every irreducible component of 
$\dim (\pi\circ p)^{-1}(x)$
has dimension $\geq t+1-n$;
thus we have $\dim T_x\geq t+1-n$.
Therefore we have $\dim V_x\geq t+2-n$.

Next we have 
$t\geq -K_P.C_1'+\dim P-\dim PGL(2,\mathbb{C})$.
Since $C_1'$ is unsplit, it follows from
Lemma~\ref{big r+1} that 
\begin{equation*}
\begin{split}
-K_P.C_1'&= rL.C_1'+(\tau -1)\det \mathcal{E}.\pi_*(C_1')\\
      &= rL.C_1'+(\tau -1)(rL.C_1+\alpha')\\
      &= \tau rL.C_1'+(\tau -1)\alpha'.
\end{split}
\end{equation*}
Therefore we have 
$t\geq \tau rL.C_1'+(\tau -1)\alpha'+\dim P-3$.
\end{proof}

\begin{cor}
If $\tau \geq 1$, $\tau r\geq (n+3)/2$, and 
$\varphi$ is of fiber type,
then $\pi(C_1)$ is a minimal extremal rational curve in $R$.
In particular, we have $l(R)\geq l(R_1)$.
\end{cor}

\section{The case $l(R_1)=n+1$}\label{Slr1=n+1}
Let $M$, $\mathcal{E}$, $\psi:M\to S$, and $C_0$
be as in \S~\ref{shinprelim}.
Suppose that $\dim S=0$.
Let $\pi:P\to M$ and $L$
be as in \S~\ref{the common setup},
and suppose that $-K_P$ is ample.
Let $R_1$, $\varphi:P\to N$, and $C_1$
be as in \S~\ref{the common setup}.
Then we have the following.
\begin{prop}\label{l(r1)=n+1new}
If $l(R_1)=n+1$, then 
$(M,\mathcal{E})\cong
(\mathbb{P}^n,\mathcal{O}(l)^{\oplus r})$
for some positive integer $l$.
\end{prop}
\begin{proof}
Since $l(R_1)=n+1$, inequality~(\ref{ub of r1}) 
in \S~\ref{the common setup}
implies that $\varphi$ is of fiber type
and that every fiber of $\varphi$ is $n$-dimensional.
Then Lemmas~\ref{comparisonshin} (1)
and \ref{C_1'=C_1} implies that $\pi(C_1)=C_0$.
Therefore we have $l(R)=l(R_1)+
(\det\mathcal{E}.C_0-rL.C_1)\geq l(R_1)$ 
by Lemma~\ref{lr1-lr=a}.
Since $l(R)\leq n+1$ by (\ref{upper bound of l(R)})
in \S~\ref{shinprelim}
and $l(R_1)=n+1$,
this implies that $l(R)=n+1$
and $rL.C_1=\det\mathcal{E}.C_0$.
Theorem~\ref{cm} then implies that $M\cong \mathbb{P}^n$;
thus $C_0$ is a line in $\mathbb{P}^n$.
Lemmas~\ref{big r+1} and \ref{C_1'to pi(C_1') birational} 
also implies that $\mathcal{E}|_{C_0}
\cong \mathcal{O}(L.C_1)^{\oplus r}$.
Set $l=L.C_1$.
Then $\mathcal{E}$ is a uniform vector bundle of 
type $(l,\dots,l)$,
so that we have $\mathcal{E}\cong \mathcal{O}(l)^{\oplus r}$
by Theorem~\ref{uniform} (1).
\end{proof}

\section{Comparison Lemma for $l(R_1)=n$}\label{S|compaforlr1=n}
In this section, we will follow the notation
in \S~\ref{section comparison}.
\begin{cor}\label{cor for l(R_1)=n}
Suppose that $\tau \geq 1$,
that $l(R_1)=n$, that $L.C_1=1$,
and that $\varphi$ has no 
$n$-dimensional fibers.
Then $\pi(C_1)$ is a minimal extremal rational curve
in $R$.
\end{cor}
\begin{proof}
By Lemma~\ref{C_1'=C_1},
it is enough to show that 
there exists 
an unsplit rational curve $C_1'$
belonging to $R_1$
and dominating  
a minimal extremal rational curve in $R$.

Assume, to the contrary, that 
there does not exist such a curve.
By abuse of notation, denote by $C_1'$ 
an unsplit rational curve dominating  
a minimal extremal rational curve $C_0$ in $R$,
whose existence is guaranteed by Lemma~\ref{existence of C_1'}.
Then $C_1'$ does not belong to $R_1$,
i.e., $\varphi$ does not contract $C_1'$,
by assumption.
Note here that 
$\varphi$ is of fiber type 
and that every fiber of $\varphi$ is $(n-1)$-dimensional
by inequality (\ref{ub of r1}) in 
\S~\ref{the common setup},
since $l(R_1)=n$
and $\varphi$ has no $n$-dimensional fibers.
Lemma~\ref{comparisonshin} (2), (3), (4), and (5)
then implies that 
\begin{equation*}
\begin{split}
0\geq \dim p_x(V_x)\cap F&\geq \dim p_x(V_x)-r\\
&\geq \tau rL.C_1'+(\tau -1)\alpha'-2\\
&\geq \tau rL.C_1'-2,
\end{split}
\end{equation*}
since $\tau \geq 1$.
Hence we have $2\geq \tau rL.C_1'\geq \tau r\geq r$.
Note here that $r\geq 2$
since $\varphi$ has no $n$-dimensional fibers.
Therefore we see that $r=2$, that $\tau =1$, 
and that $L.C_1'=1$.
Since $L.C_1=1$,
it follows from (\ref{lb of r1=}) in 
\S~\ref{the common setup}
that $n=l(R_1)=2L.C_1=2$.
Hence $M$ is a Del Pezzo surface with 
Picard number one.
Thus we have $M\cong \mathbb{P}^2$
by the classification of Del Pezzo surfaces.
Now we have $\det\mathcal{E}\cong \mathcal{O}(3)$.
Since $\mathcal{E}$ is ample,
Theorem~\ref{uniform} implies that
$\mathcal{E}\cong \mathcal{O}(1)\oplus \mathcal{O}(2)$
or $T_{\mathbb{P}^2}$.
In these cases,
$C_1$ is an unsplit rational curve 
belonging to $R_1$
and dominating a minimal extremal rational curve in $R$.
This contradicts the assumption.
\end{proof}

\section{The case where $L.C_1\geq 2$.}\label{LC1=2shin}
Let $M$, $\mathcal{E}$, $\psi:M\to S$, $C_0$,
$\pi:P\to M$, and $L$ be as in \S~\ref{the common setup}.
In particular, we assume that $\dim S=0$.
Suppose that $\tau\geq 1$,
and let
$R_1$, $C_1$, and $\varphi:P\to N$ 
be as in \S~\ref{the common setup}.
In this section, we will deal with the case $L.C_1\geq 2$.

We assume $\tau r\geq n-2$ in this section.
Since $\tau\geq 1$,
it follows from inequality~(\ref{lb of r1=}) in 
\S~\ref{the common setup}
that 
$l(R_1)\geq \tau rL.C_1$.
On the other hand, we have an upper bound 
$l(R_1)\leq n+1$
by inequality~(\ref{ub of r1}) in 
\S~\ref{the common setup}.
Hence we have 
\begin{equation}\label{r(L.C_1-1)bound}
r(L.C_1-1)
\leq \tau r(L.C_1-1)
\leq l(R_1)-
\tau r
\leq 3,
\end{equation}
since $\tau\geq 1$.
In particular, we see that $r\leq 3$ if $l(R_1)=n+1$,
that $r\leq 2$ if $l(R_1)=n$,
and that $r=1$ if $l(R_1)=n-1$.

Suppose that $l(R_1)=n+1$.
Then $(M, \mathcal{E})\cong 
(\mathbb{P}^n, \mathcal{O}(l)^{\oplus r})$,
where $l=L.C_1$,
by Proposition~\ref{l(r1)=n+1new}.
Therefore $\tau r=(n+1)/l$.
Since $\tau r\geq n-2$, we have $n+1\geq l(n-2)$,
i.e., $2l+1\geq (l-1)n$.
Since $l=L.C_1\geq 2$, we have 
\begin{equation}\label{upper bound of n in case L.C_1 big 2}
2+\frac{3}{l-1}\geq n.
\end{equation}
Since $\tau \geq 1$, we also have $n+1=\tau rl\geq rl$, i.e.,
\begin{equation}\label{lower bound of n in case L.C_1 big 2}
n\geq rl-1.
\end{equation}

Suppose moreover that $r=3$.
We have $L.C_1=2$, $\tau=1$,
and $\tau r=n-2$ by (\ref{r(L.C_1-1)bound}).
Thus $n=5$. 
Hence we have $(M, \mathcal{E})\cong 
(\mathbb{P}^5, \mathcal{O}(2)^{\oplus 3})$.
This is a special case of the case (17) of 
Theorem~\ref{Main theorem}.

Suppose moreover that $r=2$.
We have $L.C_1=2$ 
and $n-1\geq \tau r$
by (\ref{r(L.C_1-1)bound}).
Thus $n\geq 3$.
Moreover we have $n\leq 5$
by (\ref{upper bound of n in case L.C_1 big 2}).
Therefore $(M, \mathcal{E})\cong
(\mathbb{P}^n, \mathcal{O}(2)^{\oplus 2})$,
where $3\leq n\leq 5$.
This is a special case of the case (14) of 
Theorem~\ref{Main theorem}.

Suppose moreover that $r=1$.
We have $2\leq L.C_1\leq 4$ by (\ref{r(L.C_1-1)bound}).
Hence inequalities~(\ref{upper bound of n in case L.C_1 big 2})
and (\ref{lower bound of n in case L.C_1 big 2}) 
imply that 
$(M,\mathcal{E})$
is isomorphic to either
of the following;\\
a) $(\mathbb{P}^3,\mathcal{O}(4))$;\\
b) $(\mathbb{P}^n,\mathcal{O}(3))$ $(n=2, 3)$;\\
c) $(\mathbb{P}^n,\mathcal{O}(2))$ $(1\leq n\leq 5)$.\\
The case a), b), or c) is, respectively,
a special case of the case (17), (13), or (4)
of Theorem~\ref{Main theorem}.

Suppose that $l(R_1)=n$.

Suppose moreover that $r=2$.
It follows from (\ref{r(L.C_1-1)bound}) above
that $L.C_1=2$,
that $\tau =1$, and that $\tau r=n-2$.
Thus $n=4$.
This is a special case of the case (18) of Theorem~\ref{Main theorem}.

Suppose moreover that $r=1$.
We have $L.C_1=2$ or $3$
by (\ref{r(L.C_1-1)bound}) above. 
Set $e=\det \mathcal{E}.C_0=L.C_1$.
Since $n-2\leq \tau r\leq \tau e=l(R)=l(R_1)=n$,
we have $(n-2)e\leq n$.
Hence we have $(e-1)n\leq 2e=2(e-1)+2$,
and thus 
\[n\leq 2+\dfrac{2}{e-1}.\]
On the other hand, we have $e\leq \tau e=n$.
Hence we have $n=3$ if $e=3$ and $2\leq n\leq 4$ if $e=2$.
If $(n,e)=(3,3)$ or $(4,2)$, we have $K_M+(n-2)\mathcal{E}=0$,
and this is a special case of the case (17)
of Theorem~\ref{Main theorem}.
If $(n,e)=(3,2)$,
setting $A=-(K_M+\mathcal{E})$,
we see that $K_M+3A=0$.
Theorem~\ref{ko} (2) then implies that 
$(M,A)\cong (\mathbb{Q}^3,\mathcal{O}(1))$.
Thus we have $\mathcal{E}\cong \mathcal{O}(2)$.
This is a special case of the case (15) of 
Theorem~\ref{Main theorem}.
If $(n,e)=(2,2)$, 
$M$ is a Del Pezzo surface with $\rho (M)=1$,
and thus $M\cong \mathbb{P}^2$ by 
the classification.
This contradicts $l(R)=n$.
Hence this case does not occur.

Suppose that $l(R_1)=n-1$.
Then it follows from
(\ref{r(L.C_1-1)bound}) that $L.C_1=2$, that $r=1$,
that $\tau =1$,
and that $\tau r=n-2$.
Hence we have $n=3$.
Therefore we have $K_M+L=0$.
This is a special case of the case (17) of Theorem~\ref{Main theorem}.

\section{Preliminaries for the case $\dim S=0$}\label{PrelimfordimS=0}
In this section, we will recall some of the results
in \cite{yz}, \cite{fq}, and \cite{psw}
in the form useful in the following context.

First recall a remarkable argument in \cite[\S 4]{psw};
since the statement (2) in the following lemma is not stated
in this form in \cite[\S 4]{psw},
we attach its proof.
We also make a little modification applicable to a Brauer-Severi scheme,
i.e., a projective space bundle in the \'{e}tale topology.
\begin{lemma}\label{isom}
Let $P$ be a Fano manifold,
and suppose that 
$\overline{{\rm NE}}(P)$
has two different extremal rays $R_\pi$ and $R_1$.
Let $\pi:P\to M$ be the contraction morphism of $R_\pi$,
and $\varphi:P\to N$ that of $R_1$.
Suppose that every 
closed fiber of $\pi$ is isomorphic to $\mathbb{P}^{r-1}$,
and that $M$ is an $n$-dimensional Fano manifold of Picard number one.
If $\varphi$ has an $n$-dimensional closed fiber $\varphi^{-1}(z)$,
then $f:W\to A_z$ denoting the normalization of an
$n$-dimensional irreducible component $A_z$ of 
$\varphi^{-1}(z)$, the composite $(\pi|_{A_z})\circ f:W\to M$
of $\pi|_{A_z}$ and $f$ 
is an isomorphism. 
\end{lemma}
\begin{proof}
First note that 
$(\pi|_{A_z})\circ f:W\to M$ is finite by 
Lemma~\ref{finite of yz}.
Moreover it is surjective since $\dim W=n=\dim M$.
Let $T$ be the singular locus of 
$W$ and $S$ the image of $T$ via $(\pi|_{A_z})\circ f$.
We show first that $(\pi|_{A_z})\circ f$ is not ramified
over $M\setminus S$; let 
$h: W\setminus ((\pi|_{A_z})\circ f)^{-1}(S)\to M\setminus S$ 
be the restricted morphism
and $R$ the ramification divisor of $h$.
Since $(\pi|_{A_z})\circ f$ is finite and surjective, 
$S$ has codimension $\geq 2$ in $M$; 
noting $\Pic M\cong \mathbb{Z}$, 
we can apply \cite[Lemma~2]{kmmo}
and find out that a general extremal rational curve $C$
does not meet $S$. 
We may also assume that $C$ is not contained in $h(R)$;
since $\Pic M\cong \mathbb{Z}$ and $C\cap S=\emptyset$,
we see that $C\cap h(R)$ is a finite set
and not empty if $R$ is not.
A fiber of $h$ over a point in $C\cap h(R)$ consists of
fewer points than 
a fiber of $h$ over a point in $C\setminus h(R)$ does.
On the other hand, 
any morphism from a rational scroll,
morphism which
is not finite and 
does not contract
any fiber of the projection of the scroll,
has a positive dimensional fiber consisting 
of a disjoint union of sections and isolated points.
Let $\mathbb{P}^1\to M$ be the composite of the normalization
$\mathbb{P}^1\to C$ and the inclusion $C\hookrightarrow M$.
Then $\mathbb{P}^1\times_MP\to \mathbb{P}^1$
is a scroll, i.e.,  a $\mathbb{P}^{r-1}$-bundle in the 
Zariski topology, since $H^2(\mathbb{P}^1,\mathcal{O}^{\times})=0$.
Hence the function $C\ni x\mapsto \# \{\pi^{-1}(x)\cap A_z\}$
is upper-semicontinuous on $C$.
Thus $C\ni x\mapsto \# h^{-1}(x)$ is also upper-semicontinuous,
and therefore $R$ is empty. Hence $h$ is \'{e}tale.

By the purity of the branch locus, $(\pi|_{A_z})\circ f$
is also \'{e}tale and $W$ is smooth. 
Since $M$ is simply connected
by \cite{kmmorc}, we infer that $(\pi|_{A_z})\circ f$ is birational;
by Zariski Main Theorem,
we conclude that $(\pi|_{A_z})\circ f$ is an isomorphism.
\end{proof}

The following lemma is 
a slight modification of \cite[Prop. 4.2]{psw}.

\begin{lemma}\label{standard}
Let $M$ be a Fano manifold of Picard number one,
and $\mathcal{E}$ an ample vector bundle of rank $r$ on $M$.
Denote by $P$ the projective space bundle
$\mathbb{P}(\mathcal{E})$,
by $\pi:P\to M$ the projection,
and  by $L$ the tautological line bundle $H(\mathcal{E})$.
Let $R_{\pi}$ be the extremal ray corresponding to $\pi$.
Suppose that there exists a contraction morphism
$\varphi :P\ra N$ of an extremal ray $R_1$
different from $R_{\pi}$,
and that general fibers of $\varphi$ have dimension $<n$,
and that $\varphi$ has also an $n$-dimensional fiber $\varphi^{-1}(z)$.
\begin{enumerate}
\item If $-K_P-(n-1)L$ is $\varphi$-nef, then
      $M\cong \mathbb{P}^n$ or $\mathbb{Q}$.
\item Moreover if $-K_P-nL$ is $\varphi$-nef, 
      then $M\cong \mathbb{P}^n$.
\end{enumerate}
\end{lemma}
\begin{proof}
As in Lemma~\ref{isom},
let $f:W\to A_z$ denote the normalization of an
$n$-dimensional irreducible component $A_z$ of 
$\varphi^{-1}(z)$.
Then $W\cong M$ by Lemma~\ref{isom}
and thus $W$ is smooth.
Hence we have $h^n(W,t{f}^*(L|_{A_{z}}))=0$ for all $t\geq -(n-1)$
by \cite[Lemma 4]{yz},
since $-K_P-(n-1)L$ is $\varphi$-nef.
By Kodaira vanishing
we also have $h^i(W,t{f}^*(L|_{A_{z}}))=0$ for all $i<n$ and $t<0$.
Furthermore we have $h^i(W,\mathcal{O})=h^i(M,\mathcal{O})=0$
for $i>0$ since $M$ is Fano. Hence 
the Hilbert polynomial $\chi (W, t{f}^*(L|_{A_{z}}))$
is of the form 
\[\chi (W, t{f}^*(L|_{A_{z}}))=
\frac{d_0}{n!}(\prod_{k=1}^{n-1}(t+k))(t+\frac{n}{d_0})\]
where $d_0=\deg {f}^*(L|_{A_{z}})$;
thus we have $\chi (W, {f}^*(L|_{A_{z}}))=d_0+n$.
Since $W$ is Fano and hence $\chi (W, {f}^*(L|_{A_{z}}))
=h^0(W, {f}^*(L|_{A_{z}}))$, we infer that 
Fujita's $\Delta$-genus of $(W, {f}^*(L|_{A_{z}}))$ is zero.
Note here that the Picard number of $W$ is one
since that of $M$ is so; 
we therefore conclude that $W$ is isomorphic to either
$\mathbb{P}^n$ or $\mathbb{Q}^n$.
Thus $M\cong \mathbb{P}^n$ or $\mathbb{Q}^n$

If $-K_P-nL$ is $\varphi$-nef, we see that
$\chi (W, -n{f}^*(L|_{A_{z}}))=0$ by \cite[Lemma 4]{yz};
thus the same argument as above implies that $d_0=1$
and $M\cong \mathbb{P}^n$.
\end{proof}

Recall the following lemma of 
Fujita~\cite[(2.12) Lemma]{fs}.
\begin{lemma}\label{fujita lemma}
Let $\varphi :P\to N$ be a proper surjective morphism
from a manifold $P$ onto a normal variety $N$
with equidimensional fibers.
Let $L$ be a $\varphi$-ample line bundle on $P$
and suppose that $(F,L|_{F})\cong (\mathbb{P}^r,\mathcal{O}(1))$
for a general fiber $F$ of $\varphi$.
Then $N$ is nonsingular and $\varphi$ makes $(P,L)$
a scroll over N.
\end{lemma}

The following is due to Fujita~\cite[(2.2) Theorem]{fq}
and Ye-Zhang~\cite[Lemma 4]{yz}.
\begin{lemma}\label{Fujita Ye-Zhang}
Let $\varphi :P\to N$ be a contraction morphism
of an extremal ray $R_1$ from 
a manifold $P$ onto a normal variety $N$.
Let $L$ be a $\varphi$-ample line bundle on $P$
and suppose that $(K_P+sL).R_1=0$ for some
positive integer $s$.
Let $A_z$ be an $s$-dimensional irreducible component of 
an $s$-dimensional fiber $\varphi^{-1}(z)$ of $\varphi$
and suppose that a general fiber of $\varphi$
has dimension $<s$.
Then we have $(W,f^*(L|_{A_z}))
\cong (\mathbb{P}^s,\mathcal{O}(1))$,
where $f:W\to A_z$ is the normalization.
\end{lemma}
\begin{proof}
Let $\tilde{W}\to W$ be a desingularization
and let $g:\tilde{W}\to A_z$ be the composite of $f$ and this 
desingularization.
Then we have, by \cite[Lemma 4]{yz},
that 
$h^{s}(\tilde{W},-t{g}^*(L|_{A_{z}}))=0$ for all $t\leq s$
since $(K_P+sL).R_1=0$.
Hence, applying Fujita's theorem \cite[(2.2) Theorem]{fq}, 
we conclude that
$(W,{f}^*(L|_{A_{z}}))\cong (\mathbb{P}^{s},\mathcal{O}(1))$.
\end{proof}

\section{The case $l(R_1)=n$}\label{Slr1=n}
To deal with the case $l(R_1)=n$,
we will use the following:
\begin{lemma}\label{key}
Let $M$ be an $n$-dimensional projective manifold,
and $\mathcal{E}$ an ample vector bundle of rank $r$ on $M$.
Denote by $\pi :{\mathbb P}(\mathcal{E})\ra M$ the projection,
and by $L$ the tautological line bundle $H(\mathcal{E})$ on 
${\mathbb P}(\mathcal{E})$.
Assume that $\Pic M\cong {\mathbb Z}$.
Suppose that there exists a 
${\mathbb P}^{n-1}$-bundle
$\varphi :{\mathbb P}(\mathcal{E})\ra N$ onto
an $r$-dimensional 
projective manifold $N$
such that $\varphi |_{\pi^{-1}(x)}$ is finite for every
point $x\in M$ and $L|_F\cong \mathcal{O}_{\mathbb{P}^{n-1}}(1)$ for
every fiber $F$ of $\varphi$.
Then $r=n$ and
$(M,\mathcal{E})\cong ({\mathbb P}^n,T_{\mathbb{P}^n})$.
\end{lemma}
\begin{proof}
Let ${\mathcal F}$ denote $\varphi _*L$.
Then ${\mathcal F}$ is a vector bundle of rank $n$.
Moreover ${\mathcal F}$ is ample because $H({\mathcal F})=L$.
Note that $\Pic N\cong \mathbb{Z}$.
Denote by $P$ the projective space bundle
$\mathbb{P}(\mathcal{E})$ over $M$.

Since 
\[
-rL+\pi^*(K_M+\det \mathcal{E})=K_P=-nL+\varphi^*(K_N+\det {\mathcal F}),
\]
we have $n-r=\varphi^*(K_N+\det {\mathcal F}).l=
(K_N+\det {\mathcal F}).\varphi_*(l)$,
where $l$ denotes a line in a fiber of $\pi$.
Note that $l\ra \varphi (l)$ is birational because $L.l=1$.
Thus $-K_N.\varphi (l)=\det {\mathcal F}.\varphi (l)+r-n\geq r$.

Since the condition is symmetric with respect to $\pi$
and $\varphi$,
we may assume that $r\leq n$.
Denote by $\tau$ the nef value $\tau (M,\det\mathcal{E})$
of the polarized manifold $(M,\det\mathcal{E})$.
The condition $r\leq n$ implies that
$\tau \geq 1$,
since $L|_F=\mathcal{O}(1)$
for a fiber $F$ of $\varphi$.

We will assume that $r\leq n-1$,
and derive a contradiction.

We claim here that $-K_N.\varphi (l)\leq r+1$.
Assume, to the contrary, that $-K_N.\varphi (l)\geq r+2$.
Then $\varphi (l)$ can be deformed to a sum $\sum_{i=1}^{\delta} l_i$
of at least two rational curves $l_i$'s (some of which may be equal)
($i=1,\ldots,\delta, \delta\geq 2$) 
such that $-K_N.l_i\leq r+1$
by Mori's theorem \cite[Theorem 4]{mo1}.
Thus 
\[n-r=(K_N+\det {\mathcal F}).\varphi (l)=
\sum_{i=1}^{\delta} (K_N+\det {\mathcal F}).l_i\geq \delta (-r-1+n).\]
Hence $(\delta -1)(n-r)\leq \delta$.
Since $r\leq n-1$,
we have $1\leq n-r\leq 1+(1/(\delta -1))\leq 2$.
If $n-r=1$, then $1=(K_N+\det {\mathcal F}).\varphi (l)=
\sum_{i=1}^{\delta} (K_N+\det {\mathcal F}).l_i$, which
is a contradiction because $\Pic N\cong \mathbb{ Z}$ and so
$K_N+\det {\mathcal F}$ is ample.
Hence $n-r=2$, $\delta =2$, $(K_N+\det {\mathcal F}).l_i=1$,
$n=\det {\mathcal F}.l_i$, and $-K_N.l_i=r+1$;
thus we have $K_N+(r+1)(K_N+\det {\mathcal F})=0$.
Applying Theorem~\ref{ko} (1),
we infer that 
$(N, K_N+\det {\mathcal F})\cong (\mathbb{P}^r, \mathcal{O}(1))$.
Therefore $\det {\mathcal F}\cong \mathcal{O}(r+2)=\mathcal{O}(n)$
and ${\mathcal F}\cong \mathcal{O}(1)^{\oplus n}$
by Theorem~\ref{uniform} (1).
This means that $\pi$ is $\mathbb{P}^r$-bundle, which
contradicts the assumption that $\mathcal{E}$ has rank $r$.

By the claim above, we have two cases:
$(-K_N.\varphi (l),\det {\mathcal F}.\varphi (l))=(r+1,n+1)$ and 
$(-K_N.\varphi (l),\det {\mathcal F}.\varphi (l))=(r,n)$.
Let $X$ denote $P\times _N l$,
$\varphi_l:X\ra l$ the projection,
and $\pi_X$ the composite of $\pi$ and the projection $X\ra P$.
Let $g:X\ra Y$ be the projective morphism with connected fibers
onto a normal projective variety $Y$ 
determined by
$|H(\mathcal{F}\otimes \mathcal{O}_{l}(-1))|$.
We have the following commutative diagram
\[
\begin{CD}
Y                @<{g}<<X         @>{\varphi_l}>>l                  \\
@.                      @VVV                     @VVV               \\
                   @.   P         @>{\varphi}>>  N                  \\
@.                      @VV{\pi}V                @.                 \\
                   @.   M.         @.            
\end{CD}
\]
Since 
there exists a section $\tilde{l}$ of $\varphi_l$
such that $H(\mathcal{F}\otimes \mathcal{O}_{l}(-1)).\tilde{l}=0$
and $\pi_X(\tilde{l})$ is a point,
we obtain a unique finite morphism $h:Y\ra M$
such that $\pi_X=h\circ g$.

We will show that $M\cong \mathbb{P}^n$.

Suppose that $\det {\mathcal F}.\varphi(l)=n+1$.
Then 
$Y=\mathbb{P}^n$.
Hence we have
$M\cong \mathbb{P}^n$ by Lazarsfeld's theorem
\cite[Theorem~4.1]{l}.

Suppose that $\det {\mathcal F}.\varphi(l)=n$;
denoting by $F$ a fiber of $\pi$,
we infer that 
$\mathcal{F}\otimes \mathcal{O}_F\cong \mathcal{O}(1)^{\oplus n}$
by Theorem~\ref{uniform} (1).
The following argument is inspired by \cite[\S 4 (b.2)]{fn}.
Set $D_N=\varphi (F)$ and $D_P=\varphi^{-1}(D_N)$.
Note here that both $D_N$ and $D_P$ are prime divisors.
Since 
$\mathbb{P}(\mathcal{F}\otimes \mathcal{O}_F)\to F$ has a section
contracted by the composite
$\mathbb{P}(\mathcal{F}\otimes \mathcal{O}_F)\to P\to M$,
we see that $\pi(D_P)$ is also a prime divisor.
This implies that 
$0=D_P.l=\varphi^*D_N.l=D_N.\varphi (l)$.
On the other hand, since $\Pic N\cong \mathbb{Z}$, $D_N$ is ample;
thus $D_N.\varphi(l)>0$. This is a contradiction,
and this case does not occur.

Let $C_0$ be a line on $M\cong \mathbb{P}^n$.
Set $\det \mathcal{E}\cong \mathcal{O}(r+\alpha)$ 
for some integer $\alpha\geq 0$.
Let $C_1$ be a line in a fiber $\mathbb{P}^{n-1}$ of $\varphi$.
Since $L.C_1=1$, $C_1\to \pi(C_1)$ is birational.
Moreover it follows from 
Corollary~\ref{cor for l(R_1)=n}
that we may assume that 
$\pi(C_1)$ is a line in $M\cong \mathbb{P}^n$.
Proposition~\ref{lr1-lr=a} then implies that $\alpha=1$.
Therefore $\mathcal{E}$ is a uniform vector bundle,
and we infer that $\mathcal{E}\cong 
\mathcal{O}(1)^{\oplus (r-1)}\oplus \mathcal{O}(2)$
by Theorem~\ref{uniform}.
This implies that $\varphi$ is birational.
This is a contradiction.
Therefore we conclude that $r=n$.

Finally, if $r=n$, we see that $K_M+\det\mathcal{E}=0$.
Hence we infer that $(M,\mathcal{E})
\cong (\mathbb{P}^n,T_{\mathbb{P}^n})$
by the argument~\cite[\S 4]{fn} 
or \cite[\S 2, Main case 1, Subcase A]{p}.
This completes the proof of the lemma.
\end{proof}

\begin{rmk}
Note that we do not use Theorem~\ref{cm}
in the proof above;
the proof above remains the original argument in \cite{o}.
If we apply Theorem~\ref{cm},
we can give a much shorter proof of Lemma~\ref{key} 
by the similar argument as 
in the proof of Proposition~\ref{lr1=n-1&n-2fiber}.
\end{rmk}
Applying Lemma~\ref{key},
we obtain the following theorem.
\begin{thm}\label{l(r1)=n}
Let $M$ be a Fano manifold of Picard number one
and $\mathcal{E}$ an ample vector bundle of rank $r$ on $M$.
Denote by $\tau$ the nef value 
$\tau (M,\det\mathcal{E})$
of the polarized manifold
$(M,\det\mathcal{E})$.
Suppose that $\tau \geq 1$.
Let $R$, $L$, $R_1$, $C_1$, and 
$\varphi:P\to N$
be as in \S~\ref{the common setup},
and suppose that $L.C_1=1$.
Suppose that 
the length $l(R_1)$ of $R_1$ is equal to $n$.
Then we have one of the following:
\begin{enumerate}
\item $(M,\mathcal{E})\cong (\mathbb{Q}^n,
      \mathcal{O}(1)^{\oplus r})$
      and 
      $\tau r=n$;
\item $(M,\mathcal{E})\cong (\mathbb{P}^n,
      \mathcal{O}(1)^{\oplus (r-1)}\oplus \mathcal{O}(2))$;
\item $(M,\mathcal{E})\cong (\mathbb{P}^n,
      T_{\mathbb{P}^n})$
      and 
      $\tau r=n$.
\end{enumerate}
\end{thm}
\begin{proof}
Suppose that $\varphi$ has an $n$-dimensional fiber $\varphi^{-1}(z)$.
Then Lemmas~\ref{comparisonshin} (1)
and \ref{C_1'=C_1} implies that $\pi(C_1)$
is a minimal extremal rational curve in $R$.
Therefore we have $l(R)=l(R_1)+
(\det\mathcal{E}.\pi(C_1)-rL.C_1)$ 
by Lemma~\ref{lr1-lr=a}.
Let $W$ be the normalization of an $n$-dimensional
irreducible component of $\varphi^{-1}(z)$.
Then $W\cong M$ via $\pi$ by Lemma~\ref{isom}.

Suppose moreover that a general fiber $F$ of $\varphi$ is $n$-dimensional.
By taking $\varphi^{-1}(z)$ as a general fiber,
we may assume that $F=W$.
Since $l(R_1)=n$ and $L.C_1=1$,
we have $(K_P+nL).C_1=0$; thus $K_F+nL|_F=0$.
Hence $(F,L|_F)\cong (\mathbb{Q}^n,\mathcal{O}(1))$ by 
Theorem~\ref{ko} (2).
Therefore we have $\mathbb{Q}^n=F=W\cong M$ via $\pi$;
thus $\pi(C_1)$ is a line in $M\cong \mathbb{Q}^n$.
Since $l(R)=n$ and $l(R_1)=n$,
we have $\det\mathcal{E}.\pi(C_1)-rL.C_1=0$.
Thus we have $\det\mathcal{E}.\pi(C_1)=r$
since $L.C_1=1$.
Hence $\mathcal{E}\cong \mathcal{O}(1)^{\oplus r}$
by \cite[Lemma 3.6.1]{w3}.
This is the case (1) of the theorem.

Suppose moreover that
a general fiber of $\varphi$ has dimension $<n$.
We have $M\cong \mathbb{P}^n$ by Lemma~\ref{standard} (2);
thus $\pi(C_1)$ is a line in $M\cong \mathbb{P}^n$.
Since $l(R)=n+1$,
we have $\det\mathcal{E}.\pi(C_1)-rL.C_1=1$.
The assumption $L.C_1=1$ then
implies that $\det\mathcal{E}.\pi(C_1)=r+1$.
Note here that $n=l(R_1)\geq \tau r\geq r$
by (\ref{lb of r1=}) in 
\S~\ref{the common setup},
since $\tau \geq 1$ and $L.C_1=1$.
Now that $\varphi$ has an $n$-dimensional fiber,
this implies that  $\mathcal{E}\cong 
\mathcal{O}(1)^{\oplus (r-1)}\oplus \mathcal{O}(2)$
by Theorem~\ref{uniform}.
This is the case (2) of the theorem.

Suppose that $\varphi$ has no $n$-dimensional fibers.
Then it follows from inequality~(\ref{ub of r1})
that $\varphi$ is of fiber type
and that every fiber of $\varphi$ is $(n-1)$-dimensional.
Since $(K_P+nL).C_1=0$,
we have $(F,L|_F)\cong 
(\mathbb{P}^{n-1},\mathcal{O}(1))$ for a general fiber
of $\varphi$.
Moreover Lemma~\ref{fujita lemma} implies 
that $\varphi$ makes $(P,L)$ a scroll
over an $r$-dimensional manifold $N$.
Now it follows from Lemma~\ref{key} 
that $(M,\mathcal{E})
\cong (\mathbb{P}^n,T_{\mathbb{P}^n})$.
This is the case (3) of the theorem.
\end{proof}
Finally we pose the following conjecture.
\begin{conj}
Let $M$ be a Fano manifold of Picard number one
and $\mathcal{E}$ an ample vector bundle of rank $r$ on $M$.
Let $L$, $R_1$, and $C_1$ be as in \S~\ref{the common setup},
and set $L.C_1=l$.
If $l(R_1)=n$, then
$(M,\mathcal{E})$ is one of the following:
\begin{enumerate}
\item $(\mathbb{Q}^n, \mathcal{O}(l)^{\oplus r})$;
\item $(\mathbb{P}^n, 
      \mathcal{O}(l)^{\oplus (r-1)}\oplus \mathcal{O}(l+1))$;
\item $(\mathbb{P}^n, T_{\mathbb{P}^n}(l-1)\oplus
      \mathcal{O}(l)^{\oplus (r-n)})$ $(r\geq n)$.
\end{enumerate}
\end{conj}

\section{Some modifications}\label{shin prelim4}
The following is a slight modification of 
\cite[Proposition~3.5]{psw}.

\begin{prop}\label{pswa}
Let $\mathcal{E}$ be an ample vector bundle of rank $r$
on a projective manifold $M$ of dimension $n$. 
Suppose that $r\leq n-1$.
\begin{enumerate}
\item
   If $M=\mathbb{P}^n$ and $\det \mathcal{E}=\mathcal{O}(r+2)$,
   then one of the following holds:
   \begin{enumerate}
     \item
          $\mathcal{E}\cong
          \mathcal{O}(1)^{\oplus (r-1)}\oplus \mathcal{O}(3)$;
     \item
          $r\geq 2$ and 
          $\mathcal{E}\cong
          \mathcal{O}(1)^{\oplus (r-2)}\oplus \mathcal{O}(2)^{\oplus 2}$;
     \item
          $n=3$ and $\mathcal{E}\cong
          N(2)$ where $N$ is a null-correlation bundle.
   \end{enumerate}
\item
   If $M=\mathbb{Q}^n$ and $\det \mathcal{E}=\mathcal{O}(r+1)$,
   then one of the following holds:
   \begin{enumerate}
     \item 
        $\mathcal{E}\cong
        \mathcal{O}(1)^{\oplus (r-1)}\oplus \mathcal{O}(2)$;
     \item
        $n=3$, $r=2$ and 
        $\mathcal{E}\cong
        \mathbf{E}(2)$ where $\mathbf{E}$ is a spinor bundle over 
        $\mathbb{Q}^3$;
     \item
        $n=4$ and $\mathcal{E}\cong
        \mathbf{E}(2)\oplus \mathcal{O}(1)^{\oplus i}$ $(i=0,1)$
        where $\mathbf{E}$ is a spinor bundle of rank $2$ over
        $\mathbb{Q}^4$.
    \end{enumerate}
\end{enumerate}
\end{prop}
\begin{proof}
We owe to \cite{ma1} the idea of the following proof.
Let $\mathcal{E}'$ be 
$\mathcal{E}\oplus \mathcal{O}(1)^{\oplus (n-1-r)}$.
After applying \cite[Proposition~3.5]{psw} to $\mathcal{E}'$,
we recover $\mathcal{E}$ from $\mathcal{E}'$
and we obtain the proposition.
\end{proof}

Next, we will improve Proposition~(1.1) 
in \cite{psw}.
\begin{prop}\label{pswprop}
Let $P$ be a proper normal variety and 
$L$ a line bundle on $P$.
Suppose that there exist coprime positive integers
$p$ and $q$
such that $L^{\otimes p}$ and $L^{\otimes q}$ are 
spanned, and assume that
the image of the map 
$\Phi_{|L^{\otimes p}|}:P\ra \mathbb{P}^N$
is at most of dimension $k$
and that $\dim H^i(P,L^{\otimes t})=0$
for $i>0$ and $t\geq -k+1$.
Then $L$ is spanned and defines a map
$\Phi_{|L|}:P\ra Y\subseteq \mathbb{P}^N$ with connected fibers
onto a normal polarized variety $(Y,\mathcal{M})$,
i.e., $\Phi_{|L|}^*\mathcal{M}=L$,
of Fujita's $\Delta$-genus $\Delta (Y,\mathcal{M})=0$.
Moreover $Y$ has dimension either $k-1$ or $k$, and 
\begin{enumerate}
\item if $\dim Y=k$ then the Hilbert polynomial of 
      $(P,L)$ is of the form 
      \[\frac{1}{k!}\{\prod_{i=1}^{k-1}(t+i)\}(dt+k),\]
      where $d=\mathcal{M}^k$;
\item if $\dim Y=k-1$ then $(Y,\mathcal{M})\cong 
      (\mathbb{P}^{k-1},\mathcal{O}(1))$.
\end{enumerate}
\end{prop}
\begin{proof}
Put $Y=\proj \oplus_{t\geq 0}H^0(L^{\otimes t})$.
Since $L$ is semi-ample,
$L$ induces a natural morphism 
$\phi: P=\proj \oplus_{t\geq 0}L^{\otimes t}\ra
\proj \oplus_{t\geq 0}H^0(L^{\otimes t})$.
Moreover $Y$ is a normal projective variety
and $\phi$ has connected fibers 
because $P$ is normal
(see, for example, \cite[Prop~1.4]{mo3}).
Since $L^{\otimes p}$ 
is spanned, $\mathcal{O}_Y(p)$ is an ample line bundle
on $Y$ and we have $L^{\otimes p}=\phi^*\mathcal{O}_Y(p)$.
We also have $L^{\otimes q}=\phi^*\mathcal{O}_Y(q)$
because $L^{\otimes q}$
is also spanned. 
(Note that $\oplus_{t\geq 0}H^0(L^{\otimes t})$
is not a priori generated by 
elements in $H^0(L)$ of degree one, and thus
we cannot assume at first that $\mathcal{O}_Y(1)$ 
is invertible.)
Let $a$ and $b$ be integers
such that $ap+bq=1$, and 
let $\mathcal{M}$ denote the line bundle 
$\mathcal{O}_Y(ap)\otimes\mathcal{O}_Y(bq)$.
Then $L=\phi^*\mathcal{M}$. 
Since $\phi_*\mathcal{O}_P=\mathcal{O}_Y$,
we have $\mathcal{O}_Y(p)=\mathcal{M}^{\otimes p}$.
Hence we infer that $\mathcal{M}$ is ample. 

Let us consider the Hilbert polynomial $\chi(t)
:=\chi(P,tL)$.
By assumption we see that $Y$ is of dimension at most $k$.
Hence $c_1(L)^{k+1}$ is numerical trivial,
and therefore the degree of $\chi(t)$ is at most $k$.

Suppose that $\deg \chi(t)=k$.
Then, from the vanishing of $\chi(t)$
for $t=-k+1,\cdots,-1$ and from $\chi(0)=1$ we have
\[\chi(t)=\frac{d}{k!}\{\prod_{i=1}^{k-1}(t+i)\}(t+\frac{k}{d})\]
for some positive integer $d$.
Also by the vanishing of the higher cohomology,
we see that $\dim H^0(P,L^{\otimes t})=\chi(t)$ for $t\geq 0$.
Since $H^0(P,L^{\otimes t})\cong H^0(Y,\mathcal{M}^{\otimes t})$
for all $t$, $\chi(t)$ is also the 
Hilbert polynomial of $(Y,\mathcal{M})$; $\dim Y=k$
and $d=\mathcal{M}^k$.
Now we have $h^0(Y, \mathcal{M})=h^0(P, L)=d+k$.
Hence $\Delta (Y,\mathcal{M})=0$.

Suppose that $\deg \chi(t)<k$.
Then, from the vanishing of $\chi(t)$
for $t=-k+1,\cdots,-1$ and from $\chi(0)=1$ we have
\[\chi(t)=\frac{1}{(k-1)!}\prod_{i=1}^{k-1}(t+i).\]
In particular, we see that $\deg \chi(t)=k-1$.
For the same reason as above, 
$\chi(t)$ is also the Hilbert polynomial of $(Y,\mathcal{M})$; 
in this case we have $\dim Y=k-1$ and $1=\mathcal{M}^{k-1}$.
Now we have $h^0(Y, \mathcal{M})=h^0(P, L)=k$.
Hence $\Delta (Y,\mathcal{M})=0$.
Therefore $(Y,\mathcal{M})\cong (\mathbb{P}^{k-1},\mathcal{O}(1))$. 

In both cases, we see that $\mathcal{M}$ is very ample.
Hence $L$ is spanned and the statement follows.
\end{proof}

The following Corollary~\ref{pswcor} (1)
is nothing but Corollary~(1.2) in \cite{psw}.
Corollary~\ref{pswcor} (2)
is a modification of Corollary~(1.3) in \cite{psw}.
Corollary~\ref{pswcor} (2) will be used 
in the proof of Proposition~\ref{no n-1 fiber}
and in \S~\ref{-K=n-2shin}.
Corollary~\ref{pswcor} (3) will 
not be used in this paper.
\begin{cor}\label{pswcor}
Let $M$ be an $n$-dimensional projective manifold,
and $\mathcal{E}$ a rank-$r$ vector bundle on $M$.
Assume that $H^i(M,\mathcal{O}_M)=0$ for $i>0$,
that $H(\mathcal{E})^{\otimes m}$ is a spanned 
line bundle on $\mathbb{P}(\mathcal{E})$ for 
\textsf{all} $m\gg 0$,
and that $H^i(\mathbb{P}(\mathcal{E}),H(\mathcal{E})^{\otimes t})=0$
for $i>0$, $t>0$.
\begin{enumerate}
\item If $c_1(H(\mathcal{E}))^r=0$, then 
      $\mathcal{E}\cong \mathcal{O}^{\oplus r}$.
\item If $c_1(H(\mathcal{E}))^{r+1}=0$,
      then $\mathcal{E}$ is spanned.
      Suppose moreover that 
      $h^i(M,K_M\otimes \det \mathcal{E})=0$ for $i<n$.
      Then $h^n(M,K_M\otimes \det \mathcal{E})=1$ or $0$
      and $h^n(M,K_M\otimes \det \mathcal{E})=1$ if and only if 
      $c_1(H(\mathcal{E}))^r=0$. Furthermore in case 
      $h^n(M,K_M\otimes \det \mathcal{E})=0$,
      we have $h^0(\mathcal{E})=r+1$, and thus
      $\mathcal{E}$ fits into the following exact sequence
      \[0\to \mathcal{E}^*\to \mathcal{O}_M\otimes 
      H^0(\mathcal{E})^*\to \det\mathcal{E}\to 0.\]
\item If $c_1(H(\mathcal{E}))^{r+2}=0$
      and $h^i(M,K_M\otimes \det \mathcal{E})=0$ for all $i$,
      then $\mathcal{E}$ is spanned.
\end{enumerate}
\end{cor}
\begin{proof}
Put $P=\mathbb{P}(\mathcal{E})$ and $L=H(\mathcal{E})$.

Suppose that we are in (1) or (2).
In these cases, we set $k=r$.
The vanishing of $H^i(P,L^{\otimes t})$ for $-k+1\leq t\leq -1$ and $i>0$
follows from Leray's spectral sequence.
Since $L^{\otimes m}$ is spanned for all $m\gg 0$
and $L^{r+1}=0$ in both cases (1) and (2),
we now apply Proposition~\ref{pswprop}
to see that $L$ is spanned and 
that the image $Y$ of the morphism
determined by $L$ has dimension $k-1$ or $k$.
Therefore $\mathcal{E}$ is spanned.

(1) If $L^r=0$, then $\dim Y<k$, and 
again by Proposition~\ref{pswprop},
we infer that $h^0(L)=k$ and thus $h^0(\mathcal{E})=r$.
Therefore $\mathcal{E}$ is trivial.

(2) Let $\chi (t)$ denote the Hilbert polynomial $\chi (P, tL)$ 
of $(P,L)$.
Suppose that $h^i(M,K_M\otimes \det \mathcal{E})=0$ for $i<n$.
By Serre duality, we have
\[
\begin{split}
h^i(P,L^{\otimes (-k)})&=h^{n+r-1-i}(P,K_P\otimes L^{\otimes k})
=h^{n+r-1-i}
(M,K_M\otimes \det \mathcal{E})\\
&=
\begin{cases}
0 & \text{ if $i\neq r-1$},\\
h^n(M,K_M\otimes \det \mathcal{E}) & \text{ if $i=r-1$}.
\end{cases}
\end{split}
\]
Hence $\chi (-k)=(-1)^{r-1}
h^n(M,K_M\otimes \det \mathcal{E})$.
On the other hand, we have,   
by Proposition~\ref{pswprop} and its proof,
\[\chi (-k)=
\begin{cases}
(-1)^{k-1}(-d+1) & \text{ if $\dim Y=k$},\\
(-1)^{k-1} & \text{ if $\dim Y=k-1$}.
\end{cases}
\]
where $d=\deg Y$.
Therefore $h^n(M,K_M\otimes \det \mathcal{E})=1$
if and only if $\dim Y=k-1$, i.e., $L^k=0$.
Moreover if $h^n(M,K_M\otimes \det \mathcal{E})=0$
then $d=1$ and $\dim Y=k$, and thus
Proposition~\ref{pswprop} implies  
$(Y,\mathcal{M})\cong (\mathbb{P}^k,\mathcal{O}(1))$.
Hence $k+1=h^0(\mathcal{M})=h^0(L)=h^0(\mathcal{E})$.

Suppose that we are in (3).
In this case we set $k=r+1$.
The vanishing of $H^i(P,L^{\otimes t})$ for 
$-k+2\leq t\leq -1$ and $i>0$
follows from Leray's spectral sequence.
By Serre duality and the assumption, we have
\[
\begin{split}
h^i(P,L^{\otimes (-k+1)})&=h^{n+r-1-i}(P,K_P\otimes L^{\otimes k-1})\\
&=h^{n+r-1-i}(M,K_M\otimes \det \mathcal{E})=0.
\end{split}
\]
Since $L^{\otimes m}$ is spanned for all $m\gg 0$
and $L^{r+2}=0$,
we apply Proposition~\ref{pswprop}
to see that $L$ is spanned.
Therefore $\mathcal{E}$ is spanned.
This completes the proof.
\end{proof}

\section{Lemmas for the case 
$\varphi$ is of fiber type}\label{Sl(r1)=s'}
\begin{lemma}\label{l(r1)=s'}
Let $M$ be an $n$-dimensional
Fano manifold of Picard number one
and $\mathcal{E}$ an ample vector bundle of rank $r$ on $M$.
Let $\tau$, $\pi:P\to M$, $L$, $R_1$, $C_1$, and 
$\varphi:P\to N$
be as in \S~\ref{the common setup}.
In particular assume that $-K_P$ is ample.
Suppose that a general fiber of $\varphi$ is $n$-dimensional.
Then $\mathcal{E}\cong \mathcal{O}_M(D)^{\oplus r}$,
where $\mathcal{O}_M(D)$ is an ample line bundle
on $M$.
Moreover if $L.C_1=1$ then
$\mathcal{O}_M(D)$ is the ample generator
$\mathcal{O}_M(1)$ of $\Pic M$,
the Fano index of $M$ is $l(R_1)$,
and $l(R_1)=\tau r$. 
\end{lemma}
\begin{proof}
Let $F$ be a general fiber $F$ of $\varphi$.
It follows from Lemma~\ref{isom} that 
$F\cong M$ via $\pi$.
Let $\mathcal{O}_M(D)$ be the line bundle 
corresponding to $L|_F$ via $F\cong M$.
For a curve $l$ on $M$,
letting $\tilde{l}$ be a curve in $F$ corresponding to $l$,
we see that $H(\mathcal{E}(-D)).\tilde{l}=0$.
Hence $H(\mathcal{E}(-D))$ is nef
and is a supporting divisor for $R_1$ (or $\varphi$).
This implies that there exists a positive integer $m_0$
such that $H(\mathcal{E}(-D))^{\otimes m}$ is spanned for 
all integers $m\geq m_0$;
by Corollary~\ref{pswcor} (1), we infer that 
$\mathcal{E}(-D)\cong \mathcal{O}^{\oplus r}$.
Hence $\mathcal{E}\cong \mathcal{O}_M(D)^{\oplus r}$.
Now we may think that $C_1$ is contained in $F$.
If $L.C_1=1$,
then $\mathcal{O}_M(D)$
is the ample generator $\mathcal{O}_M(1)$ of $\Pic M$.
Moreover if $L.C_1=1$,
then $K_F+l(R_1)L|_F=0$, and thus $K_M+l(R_1)\mathcal{O}_M(1)=0$.
Hence the Fano index of $M$ is $l(R_1)$.
Note finally that this also implies that $l(R_1)=\tau r$
since $\mathcal{E}\cong \mathcal{O}_M(1)^{\oplus r}$.
\end{proof}

The idea of the following lemma comes from \cite[\S 4 (b.2)]{fn}.
\begin{lemma}\label{kantanlemma}
Let $M$ be an $n$-dimensional 
Fano manifold of Picard number one
and $\mathcal{E}$ an ample vector bundle of rank $r$ on $M$.
Let $\pi:P\to M$, $L$, and $\varphi:P\to N$
be as in \S~\ref{the common setup}.
Suppose that $\varphi$ is of fiber type.
If a general fiber $W$ of $\varphi$ has dimension $<n$,
then $\mathcal{E}_W$ cannot be isomorphic to $L_W^{\oplus r}$.
\end{lemma}
\begin{proof}
Suppose that $\mathcal{E}_W$ is isomorphic to 
$L_W^{\oplus r}$.
Note that $\mathbb{P}(\mathcal{E}_W)\to W$
has a section $W_s$ 
induced from the inclusion $W\hookrightarrow P$.
Set $P_W=\mathbb{P}(\mathcal{E}_W)$.
Since $W_s$ is contracted by 
$P_W\hookrightarrow P\to N$,
it follows from $\mathcal{E}_W\cong L_W^{\oplus r}$
that $P_W$ is contracted to a variety of dimension $r-1$.
Fix one $W$, and denote it by $W_0$.
Set $M_0=\pi (W_0)$ and $P_0=\pi^{-1}(M_0)$. 
We have $\dim \varphi(P_0)=r-1$.

Set $t=\dim M_0$.
Note that $t=\dim W_0$
and thus  $\dim N=n+r-1-t$.
Hence we have $\codim (M_0,M)=n-t=\codim (\varphi(P_0),N)$.
Since $\dim W_0<n$ by assumption,
we have $t<n$. Therefore 
$\varphi(P_0)$ is a proper subset of $N$.

If $t=n-1$,
then $M_0$ is an ample divisor since $\Pic M\cong \mathbb{Z}$.
Hence we have $0<\pi^*M_0.e=P_0.e$
for any curve $e$ in a fiber of $\varphi$.
On the other hand, 
we have $P_0.e=0$ if $e$ lies in a fiber
over $N\setminus \varphi(P_0)$.
This is a contradiction.

Suppose that $t\leq n-2$.
Set $N_0=\varphi(P_0)$.
Take a general curve $Y_1$ such that
$Y_1$ is not contained in $N_0$.
Denote by $W_1$ the irreducible component
of $\varphi^{-1}(Y_1)$ which dominates $Y_1$.
Set $M_1=\pi(W_1)$. Then $\dim M_1=\dim W_0+1$.
Set $P_1=\pi^{-1}(M_1)$ and $N_1=\varphi(P_1)$.
Then a general fiber of $P_1\to N_1$ is $t$-dimensional.
Hence $\codim (N_1,N)=\codim (P_1,P)=\codim (M_1,M)=n-t-1$.
Next take a general surface $Y_2$
such that $Y_2$ is not contained in $N_1$.
Denote by $W_2$ the irreducible component
of $\varphi^{-1}(Y_2)$ which dominates $Y_2$.
Set $M_2=\pi(W_2)$. Then $\dim M_2=\dim W_0+2$.
Set $P_2=\pi^{-1}(M_2)$ and $N_2=\varphi(P_2)$.
Then a general fiber of $P_2\to N_2$ is $t$-dimensional.
Hence $\codim (N_2,N)=\codim (P_2,P)=\codim (M_2,M)=n-t-2$.
Likewise we can repeat this procedure inductively,
so that if we take a general $(n-t-1)$-dimensional variety $Y_{n-t-1}$,
letting $W_{n-t-1}$ be the irreducible component
of $\varphi^{-1}(Y_{n-t-1})$ which dominates $Y_{n-t-1}$,
and setting $M_{n-t-1}=\pi(W_{n-t-1})$, 
$P_{n-t-1}=\pi^{-1}(M_{n-t-1})$, and $N_{n-t-1}=\varphi(P_{n-t-1})$,
we see that  $\codim (M_{n-t-1},M)=1=\codim (N_{n-t-1},N)$.
Since $M_{n-t-1}$ is ample,
we have $0<\pi^*M_{n-t-1}.e=P_{n-t-1}.e=0$
for a curve $e$ in a fiber over $N\setminus N_{n-t-1}$.
This is a contradiction.
\end{proof}

\section{Comparison Lemma for $l(R_1)=n-1$}\label{Scompaforlr1n-1}
In this section, we will follow the notation
in \S~\ref{section comparison}.
\begin{cor}\label{comparison for lr1=n-1}
Suppose that $l(R_1)=n-1$,
that $L.C_1=1$,
that $\tau \geq 1$,
and that $\varphi$ has no 
$n$-dimensional fibers.
Then $\pi(C_1)$ is a minimal extremal rational curve in $R$
unless the following holds:\\
(*)
$n=4$, $r=3$, $\tau =1$, 
$\varphi$ is of fiber type,
every fiber of $\varphi$ has a
two dimensional irreducible component $F$
such that $\dim p_x(V_x)\cap F=0$
(some special fiber of $\varphi$ might have
a $3$-dimensional irreducible component $F'$
such that $p_x(V_x)\cap F'=\emptyset$).
\end{cor}
\begin{proof}
Note first that $r\geq 2$
since $\varphi$ has no $n$-dimensional fibers.

By Lemma~\ref{C_1'=C_1},
it is enough to show that 
there exists an unsplit rational curve $C_1'$
belonging to $R_1$
and dominating a minimal extremal rational curve $C_0$.
Assume, to the contrary,
that no unsplit rational curve
dominating a minimal extremal rational curve in $R$
belongs to $R_1$,
i.e., is contracted by $\varphi$.
Let $C_1'$ be an unsplit 
rational curve
dominating a minimal extremal rational curve in $R$.
Then, by Lemma~\ref{comparisonshin} (2),
we may assume that $p_x(V_x)\to \varphi(p_x(V_x))$
is finite.
Hence we have
\[
\dim N\geq 
\dim \varphi(p_x(V_x))=
\dim p_x(V_x)=\dim V_x\geq \tau rL.C_1'+r-2
\]
by Lemma~\ref{comparisonshin} (4) and (5),
since $\tau\geq 1$.

Suppose that $\varphi$ is of fiber type.
Since $l(R_1)=n-1$, 
inequality (\ref{ub of r1}) in \S~\ref{the common setup}
implies that every fiber of $\varphi$
has dimension $\geq n-2$.
Hence we have $r+1\geq \dim N$,
and thus $3\geq \tau rL.C_1'\geq 2$.
Therefore $L.C_1'=1$, $3\geq \tau r\geq 2$, and $r=2$ or $3$.

Suppose that $r=2$ and that $\varphi$ is of fiber type.
If $\tau =1$, we have $n-1=l(R_1)=rL.C_1=2$
by (\ref{lb of r1=}) in \S~\ref{the common setup},
since $L.C_1=1$. Hence $n=3$.
In this case, we have $\pi(C_1)=C_0$
by \cite{w2}, 
and thus there exists 
an unsplit rational curve $C_1$
dominating $C_0$
and contracted by $\varphi$,
which contradicts the assumption.
Hence $\tau >1$.
Then we have $\tau r>2$,
and thus $r+1=\dim N=\dim \varphi(p_x(V_x))$.
Hence a general fiber of $\varphi$ has dimension $n-2$.
Let $F$ be a general fiber of $\varphi$.
Since $L.C_1=1$, we have $(K_P+(n-1)L).C_1=0$,
and thus $K_F+(n-1)L_F=0$.
Theorem~\ref{ko} (1) then implies 
that $(F,L_F)\cong (\mathbb{P}^{n-2},\mathcal{O}(1))$.
Note here that $n-1=l(R_1)\geq \tau r>2$
by (\ref{lb of r1=}) in \S~\ref{the common setup},
and thus $n\geq 4$ and $\dim F\geq 2$.
Since $L_F=\mathcal{O}(1)$,
we have the following exact sequence
\[0\to \mathcal{O}(\alpha)\to \mathcal{E}_F(-1)\to \mathcal{O}\to 0\]
for some integer $\alpha$.
Since $F$ is a projective space of dimension $\geq 2$,
this exact sequence must split: 
$\mathcal{E}_F(-1)\cong \mathcal{O}(\alpha)\oplus\mathcal{O}$.
Since $\mathcal{E}_F$ is ample, we have $\alpha\geq 0$.
Moreover Lemma~\ref{kantanlemma} implies that $\alpha>0$.
Hence the restriction of $\varphi$ to $\pi^{-1}(\pi(F))$
is the composite of a birational morphism 
and a finite morphism, 
and hence $r+1=\dim N\geq \dim \pi^{-1}(\pi(F))
=n-1$. Since $n\geq 4$, this implies that $n=4$.
Set $\det\mathcal{E}.C_0=r+a$
for a minimal extremal rational curve $C_0$ in $R$.
Then inequality 
(\ref{bound of detE.C_0}) in \S~\ref{overview}
shows that $0\leq a\leq l(R)-\tau r<5-2$
since $\tau>1$. 
Hence $0\leq a\leq 2$.
Moreover if $a=2$ then $l(R)=n+1$
and thus $M\cong \mathbb{P}^4$
by Theorem~\ref{cm}.
Then Proposition~\ref{pswa} implies
that $\mathcal{E}$ is either 
$\mathcal{O}(2)^{\oplus 2}$ or
$\mathcal{O}(1)\oplus \mathcal{O}(3)$.
This, however, contradicts the fact
that $\varphi$ has no $n$-dimensional fibers.
Hence $a\leq 1$.
Suppose that $a=1$.
If $l(R)=n+1$, then we see that
$(M,\mathcal{E})\cong (\mathbb{P}^4,
\mathcal{O}(1)\oplus \mathcal{O}(2))$
by Theorem~\ref{cm},
which contradicts that $\varphi$ has no $n$-dimensional fibers.
Hence $l(R)=4$.
Putting $A=-K_M-\det \mathcal{E}$, 
we see that $A.C_0=1$.
Thus $A$ is ample
and $K_M+4A=0$ 
because the Picard number of $M$ is one.
Hence $(M,A)\cong (\mathbb{Q}^4,\mathcal{O}(1))$ 
by Theorem~\ref{ko} (2).
Now we apply Proposition~\ref{pswa}
to know that $\mathcal{E}\cong 
\mathbf{E}(2)$ where $\mathbf{E}$ is a spinor bundle
because $\varphi$ is of fiber type.
In this case, $\varphi$ is a $\mathbb{P}^2$-bundle
and every fiber of $\varphi$ is mapped to a 
plane in $\mathbb{Q}^4$ by $\pi$.
Hence there exists an unsplit rational curve $C_1$
dominating $C_0$ and belonging to $R_1$,
which contradicts the assumption.
Suppose that $a=0$.
Then $l(R)=\tau \det\mathcal{E}.C_0=\tau r$.
Since $3\geq \tau r>2$, we have $\tau r=3$ and $\tau =3/2$.
Now it follows from (\ref{-K_P.C_1alpha siki})
in \S~\ref{the common setup}
that $\det\mathcal{E}.\pi_*(C_1)=2=r$,
which implies that $\pi(C_1)$
is a minimal extremal rational curve in $R$.
Therefore there exists an unsplit rational curve $C_1$
belonging to $R_1$ and dominating 
a minimal extremal rational curve in $R$.
This contradicts the assumption.

Suppose that $r=3$ and that $\varphi$ is of fiber type.
Since $3\geq \tau r$, we have $\tau =1$. 
Then (\ref{lb of r1=}) in \S~\ref{the common setup}
implies that $n-1=l(R_1)=rL.C_1=3$.
Hence $n=4$.
Since $\tau r>2$, we have 
$r+1=\dim N=\dim \varphi(p_x(V_x))$ as above,
and thus every fiber of $\varphi$ has 
an irreducible component $F$ of dimension $n-2$
such that $\dim p_x(V_x)\cap F=0$.
This is the case (*) of the corollary.

Suppose that $\varphi$ is birational.
Since $l(R_1)=n-1$ and $\varphi$ has no
$n$-dimensional fibers,
we see,
by inequality (\ref{ub of r1}) in \S~\ref{the common setup},
that $\varphi$ is divisorial
and that every positive dimensional fiber has 
dimension $n-1$.
Denote by $E$ the exceptional divisor of $\varphi$.
Theorem~\ref{Andreatta Occhetta}
then implies that $N$ is smooth
and $\varphi$ is the blowing up
along a submanifold $\varphi(E)$ of $N$.
Therefore we have $K_P=\varphi^*(K_N)+(n-1)E$.

Note here that $\pi(E)=M$
since every non-zero effective divisor on $M$
is ample.
Therefore a general fiber $\pi|_E$
has dimension $r-2$.

Let $l$ be a general line in a general fiber of $\pi$;
note that $l$ is not contained in $E$.
We claim here that $E$ intersects $l$ in 
only one point.
The idea of this claim is due to \cite[Lemma (7.2)]{psw}.
Since $E$ defines an ample divisor
of a general fiber of $\pi$,
we see first that 
$E\cap l$ is a non-empty finite set.
Second we show that
$-K_N.\varphi (l)\leq \dim N+1$;
suppose, to the contrary, 
that $-K_N.\varphi (l)\geq \dim N+2$.
Then for general two points $p$, $q$ of $\varphi (l)$,
we can deform $\varphi (l)$ to form a $1$-dimensional 
family of curves all of which pass through $p$ and $q$.
Since $\varphi$ is birational, this family can be lifted up
to a family on $P$, which implies that $l$ can break.
This is a contradiction.
Hence  $-K_N.\varphi (l)\leq \dim N+1$.
Now we have
\[E.l=\frac{K_P.l-\varphi^*(K_N).l}{n-1}\leq \frac{-r+n+r-1+1}{n-1}
=1+\frac{1}{n-1}.\]
Note here that $n-1=l(R_1)\geq \tau rL.C_1\geq r\geq 2$
by (\ref{lb of r1=}) in \S~\ref{the common setup},
since $\tau \geq 1$.
Hence we have $E.l=1$.
Therefore a general fiber of $\pi|_E$
is a hyperplane of a fiber of $\pi$.

Suppose that $\pi|_E$ has an $(r-1)$-dimensional fiber
$\pi^{-1}(x_1)$.
Then it dominates $\varphi(E)$.
Since $\varphi(E)$ is smooth,
Lazarsfeld's theorem~\cite[Theorem~4.1]{l}
implies that $\varphi(E)\cong \mathbb{P}^{r-1}$.
Lemma~\ref{isom} then implies that
the $\pi^{-1}(x_1)$ is isomorphic to $\varphi(E)$
via $\varphi|_E$.
Denote by $\mathcal{F}$ 
the vector bundle $(\varphi|_E)_*(L|_E)$ of rank $n$
on $\varphi(E)$.
We have $E\cong \mathbb{P}(\mathcal{F})$ over 
$\varphi(E)$ and $L|_E\cong H(\mathcal{F})$,
and thus $\mathcal{F}$ is ample.
Let $l$ denote a line in the fiber $\pi^{-1}(x_1)$
by abuse of notation.
We see that $K_E.l=-(r-1)$
and that 
$K_E=-nL|_E+(\varphi|_E)^*(K_{\varphi(E)}+\det\mathcal{F})$.
Since $l$ is mapped isomorphically
onto a line in $\varphi(E)\cong \mathbb{P}^{r-1}$,
we infer that $\det\mathcal{F}.\varphi(l)=n+1$.
Since $\mathcal{F}$ is ample,
$\mathcal{F}$ is a uniform vector bundle
of type $(1,\dots,1,2)$,
and therefore is isomorphic to 
either $\mathcal{O}(1)^{\oplus (n-1)}\oplus \mathcal{O}(2)$
or $T_{\mathbb{P}^{r-1}}\oplus 
\mathcal{O}(1)^{\oplus (n-r+1)}$
by Theorem~\ref{uniform} (2).
Since $\pi|_E$ has connected fibers,
this implies that $M\cong \mathbb{P}^n$
and that $r=2$ 
if 
$\mathcal{F}\cong 
\mathcal{O}(1)^{\oplus (n-1)}\oplus \mathcal{O}(2)$.
Moreover we see that every fiber of $\varphi|_E$
is mapped isomorphically
to a hyperplane of $M\cong \mathbb{P}^n$ via $\pi|_E$.
Therefore we may assume that $\pi(C_1)=C_0$,
which contradicts the assumption.

Suppose that every fiber of $\pi|_E$
is $(r-2)$-dimensional,
so that $\pi|_E$ also makes $(E,L|_E)$
a scroll over $M$.
Now we apply Lemma~\ref{key} to see
that $n=r-1$.
Since $n-1=l(R_1)\geq \tau rL.C_1\geq r$
by (\ref{lb of r1=}) in \S~\ref{the common setup},
this is a contradiction.
Therefore this case does not happen either.
\end{proof}

\section{The case $l(R_1)=n-1$}\label{Slr1=n-1}
Let $M$, $\mathcal{E}$, $\tau$, $\pi:P\to M$,
$L$, $R_1$, $C_1$, and $\varphi:P\to N$
be as in \S~\ref{the common setup}.

\begin{prop}\label{l(r1)=n-1}
Suppose that $\tau \geq 1$
and that $L.C_1=1$.
If the length $l(R_1)$ of $R_1$ is equal to $n-1$
and $\varphi$ has an $n$-dimensional fiber,
then we have one of the following:
\begin{enumerate}
\item 
$K_M+(n-1)A=0$
and $\mathcal{E}\cong A^{\oplus r}$
where $A$ is the ample generator of $\Pic M$.
In this case we have $\tau r=n-1$;
\item 
$(M,\mathcal{E})\cong (\mathbb{P}^n,
\mathcal{O}(1)^{\oplus (r-2)}\oplus \mathcal{O}(2)^{\oplus 2})$
$(r\geq 3)$;
\item
$(M,\mathcal{E})\cong (\mathbb{P}^n,
\mathcal{O}(1)^{\oplus (r-1)}\oplus \mathcal{O}(3))$
$(r\geq 2)$;
\item $(M,\mathcal{E})\cong (\mathbb{Q}^n,
\mathcal{O}(1)^{\oplus (r-1)}\oplus \mathcal{O}(2))$
$(r\geq 2)$;
\item  
$(M,\mathcal{E})\cong (\mathbb{Q}^4,
\mathbf{E}(2)\oplus\mathcal{O}(1))$,
where $\mathbf{E}$ is a spinor bundle on $\mathbb{Q}^4$.
In this case we have $\tau r=n-1$ and $\tau =1$.
\end{enumerate}
\end{prop}
\begin{proof}
Suppose that a general fiber $F$ of $\varphi$ is $n$-dimensional.
Since $l(R_1)=n-1$ and $L.C_1=1$,
it follows from Lemma~\ref{l(r1)=s'}
that $(M,A)$ is a Del Pezzo manifold
and that $\mathcal{E}\cong A^{\oplus r}$,
where $A$ is the ample generator of $\Pic M$.
This is the case $(1)$ of the proposition.

We assume that a general fiber of $\varphi$ has dimension $<n$ in 
the following.

Let $\varphi^{-1}(z)$ be an $n$-dimensional fiber 
of $\varphi$,
and let $W$ be the normalization of an $n$-dimensional
irreducible component of $\varphi^{-1}(z)$.
Then $W\cong M$ via $\pi$ by Lemma~\ref{isom}.
Moreover $(W,L_W)$ is isomorphic to either
$(\mathbb{P}^n,\mathcal{O}(1))$ 
or $(\mathbb{Q}^n,\mathcal{O}(1))$ 
by Lemma~\ref{standard} (1)
and its proof.
Let $\tilde{C}_1$ be a line in $W$;
we have $L.\tilde{C}_1=1$ and the image of $\tilde{C}_1$ in 
$\varphi^{-1}(z)$ is $C_1$.
Since $W\cong M$,
we see that $\tilde{C}_1\to \pi (C_1)$ is an isomorphism
and that $\pi (C_1)$ is a line in $M$;
in particular $C_1\to \pi (C_1)$ is an isomorphism.
Set $\det \mathcal{E}.\pi (C_1)=r+a$ for some integer $a$.
We have $n-1=l(R_1)=-K_M.\pi(C_1)-a$ by Proposition~\ref{lr1-lr=a}.
Thus $a=2$ if $W\cong M\cong \mathbb{P}^n$,
and $a=1$ if $W\cong M\cong \mathbb{Q}^n$.
Note here that $n-1\geq r$ by (\ref{lb of r1=})
in \S~\ref{the common setup},
since $\tau \geq 1$.

Suppose that $M\cong \mathbb{P}^n$
and $\det \mathcal{E}\cong \mathcal{O}(r+2)$.
Since $L.C_1=1$, it follows from 
Proposition~\ref{pswa} (1)
that 
$\mathcal{E}$ is isomorphic to 
$\mathcal{O}(1)^{\oplus (r-1)}\oplus \mathcal{O}(3)$ 
$(r\geq 2)$ or 
$\mathcal{O}(1)^{\oplus (r-2)}\oplus \mathcal{O}(2)^{\oplus 2}$ 
$(r\geq 3)$.
The reason why $\mathcal{E}$ is not isomorphic to $N(2)$,
where $N$ is the null correlation bundle on $\mathbb{P}^3$,
is as follows:
First recall (see, e.g., \cite[I, 4.2]{oss}) an exact sequence
\[0\to \mathcal{O}_{\mathbb{P}^3}\to \Omega(2)\to \check{N}(1)\to 0.\]
Note that $\check{N}(1)\cong N(1)$ on $\mathbb{P}^3$.
Since a natural morphism 
$\mathbb{P}(\Omega(2)) \to \mathbb{P}(H^0(\Omega(2)))$
is a $\mathbb{P}^1$-bundle onto the Grassmannian 
$\mathbb{G}(1,\mathbb{P}^3)$ of lines in $\mathbb{P}^3$,
the exact sequence above with $\check{N}(1)\cong N(1)$
shows that
a natural morphism $\mathbb{P}(N(1)) \to \mathbb{P}(H^0(N(1)))$,
which is $\varphi$ if $\mathcal{E}\cong N(2)$,
is a $\mathbb{P}^1$-bundle onto a three-dimensional hyperquadric.
This contradicts the assumption that $\varphi$ has 
an $n$-dimensional fiber.

Suppose that $M\cong \mathbb{Q}^n$
and that $\det \mathcal{E}\cong \mathcal{O}(r+1)$.
Since $L.C_1=1$, it follows from 
Proposition~\ref{pswa} (2)
that $\mathcal{E}$ is isomorphic to
$\mathcal{O}(1)^{\oplus (r-1)}\oplus \mathcal{O}(2)$
$(r\geq 2)$ or
$\mathbf{E}(2)\oplus \mathcal{O}(1)$
where $\mathbf{E}$ is a spinor bundle over $\mathbb{Q}^4$.
Here, the reason why $\mathcal{E}$ is not isomorphic to $\mathbf{E}(2)$,
where $\mathbf{E}$ is a spinor bundle on 
$\mathbb{Q}^4$ or $\mathbb{Q}^3$,
is as follows:
First recall that
$\mathbb{Q}^4$ is isomorphic to 
the Grassmannian $\mathbb{G}(1,\mathbb{P}^3)$ 
of lines in $\mathbb{P}^3$
and that, by changing this isomorphism if necessary,
$\mathbf{E}(1)$ can be, via this isomorphism, 
isomorphic to the universal quotient bundle over 
$\mathbb{G}(1,\mathbb{P}^3)$
(see \cite[Examples~1.5]{ot});
therefore the natural map
$\mathbb{P}(\mathbf{E}(1)) \to \mathbb{P}(H^0(\mathbf{E}(1)))=\mathbb{P}^3$
is a $\mathbb{P}^2$-bundle,
and every fiber of this map is isomorphic to a plane in $\mathbb{Q}^4$
via the projection.
Second recall that the restriction of a spinor bundle on $\mathbb{Q}^4$
to $\mathbb{Q}^3$ is the spinor bundle on $\mathbb{Q}^3$
by \cite[Theorem~1.4]{ot}; 
hence the natural map 
$\mathbb{P}_{\mathbb{Q}^3}(\mathbf{E}(1)) \to \mathbb{P}^3$
is a $\mathbb{P}^1$-bundle.
Therefore we conclude that $\varphi$ has no $n$-dimensional fibers
if $\mathcal{E}\cong \mathbf{E}(2)$ on $\mathbb{Q}^4$
or $\mathbb{Q}^3$.
This contradicts our assumption.
\end{proof}

In the proof of the following proposition,
applying Theorem~\ref{cm},
we adapt the argument in \cite[\S 5]{psw}
to our case.

\begin{prop}\label{no n-1 fiber}
Suppose that $\tau \geq 1$,
that $L.C_1=1$,
that $l(R_1)=n-1$,
and that $\varphi$ has no $n$-dimensional fibers.
If $\varphi$ has an $(n-1)$-dimensional fiber,
then $n=4$, $r=3$, $\tau =1$, 
$\varphi$ is of fiber type,
every fiber of $\varphi$ has a
two dimensional irreducible component $F$
such that $\dim p_x(V_x)\cap F=0$,
and some special fiber of $\varphi$ has
a $3$-dimensional irreducible component $F'$
such that $p_x(V_x)\cap F'=\emptyset$.
\end{prop}
\begin{proof}
First note that $n-1=l(R_1)\geq r$
by (\ref{lb of r1=}) in \S~\ref{the common setup},
since $\tau \geq 1$.

Since every fiber of $\varphi$ has 
dimension $\leq n-1$,
we see that $r\geq 2$,
and we have two cases:\\
A) a general fiber of $\varphi$ has dimension $n-1$, 
and thus $\dim N=r$;\\
B) a general fiber of $\varphi$ has dimension $<n-1$,
and thus $\dim N>r$.\\
In the case A), 
since $(K_P+(n-1)L)_W=\mathcal{O}_W$
for a general fiber $W$ of $\varphi$,
we have
$(W,f^*(L))\cong (\mathbb{Q}^{n-1},\mathcal{O}(1))$
by Theorem~\ref{ko} (2),
where $f:W\to P$ is the inclusion.
In the case B), 
for any $(n-1)$-dimensional irreducible component
$A_z$ of an $(n-1)$-dimensional fiber $\varphi^{-1}(z)$,
we have
$(W,{f}^*(L))\cong (\mathbb{P}^{n-1},\mathcal{O}(1))$
by Lemma~\ref{Fujita Ye-Zhang},
where 
$f:W\to P$ is the composite of the normalization 
$W\to A_z$ and the inclusion $A_z\to P$.

Set $\mathcal{E}_W=(\pi\circ f)^*\mathcal{E}$
and $P_W=\mathbb{P}(\mathcal{E}_W)$.
Let $\rho : P_W\ra P$
be the morphism induced from $\pi \circ f:W\ra M$
by base change $\pi : P\ra M$.
Set $\varphi_W=\varphi\circ\rho$.
We have the following commutative diagram
\[
\begin{CD}
P_W@>{\rho}>>P@>{\varphi}>>N\\
@V{\pi_W}VV@VV{\pi }V  @. \\
W@>{\pi\circ f}>>M.@.
\end{CD}
\]
Let $W_s$ be the section of $\pi_W:P_W\to W$,
the section induced from $f:W\to P$.
Since $\pi\circ f$
is finite by Lemma~\ref{finite of yz},
we see that $\rho $ is finite
and that $\mathcal{E}_W$ is ample.

If $n=3$, then it follows from \cite{w2}
that $\varphi$ has no $(n-1)$-dimensional fibers.
Therefore we will assume that $n\geq 4$;
thus $\dim W\geq 3$ and we infer that 
$\Pic W\cong \mathbb{Z}$
in both cases A) and B).
We have
$\det \mathcal{E}_W=\mathcal{O}_W(r+\alpha)$
for some non-negative integer $\alpha$.

We claim here that $H(\mathcal{E}_W(-1))$ is nef
and that $\mathcal{E}_W(-1)$ is nef.
Observe that 
$H(\mathcal{E}_W(-1))|_{W_s}=\mathcal{O}_{W_s}$.
Since $\Pic P_W\cong \mathbb{Z}^{\oplus 2}$,
$\overline{{\rm NE}}(P_W)$
is spanned by two rays;
one is spanned by a curve in a fiber of $\pi$,
and the other is spanned by a curve in a fiber,
e.g., $W_s$, of $\varphi_W$.
Hence we have $\overline{{\rm NE}}(P_W)={\rm NE}(P_W)$,
and thus $H(\mathcal{E}
_W(-1))$ is nef.

Next we claim that $\alpha$ is positive.
Assume to the contrary that $\alpha=0$:
$\det \mathcal{E}_W=\mathcal{O}_W(r)$.
We have
$\mathcal{E}_W\cong \mathcal{O}_W(1)^{\oplus r}$
by \cite[Lemma 3.6.1]{w3} in the case A)
and by Theorem~\ref{uniform} (1) in the case B).
Since the projection $P_W\ra W$ has a section $W_s$
contracted to a point by $\varphi_W$,
the image $\im \varphi_W$
of $\varphi_W$ has dimension $r-1$,
and $\varphi_W$ is of fiber type.
If $\varphi$ is birational, this implies that 
$\im \rho$ is the exceptional locus of $\varphi$.
On the other hand, the divisor $\im \rho$
is nef because it is the pull back of a non-zero effective 
divisor $\im \pi\circ f$ on $M$ with $\Pic M\cong \mathbb{Z}$.
This is a contradiction.
Therefore $\varphi$ is of fiber type.
Let $F$ be a general fiber of $\varphi$.
Then we find that 
$\pi(F)\cap \im \pi\circ f$ is not empty
since $F$ is positive dimensional
and $\im \pi\circ f$ is an ample divisor on $M$.
Hence $F\cap \pi^{-1}(\im \pi\circ f)\neq \emptyset$,
so that $\rho^{-1}(F)\neq \emptyset$. 
This implies that $\varphi_W$ is dominant,
and thus $\dim N=r-1$.
This is a contradiction.
Therefore $\alpha$ is positive.

Let $R$ and $C_0$ be as in \S~\ref{the common setup}.
We will show that if 
we may assume that $\pi(C_1)=C_0$,
then $\varphi$ has no $(n-1)$-dimensional fibers;
Proposition~\ref{no n-1 fiber}
then follows from
Corollary~\ref{comparison for lr1=n-1}.

Suppose that $\pi(C_1)=C_0$.
Since $C_1$ is an image in $P$ of a line in $W$,
we have $r+\alpha=\pi^*\det\mathcal{E}.C_1
=\det\mathcal{E}.C_0$
by Corollary~\ref{C_1'to pi(C_1') birational}.
We see moreover that 
\[\alpha=l(R)-l(R_1)\leq 2\]
by Lemma~\ref{lr1-lr=a}.

Suppose that the case A) holds.
Since $\det\mathcal{E}_W(-1)=\mathcal{O}_W(\alpha)$,
we have
\[K_W+ \det \mathcal{E}_W(-1)=
\mathcal{O}_W(-n+1+\alpha).
\]
Since $n\geq 4$ and $\alpha\leq 2$,
we have 
$-n+1+\alpha\leq -1$.
Since we also have $1\leq \alpha$,
we have $H^i(W,K_W+ \det \mathcal{E}_W(-1))=0$ for all $i$.
Since $-(K_W+\det \mathcal{E}_W(-1))$ is ample,
it follows from
$-K_{P_W}=rH(\mathcal{E}_W(-1))-
\pi_W^*(K_W+\det \mathcal{E}_W(-1))$
that $-K_{P_W}$ is ample.
Since $H(\mathcal{E}
_W(-1))$ is nef,
this implies that $H(\mathcal{E}_W(-1))^{\otimes m}$
is spanned for all $m\gg 0$.

Since $\dim N= r$, we apply 
Corollary~\ref{pswcor} (2) to the pair 
$(W,\mathcal{E}_W(-1))$
to see that $\mathcal{E}_W(-1)$ is spanned,
that $h^0(\mathcal{E}_W(-1))=r+1$,
and that $\mathcal{E}_W(-1)$ fits into 
the following exact sequence
\[0\to (\mathcal{E}_W(-1))^*\to 
\mathcal{O}\otimes H^0(\mathcal{E}_W(-1))^*
\to \mathcal{O}(\alpha)\to 0.
\]
Since $\dim W=n-1$, this implies that $n\leq r+1$.
On the other hand, we have $n-1\geq r$ 
as we noticed
at the beginning of the proof.
Therefore we have $n=r+1$,
and moreover 
$H^0(\mathcal{E}_W(-1))^*\to H^0(\mathcal{O}(\alpha))$ is injective.
Hence $H^0((\mathcal{E}_W(-1))^*)=0$.
Note here that the natural surjection
$\pi^*\mathcal{E}\to L$ induces
a surjection $\mathcal{E}_W(-1)\to \mathcal{O}_W$.
Since $\mathcal{E}_W(-1)$ is spanned,
this implies that $\mathcal{E}_W(-1)\cong
\mathcal{O}_W\oplus \mathcal{F}$
for some vector bundle $\mathcal{F}$
of rank $r-1$ on $W$.
Hence we have $H^0((\mathcal{E}_W(-1))^*)\neq 0$
since $(\mathcal{E}_W(-1))^*\cong
\mathcal{O}_W\oplus \mathcal{F}^*$.
This is a contradiction.

Suppose that the case B) holds.
If $\alpha=2$, then $l(R)=n+1$,
and $M\cong \mathbb{P}^n$ by Theorem~\ref{cm}.
Thus $C_0$ is a line in $\mathbb{P}^n$.
Hence $\det\mathcal{E}.C_0=r+2$ implies that 
$\det\mathcal{E}\cong\mathcal{O}(r+2)$.
Proposition~\ref{pswa} (1) then implies
that $\varphi$ has an $n$-dimensional fibers,
since $n\geq 4$. This is a contradiction.
Therefore $\alpha=1$.
Since $\mathcal{E}_W(-1)$ is nef,
this implies that $\mathcal{E}_W(-1)$
is a uniform vector bundle of type $(0,\ldots,0,1)$.
Thus $\mathcal{E}_W(-1)$ is spanned
by Theorem~\ref{uniform} (2).
Since we have a surjection $\mathcal{E}_W(-1)\to \mathcal{O}_W$,
this implies that $\mathcal{E}_W(-1)\cong
\mathcal{O}_W\oplus \mathcal{F}$
for some spanned vector bundle 
$\mathcal{F}$
of rank $r-1$ on $W$.
Moreover we have
$\mathcal{F}\cong 
\mathcal{O}^{\oplus (r-2)}\oplus \mathcal{O}(1)$
by Theorem~\ref{uniform},
since $\rank \mathcal{F}=r-1\leq n-2<n-1=\dim W$.
Therefore we have 
$\mathcal{E}_W(-1)\cong 
\mathcal{O}^{\oplus (r-1)}\oplus \mathcal{O}(1)$
and we infer that the contraction $h$ different from
the projection $P_W\to W$ is the blowing-up
of $\mathbb{P}^{n+r-2}$ along a linear
subspace of dimension $r-2$.
Since $W_s$ is contracted by $\varphi_W$,
we see that $\varphi_W$ is factored as $\varphi_W=i\circ h$
for some finite morphism $i:\mathbb{P}^{n+r-2}\to N$.
Suppose that $\varphi$ is of fiber type.
Since $l(R_1)=n-1$, it follows from 
inequality~(\ref{ub of r1})
in \S~\ref{the common setup} that 
a general fiber of $\varphi$ is $(n-2)$-dimensional
and $\dim N=r+1$.
Therefore $n+r-2\leq r+1$, i.e., $n\leq 3$.
This, however, contradicts the assumption $n\geq 4$.
Suppose that $\varphi$ is birational.
Since $l(R_1)=n-1$, it follows again from 
inequality~(\ref{ub of r1})
in \S~\ref{the common setup} that 
$\varphi$ is divisorial.
Moreover every positive dimensional fiber of $\varphi$
is $(n-1)$-dimensional.
Let $E$ be the exceptional divisor of $\varphi$.
Since any nonzero effective divisor on $M$ is ample,
we infer that $\pi(E)=M$.
This implies that there exits an $(n-1)$-dimensional fiber $F$
of $\varphi$ such that $\pi(F)\nsubseteq \pi(W)$.
Therefore $\varphi_W$ must have $(n-2)$-dimensional fibers.
This contradicts the description $\varphi_W=i\circ h$.
This completes the proof of the proposition.
\end{proof}

\begin{prop}\label{lr1=n-1&n-2fiber}
Suppose that $\tau \geq 1$,
that $L.C_1=1$,
that $l(R_1)=n-1$,
and that every fiber of $\varphi$ has 
dimension $\leq n-2$.
Then $N$ is smooth of dimension $r+1$ 
and $\varphi$ makes $(P,L)$ a scroll over $N$.
Moreover $(M,\mathcal{E})$
is one of the following:
\begin{enumerate}
\item $(\mathbb{P}^3,N(2))$,
where $N$ is a null correlation bundle
on $\mathbb{P}^3$;
\item $(\mathbb{Q}^3,\mathbf{E}(2))$,
where $\mathbf{E}$
is a spinor bundle on $\mathbb{Q}^3$;
\item $(\mathbb{Q}^4,\mathbf{E}(2))$,
where $\mathbf{E}$
is a spinor bundle on $\mathbb{Q}^4$;
\item $n=4$, $r=3$, $\tau =1$, and either of the following holds:
   \begin{enumerate}
    \item $\pi(C_1)$ is not a minimal extremal rational curve in $M$
    (i.e., in $R$);
    \item $\varphi (l)$ is not a minimal extremal rational curve in $N$,
    where $l$ is a line in a fiber of $\pi$.
   \end{enumerate}
\end{enumerate}
\end{prop}
\begin{proof}
We see immediately by inequality~(\ref{ub of r1})
in \S~\ref{the common setup}
that $\varphi$ is of fiber type
and that every fiber of $\varphi$ has dimension $n-2$.
Thus $(\varphi^{-1}(z), L|_{\varphi^{-1}(z)})\cong
(\mathbb{P}^{n-2},\mathcal{O}(1))$ for a general point $z\in N$
by Theorem~\ref{ko} (1).
Hence Lemma~\ref{fujita lemma}
shows that $\varphi$ makes
$(P,L)$ a scroll over an $(r+1)$-dimensional manifold $N$.
Set $\varphi_*L=\mathcal{G}$;
$H(\mathcal{G})=L$ and $\mathcal{G}$
is an ample vector bundle of rank $n-1$ on $N$.

Note that $n-1=l(R_1)\geq \tau r\geq r$
by (\ref{lb of r1=}) in \S~\ref{the common setup},
since $\tau \geq 1$.

Let $W$ be a fiber of $\varphi$,
and set $\mathcal{E}_W=(\pi^*\mathcal{E})_W$
and $P_W=\mathbb{P}(\mathcal{E}_W)$.
Let $\rho :P_W \ra P$
be the morphism induced from $W\hookrightarrow P\to M$
by base change $\pi :P\ra M$.
Set $\varphi_W=\varphi\circ\rho$.
We have the following commutative diagram
\[
\begin{CD}
P_W@>{\rho}>>P@>{\varphi}>>N\\
@V{\pi_W}VV @VV{\pi }V  @. \\
W@>>>  M. @.
\end{CD}
\]
Since $W\hookrightarrow P\to M$
is finite by Lemma~\ref{finite of yz},
we see that $\rho$ is finite
and that $\mathcal{E}_W$ is ample.
Moreover we see, by the same argument 
as in the proof of Proposition~\ref{no n-1 fiber},
that $\mathcal{E}_W(-1)$ is nef.
Set $\det \mathcal{E}_W(-1)\cong \mathcal{O}_W(\alpha)$
for some non-negative integer $\alpha$.

We claim here that $\alpha>0$.
Assume, to the contrary, 
that $\alpha=0$.
Then we have $\mathcal{E}_W\cong \mathcal{O}_W(1)^{\oplus r}$
for any fiber $W=\varphi^{-1}(z)$ of $\varphi$.
Since $\dim W<n$, this contradicts Lemma~\ref{kantanlemma}.
Hence we have $\alpha>0$.

Suppose that we are not in the case (4) (a)
of the proposition.
Then Corollary~\ref{comparison for lr1=n-1}
implies that we may assume
that $\pi(C_1)=C_0$.
Note that $C_1$ is a line in $W$.
We have $r+\alpha=\pi^*\det\mathcal{E}.C_1
=\det\mathcal{E}.C_0$
by Corollary~\ref{C_1'to pi(C_1') birational}.
Moreover we have
\[\alpha=l(R)-l(R_1)\leq 2\] 
by Lemma~\ref{lr1-lr=a}.

If $\alpha=2$, then $l(R)=n+1$
and Theorem~\ref{cm} implies that 
$M\cong \mathbb{P}^n$.
Since $C_0$ is now a line in $\mathbb{P}^n$,
$\det\mathcal{E}.C_0=r+2$ implies that 
$\mathcal{E}\cong\mathcal{O}(r+2)$.
Therefore $n=3$
and $\mathcal{E}\cong N(2)$
by Proposition~\ref{pswa} (1),
since $\varphi$ has no $n$-dimensional fibers.
This is the case (1) of the proposition.

We will assume that $\alpha=1$ in the following.
We have $l(R)=n$.
Since $\mathcal{E}_W(-1)$ is nef,
this implies that $\mathcal{E}_W(-1)$
is a uniform vector bundle of type $(0,\dots,0,1)$.
Hence $\mathcal{E}_W(-1)$ is spanned by Theorem~\ref{uniform} (2).
On the other hand, $\mathcal{E}_W(-1)$
has an $\mathcal{O}_W$ as a quotient,
induced from the natural surjection
$\mathcal{E}_W\to L_W$.
Since $\mathcal{E}_W(-1)$ is spanned,
this implies that $\mathcal{E}_W(-1)
\cong \mathcal{O}_W\oplus \mathcal{F}$
for some spanned vector bundle $\mathcal{F}$
of rank $r-1$ on $W$.
Since $\mathcal{F}$ is a uniform vector bundle
of type $(0,\dots,0,1)$
of rank $r-1\leq n-2=\dim W$,
we infer that $\mathcal{F}$ is either
$\mathcal{O}^{\oplus (r-2)}\oplus \mathcal{O}(1)$
or $T_W(-1)$ by Theorem~\ref{uniform}.

Suppose that $\mathcal{E}_W(-1)\cong 
\mathcal{O}^{\oplus (r-1)}\oplus \mathcal{O}(1)$.
Then we have $n+r-3=\dim \im \varphi_W
\leq \dim N=r+1$, i.e., $n\leq 4$.
Thus we have 
\[3\geq n-1\geq \tau r\geq r\geq 2.\]
If $r=3$, then $\tau =1$ and $n=4$.
Since $\dim \im \varphi_W=\dim N$ in this case,
we see, moreover, that $N\cong \mathbb{P}^{r+1}$
by Lazarsfeld's theorem~\cite[Theorem~4.1]{l}.
Since we have 
$-K_P=(n-1)L-\varphi^*(K_N+\det\mathcal{G})$
in this case,
we infer that $K_N+\det\mathcal{G}=0$.
Therefore $\det\mathcal{G}\cong \mathcal{O}(5)
\cong \mathcal{O}(r+2)$.
Now Proposition~\ref{pswa} (1) implies that
$\pi$ has a fiber of dimension $=\dim N$, a contradiction.
Therefore we have $r=2$; we have $\det\mathcal{E}.C_0
=r+\alpha=3$.
If $n=4$, then $(-K_M-\det\mathcal{E}).C_0=1$
since $l(R)=4$.
Set $A=-K_M-\det\mathcal{E}$.
Then $A$ is ample
and we have $K_M+4A=0$.
Therefore we have $(M,A)
\cong (\mathbb{Q}^4,\mathcal{O}(1))$
by Theorem~\ref{ko} (2).
Thus $\det\mathcal{E}\cong \mathcal{O}(r+1)$.
Now Proposition~\ref{pswa} (2) implies that
$\mathcal{E}\cong \mathbf{E}(2)$.
This is the case (3) of the proposition.
If $n=3$, then $\tau =1$.
Since $l(R)=n$, it follows from \cite{w2}
that $(M,\mathcal{E})
\cong (\mathbb{Q}^3, \mathbf{E}(2))$.
This is the case (2) of the proposition.

Suppose that $\mathcal{E}_W(-1)
\cong \mathcal{O}_W\oplus T_W(-1)$.
Then we have $r=n-1$ and $\tau =1$.
Let $\varphi_W=h\circ g$ be
the Stein factorization of $\varphi_W$,
where $g:P_W\to \tilde{N}$ has connected fibers,
$\tilde{N}$ is normal,
and $h:\tilde{N}\to N$ is finite.
We see that $\tilde{N}\cong\mathbb{P}^{n-1}$
since $\mathcal{E}_W(-1)
\cong \mathcal{O}_W\oplus T_W(-1)$
and the section $W_s$ of $\pi_W$
corresponding to the quotient
$\mathcal{E}_W(-1)
\to \mathcal{O}_W$ is contracted by $\varphi_W$.
Set $y=g(W_s)$.
We see that every fiber of $\pi_W$
is mapped isomorphically onto an hyperplane
in $\tilde{N}\cong\mathbb{P}^{n-1}$
passing through $y$.

Assume that we are not in the case (4) (b)
either. Then we may assume, by the consideration above,
that $\mathcal{G}_F(-1)\cong \mathcal{O}_F\oplus 
T_F(-1)$ for any fiber $F$ of $\pi$. 
Moreover we may assume that $T_F(-1)$
is not a line bundle,
i.e., that $r-1=\dim F\geq 2$.
Now $(\mathcal{G}\otimes \mathcal{O}_{\tilde{N}})|_{H_y}$
is isomorphic to $\mathcal{O}_{H_y}(1)\oplus T_{H_y}$
for any hyperplane $H_y$ in $\tilde{N}\cong\mathbb{P}^{n-1}$
passing through $y$.
Note here that for any line in $\mathbb{P}^{n-1}$
there exists an hyperplane $H_y$ containing the line
and $y$, since $n-1\geq r\geq 3$.
Therefore $\mathcal{G}\otimes \mathcal{O}_{\tilde{N}}$
is a uniform vector bundle of type
$(1,\dots,1,2)$ of rank $n-1=\dim \tilde{N}$,
and thus it is isomorphic to either
$\mathcal{O}(1)^{\oplus (n-2)}\oplus \mathcal{O}(2)$
or $T_{\mathbb{P}^{n-1}}$
by Theorem~\ref{uniform} (2).
Here consider the composite $\tilde{\pi}$ 
of the morphism
$\mathbb{P}(\mathcal{G}\otimes
\mathcal{O}_{\tilde{N}})=
P\times_N\tilde{N}\to P$ and the projection
$\pi:P\to M$.
Note that $g:P_W\to \tilde{N}$ induces a finite
morphism $P_W\to P\times_N\tilde{N}$
and that any image via this finite morphism
of any fiber of $\pi_W$
is contracted by $\tilde{\pi}$.
Thus if $\mathcal{G}\otimes \mathcal{O}_{\tilde{N}}
\cong \mathcal{O}(1)^{\oplus (n-2)}\oplus \mathcal{O}(2)$,
then $\tilde{\pi}$
contracts an $(n-1)$-dimensional section 
corresponding to a quotient
$\mathcal{G}\otimes \mathcal{O}_{\tilde{N}}(-1)
\to \mathcal{O}$,
which implies that $\pi$ has an $(n-1)$-dimensional 
fiber, a contradiction.
If $\mathcal{G}\otimes \mathcal{O}_{\tilde{N}}
\cong T_{\mathbb{P}^{n-1}}$,
then the image of $\tilde{\pi}$ has 
dimension $n-1$.
Let $D_N$ denote a divisor 
defined as the image of $\tilde{N}\to N$
and set $D_P=\varphi^*D_N$.
Since $\dim \im \tilde{\pi}=n-1$,
we have $\dim \pi(D_P)=n-1$.
On the other hand,
we see that $\pi(D_P)=M$,
since $D_P.l=D_N.\varphi_*(l)>0$
for any line $l$ in a fiber of $\pi$.
This is a contradiction.
Therefore this case does not occur either.
\end{proof}

Finally we end this section with some examples,
which show that if we remove the assumption $\tau \geq 1$
then there are other examples in 
Propositions~\ref{no n-1 fiber} and \ref{lr1=n-1&n-2fiber}.
\begin{ex}
Let $\mathcal{E}$ be an ample vector bundle of rank $r$ on
$M=\mathbb{P}^2$ defined by the following exact sequence:
\[0\to \mathcal{O}(-1)\to \mathcal{O}(1)^{\oplus (r+1)}
\to \mathcal{E}\to 0.\]
Then $P=\mathbb{P}(\mathcal{E})$ is a divisor of bidegree $(2,1)$
in $\mathbb{P}^2\times \mathbb{P}^r$,
and $\varphi$ is the composite of the embedding
$P\subset \mathbb{P}^2\times \mathbb{P}^r$
and the projection $\mathbb{P}^2\times \mathbb{P}^r
\to \mathbb{P}^r$.
We see that $\varphi$ is a conic bundle
and that its discriminant locus $D$
is a hypersurface of degree three in $\mathbb{P}^r$.
Therefore $\varphi$ has a singular fiber
and $C_1$ is a line in a singular fiber of $\varphi$.
We have $L.C_1=1$ and $l(R_1)=1=n-1$, where $n=\dim M=2$.
Finally we have $\tau =3/(r+2)\leq 3/4<1$ since $r\geq 2$.
\end{ex}
\begin{ex}
Set $M=\mathbb{P}^3$ and $\mathcal{E}=\Omega_{\mathbb{P}^3}(3)$.
Then $\tau =4/5<1$.
\end{ex}

\section{The case $l(R_1)=n-2$}\label{-K=n-2shin}
Let $M$ be a Fano manifold of Picard number one
and $\mathcal{E}$ an ample vector bundle of 
rank $r$ on $M$.
Suppose that $\tau \geq 1$,
and let $\pi:P\to M$, $L$, $R_1$, $C_1$, and 
$\varphi:P\to N$
be as in \S~\ref{the common setup}.
Suppose that $L.C_1=1$,
that 
$l(R_1)=n-2$.
We have 
\[n-2=l(R_1)\geq \tau r\geq r\]
by (\ref{lb of r1=}) in \S~\ref{the common setup},
since $\tau \geq 1$.
In this section, we assume that $\tau r\ge n-2$
and prove Theorem~\ref{Main theorem}
in this case.

If $r=n-2$, then $\tau =1$, and 
we have the case (18) of 
Theorem~\ref{Main theorem}.
In the following, we assume that $r\leq n-3$.
Since $\tau r=n-2$, this implies that $\tau>1$.

Since we have 
\[(\tau -1)\det \mathcal{E}.\pi_*(C_1)=
-\pi^*(K_M+\det\mathcal{E})=
-K_P.C_1-rL.C_1=
n-2-r,\]
it follows from $\tau > 1$ that
$\det \mathcal{E}.\pi_*(C_1)=r$;
we see therefore that 
$C_1\to \pi(C_1)$ is birational and
that $\det \mathcal{E}.\pi (C_1)=r$.
Hence we have $-K_M.\pi (C_1)=\tau \det\mathcal{E}.\pi(C_1)=n-2$
and thus $l(R)=n-2$. Set $C_0=\pi(C_1)$.

Let $A_z$ be an irreducible component of 
a positive dimensional fiber $\varphi^{-1}(z)$ of $\varphi$,
and $f:W\ra P$ the composite of the normalization
$W\ra A_z$ and the inclusion $A_z\hookrightarrow P$. 
Since $\pi^*\det\mathcal{E}.C_1=r=rL.C_1$,
we have $\det (\pi^*\mathcal{E})|_{A_z}\cong rL|_{A_z}$.
Hence we have $\det \mathcal{E}_{W}\cong rL_{W}$,
where $\mathcal{E}_W$ and $L_W$ denotes 
$(\pi\circ f)^*\mathcal{E}$ and $f^*L$ respectively.

Suppose that $\dim A_z=n$.
Then, by Lemma~\ref{isom},
we see  
that $\pi\circ f: W\ra M$ is an isomorphism.

If, moreover, every fiber of $\varphi$ is $n$-dimensional,
it follows from Lemma~\ref{l(r1)=s'}
that $M$ is a Fano manifold of index $n-2$,
i.e., a Mukai manifold,
and $\mathcal{E}\cong A^{\oplus r}$
for the ample generator $A$ of $\Pic M$.
This is the case (17) of Theorem~\ref{Main theorem}.

If, moreover, a general fiber of $\varphi$ has dimension less than $n$,
then we have $r\geq 2$.
Let $C_1'$ be the rational curve on $W$ 
corresponding to $C_0$
via the isomorphism $\pi\circ f: W\cong M$.
Since $-K_M.C_0=n-2$, we have $-K_W.C_1'=n-2$.
Let $\tilde{C_1}\to C_1$,
$\tilde{C_1'}\to C_1'$,
and $\tilde{C_0}\to C_0$
be the normalizations.
We can regard $\tilde{C_1}$ and $\tilde{C_1'}$
as sections of 
$\mathbb{P}(\mathcal{E}\otimes \mathcal{O}_{\tilde{C_0}})
\to \tilde{C_0}$.
Since $\det \mathcal{E}.C_0=r$, we have
$\mathcal{E}\otimes \mathcal{O}_{\tilde{C_0}}
\cong \mathcal{O}(1)^{r}$.
Hence $\tilde{C_1}$ and $\tilde{C_1'}$ are rationally equivalent,
so that $L.C_1'=1$.
Therefore we have $K_W+(n-2)L_W=0$.
In particular we see that $h^n(W, -(n-2)L_W)=1$.
On the other hand, since a general fiber
of $\varphi$ has dimension $<n$, we have
$h^n(W,-tL_W)=0$ for all $t\leq n-2$
by \cite[Lemma 4]{yz}.
This is a contradiction.
Hence this case does not occur.

Let us assume that $\varphi$ has no $n$-dimensional fibers
in the following; hence we may assume that 
$n\geq r+3\geq 5$.
Since $\dim A_z\leq n-1$,
it follows from 
inequality \eqref{ub of r1} in \S~\ref{the common setup}
that for a general point $x\in A_z$
there exists a rational curve $C_1'$
such that $-K_P.C_1'\leq n$.
Since $-K_P.C_1'=(n-2)L.C_1'$,
if $L.C_1'\geq 2$, then we have $2(n-2)\leq n$,
i.e., $n\leq 4$.
This is a contradiction.
Hence we have $L.C_1'=1$ for such a curve $C_1'$; 
this implies that we may regard $C_1'$ as $C_1$,
and thus we may assume that 
for a general point $x\in A_z$ there exists
a $C_1$ passing through $x$.

Suppose that $\varphi$ is birational.
Since we have 
$\det \mathcal(\pi^*\mathcal{E})|_{A_z}\cong rL|_{A_z}$,
we have $\mathcal{E}_{\tilde{C_1}}\cong
\mathcal{O}_{\tilde{C_1}}(1)^{\oplus r}$,
where $\tilde{C_1}\to C_1$ denotes the normalization.
This implies that 
$\pi^{-1}(\pi(C_1))$ is contained in 
the exceptional locus of $\varphi$.
Since an irreducible family of $C_1$ dominates $A_z$,
we infer that $\pi^{-1}(\pi(A_z))$ is contained in 
the exceptional locus of $\varphi$.
Suppose that $\varphi$ has an $(n-1)$-dimensional fiber
and that $\dim A_z=n-1$.
Then $\codim (\pi^{-1}(\pi(A_z)),P)=1$
and thus $\varphi$ is divisorial.
Since the exceptional locus of $\varphi$ is irreducible
if $\varphi$ is divisorial,
we see that $\pi^{-1}(\pi(A_z))$ is the exceptional
divisor of $\varphi$.
On the other hand, since $\Pic M\cong \mathbb{Z}$,
a nonzero effective divisor $\pi(A_z)$ is ample.
Hence $\pi^{-1}(\pi(A_z))$ is nef.
This is a contradiction.
Therefore this case does not happen.
Suppose that
every positive dimensional fiber of $\varphi$ has 
dimension $\leq n-2$.
Then it follows from inequality~(\ref{ub of r1})
in \S~\ref{the common setup}
that $\varphi$ is divisorial
and that every positive dimensional fiber of $\varphi$
is equidimensional of dimension $n-2$.
Let $E$ denote the exceptional divisor of $\varphi$,
and $A_z$ any irreducible component of 
any positive fiber $\varphi^{-1}(z)$ of $\varphi$.
Then, as we have seen as above,
$\pi^{-1}(\pi(A_z))$ is contained in
the exceptional locus $E$ of $\varphi$.
Hence we have $\pi^{-1}(\pi(E))\subseteq E$
and thus we have $\pi^{-1}(\pi(E))=E$.
This is a contradiction for the same reason as above.
Hence this case does not occur either.
Therefore $\varphi$ cannot be birational.

Suppose that $\varphi$ is of fiber type
and that a general fiber of $\varphi$ has dimension $n-1$;
hence we see that $\dim N=r$.
We will denote by $W$ 
a general fiber of $\varphi$.
Since $K_{W}+(n-2)L_{W}=0$,
$(W,L_W)$ is a Del Pezzo manifold.
The classification of Del Pezzo manifolds 
(see, e.g., \cite[(8.11)]{fb}) then 
implies that $\Pic W\cong \mathbb{Z}$ since $n\geq 5$.
Set $P_W=\mathbb{P}(\mathcal{E}_{W})$.
Note that
$-K_{P_{W}}=rH(\mathcal{E}_{W})+(n-2-r)\pi_W^*L_{W}$,
where $\pi_W:P_{W}\to W$ is the projection,
because $\det\mathcal{E}_{W}\cong rL_{W}$.
Hence we see that $P_{W}$
is a Fano manifold of Picard number two.
Let $\rho :P_W\ra P$
be the morphism induced from $\pi\circ f:W\ra M$
by base change $\pi :P\ra M$.
Observe here that 
the image $\im \varphi\circ\rho$ 
of $P_{W}$ via the composite
$\varphi\circ \rho$ of $\rho$ and $\varphi$
has dimension either $r-1$ or $r$.
Let $W_s$ denote the section of $\pi_W:P_{W}\to W$ 
induced from the inclusion $f:W\hookrightarrow P$;
$W_s$ is also a fiber of $\varphi\circ \rho$.
Note here that $H(\mathcal{E}_{W}\otimes L_{W}^{-1})|_{W_s}=0$.
Hence  $H(\mathcal{E}_{W}\otimes L_{W}^{-1})$ is nef,
and we infer that there exists a positive integer $m_0$
such that $H(\mathcal{E}_{W}\otimes L_{W}^{-1})^{\otimes m}$ 
is spanned for all integers $m\geq m_0$.
We have $H^i(P_W,H(\mathcal{E}_W\otimes L_W^{-1})^{\otimes t})=0$
for $i>0$ and $t>0$.
Moreover we see that
$H^{i}(W, K_W\otimes \det (\mathcal{E}_W\otimes L_W^{-1}))=0$ 
for $i<\dim W$
and that 
$h^{n-1}(W, K_W\otimes \det (\mathcal{E}_W\otimes L_W^{-1}))=1$.
Therefore it follows from Corollary~\ref{pswcor}
that $\mathcal{E}_W\otimes L_W^{-1}\cong \mathcal{O}_W^{\oplus r}$.
This however contradicts Lemma~\ref{kantanlemma}
because $\dim W<n$.

Suppose that a general fiber $W$
of $\varphi$ has dimension $n-2$.
Since $K_W+(n-2)L_W=0$,
it follows from Theorem~\ref{ko} (2)
that $(W, L_W)\cong (\mathbb{Q}^{n-2},\mathcal{O}(1))$.
Since $\det \mathcal{E}_{W}\cong rL_{W}$,
we see that $\mathcal{E}_W\cong \mathcal{O}_W(1)^{\oplus r}$
by \cite[Lemma 3.6.1]{w3}.
This contradicts Lemma~\ref{kantanlemma}.

Suppose that a general fiber $W$
of $\varphi$ has dimension $n-3$.
Since $K_W+(n-2)L_W=0$,
it follows from Theorem~\ref{ko} (1)
that $(W, L_W)\cong (\mathbb{P}^{n-3},\mathcal{O}(1))$.
Since $\det \mathcal{E}_{W}\cong rL_{W}$,
we see $\mathcal{E}_W\cong \mathcal{O}_W(1)^{\oplus r}$
by Theorem~\ref{uniform} (1);
this contradicts Lemma~\ref{kantanlemma}.

\section{Lemma on a quadric fibration}\label{quadric fibration}
Let $\psi:M\to S$ be a projective morphism
of smooth algebraic varieties,
and suppose that every closed fiber of $\psi$
is an $m$-dimensional smooth 
hyperquadric in $\mathbb{P}^{m+1}$.
We will call such a fibration $\psi$
a quadric fibration.
Let $p$ be an arbitrary closed point of $S$.
As an analogue of Brauer-Severi schemes,
it seems to be reasonable to consider
the following problem:
\begin{problem}
Does there exist an \'{e}tale finite covering
$U\to U'$ of some neighborhood $U'$ of $p$
such that on the fiber product $M\times_SU$, denoted by $M_U$,
there exists a line bundle $A$ such that 
$A|_{\tilde{F}}=\mathcal{O}(1)$ for any
closed fiber $\tilde{F}$ of $M_U\to U$~? 
\end{problem}
If such an \'{e}tale covering of $S$ exists,
we will say, in this paper, that 
$\psi$ is a quadric fibration in the \'{e}tale topology.
Simply replacing the \'{e}tale topology
to the complex topology,
we will also use the phrase ``a quadric fibration
in the complex topology''.

The following lemma gives a sufficient condition for
the affirmative answer of this problem. 
In general, I do not know its answer.
\begin{lemma}\label{quadricbundle}
Let $\psi$ be as above.
Assume moreover that 
the relative Picard number $\rho (M/S)$ is one.
Suppose that there exists a locally free sheaf $\mathcal{E}$
of rank $r\leq m$ on $M$ such that 
$\mathcal{E}|_F\cong \mathcal{O}(1)^{\oplus r}$
for any closed fiber $F$ of $\psi$.
Then $\psi$ is a quadric fibration in the \'{e}tale
(or complex) topology.
\end{lemma}
\begin{proof}
Denote by $P$ the projective space bundle $\mathbb{P}(\mathcal{E})$
associated to $\mathcal{E}$.
Let $\pi\colon P\to M$ be the projection
and $L$ the tautological line bundle on $P$
associated to $\mathcal{E}$.
By assumption, we see that 
$L$ is ample over $S$ and that 
the relative Picard number $\rho(P/S)$
of $P$ over $S$ is two.
Note here that $-(K_M+\det \mathcal{E})$ is $\psi$-nef
since $m\geq r$.
Thus $-K_P$ is $\psi\circ \pi$-ample.
Therefore $\overline{{\rm NE}}(P/S)$
is spanned by two extremal rays $R_{\pi}$ and $R_1$,
where $R_{\pi}$ corresponds to the projection $\pi \colon P\ra M$
and $R_1$ is the other ray.

Let $\varphi \colon P\to N$ be the contraction morphism of $R_1$
over $S$ and $\pi'\colon N\to S$ the structural morphism.
We have the following commutative diagram
\begin{equation}\label{quadric diagram}
\begin{CD}
P         @>{\varphi}>> N            \\
@VV{\pi}V               @VV{\pi '}V \\
M         @>{\psi}>>    S.
\end{CD}
\end{equation}
Let $C_1$ be 
a minimal extremal rational curve of $R_1$.
Then $\pi(C_1)$ is a rational curve contracted by $\psi$.
We have 
\[-K_P.C_1=
rL.C_1+((m/r)-1)\det\mathcal{E}.\pi_*(C_1).\]
Therefore $l(R_1)\geq m$
since $m/r\geq 1$.
Note that, although $P$ itself may not be complete,
we can apply Theorem~\ref{wisniewski}
since $P$ is smooth and projective over $S$.
Therefore,
as in \S~\ref{the common setup},
we see that $l(R_1)\leq \dim F(\varphi)+1\leq m+1$
for any irreducible component $F(\varphi)$ 
of any positive dimensional closed fiber of $\varphi$.
Denote by $F$ the closed fiber of $\psi$ containing $\pi(F(\varphi))$.

We claim here $F(\varphi)_{\textrm{red}}$ is isomorphic to $F$ via $\pi$.
Since $\mathcal{E}|_F\cong \mathcal{O}(1)^{\oplus r}$,
we see that $F(\varphi)$ is contained in 
$\mathbb{P}(\mathcal{E}|_F)=F\times \mathbb{P}^{r-1}$.
Note here that $F(\varphi)$ is contracted
by $\varphi|_{\mathbb{P}(\mathcal{E}|_F)}$
which is different from the projection 
$\mathbb{P}(\mathcal{E}|_F)\to F$.
If $m\neq 2$, then $\Pic (F)\cong \mathbb{Z}$,
so that we see immediately that 
$\varphi|_{\mathbb{P}(\mathcal{E}|_F)}$
is factored as 
$\varphi|_{\mathbb{P}(\mathcal{E}|_{F})}
=g\circ f$,
where $f:\mathbb{P}(\mathcal{E}|_{F})=F\times \mathbb{P}^{r-1}
\to \mathbb{P}^{r-1}$
is the projection 
and $g:\mathbb{P}^{r-1}\to {\pi'}^{-1}(\psi(F))$
is a finite unibranch
(i.e., set-theoretically one to one) morphism.
Hence we infer that $F(\varphi)_{\textrm{red}}\cong F$.
If $m=2$, then $F\cong \mathbb{P}^1\times \mathbb{P}^1$.
If $\dim F(\varphi)=2$, then we see, as above, that
$F(\varphi)_{\textrm{red}}\cong F$.
We will show that the case $\dim F(\varphi)=1$ does not happen.
Assume to the contrary that $\dim F(\varphi)=1$. 
Then $\varphi|_{\mathbb{P}(\mathcal{E}|_F)}$
is factored as 
$\varphi|_{\mathbb{P}(\mathcal{E}|_{F})}
=g\circ f$,
where $f:\mathbb{P}(\mathcal{E}|_{F})=\mathbb{P}^1\times
\mathbb{P}^1\times \mathbb{P}^{r-1}
\to \mathbb{P}^1\times\mathbb{P}^{r-1}$
is a projection 
and $g:\mathbb{P}^1\times\mathbb{P}^{r-1}\to {\pi'}^{-1}(\psi(F))$
is a finite unibranch morphism.
Hence 
$F(\varphi)_{\textrm{red}}\cong \mathbb{P}^1$.
Since $\varphi$ is different from $\psi$,
we see that $r\geq 2$.
This implies that $r=2$,
since $r\leq m=2$.
Moreover we see, by  Theorem~\ref{wisniewski}, that
$l(R_1)=m$, that $L.C_1=1$,
and that $\varphi$ is of fiber type.
If $\varphi$ has a two-dimensional 
irreducible component $F_{\varphi}$
(with reduced structure) 
of a positive dimensional fiber of $\varphi$,
then it follows from Lemma~\ref{Fujita Ye-Zhang}
that $\tilde{F_{\varphi}}\cong \mathbb{P}^2$
where $\tilde{F_{\varphi}}\to F_{\varphi}$ is the normalization.
On the other hand, we must have 
$F_{\varphi}\cong F\cong \mathbb{P}^1\times \mathbb{P}^1$
as in the case $\dim F(\varphi)=2$.
This is a contradiction.
Hence every fiber of $\varphi$ has dimension one.
Note here that $\mathcal{E}_{F(\varphi)_{\textrm{red}}}
\cong \mathcal{O}(1)^{\oplus 2}$.
Denote by $z$ the point $\varphi(F(\varphi))$.
We see that 
$\varphi(\pi^{-1}(\pi(\varphi^{-1}(z))))$ is 
of codimension one in ${\pi'}^{-1}(\psi(F))$.
Now take a subvariety $S'$ of $N$ such that 
$S'\to S$ is finite over
a Zariski dense open subset $S_0$ of $S$ containing ${\pi'}(z)$;
we can take such a subvariety $S'$ since $\pi'$ is projective.
Set $M_0=\psi^{-1}(S_0)$.
Since $M_0$ is a Zariski dense open subset of $M$
and $M$ is smooth and algebraic,
every Cartier divisor on $M_0$ can be extended to 
a Cartier divisor on $M$.
Hence we see that the restriction map
$\Pic (M/S)\to \Pic (M_0/S_0)$ is surjective.
Therefore the relative Picard number $\rho (M_0/S_0)$ is also one.
Replacing $S$ with $S_0$, we may assume that $S'\to S$ is finite.
Set $D_M=\pi (\varphi^{-1}(S'))$ with the reduced structure.
Then $D_M$ is a prime divisor on $M$.
Set $D_P=\pi^*D_M$; $D_P$ is also a prime divisor on $P$.
Now set $D_N=\varphi(D_P)$.
Then we see, by the consideration above, that
$D_N$ is also a prime divisor on $N$
and that $\varphi^*(D_N)=D_P$.
Since $\rho (P/S)=2$, this implies that $D_P=0$,
i.e., $D_P=\emptyset$ as sets, which is a contradiction.
Therefore the case $\dim F(\varphi)=1$ does not happen,
and thus we conclude that 
$F(\varphi)_{\textrm{red}}\cong F$ via $\pi$;
as consequences, we see also that 
$L|_{F(\varphi)_{\textrm{red}}}\cong \mathcal{O}(1)$
since $\mathcal{E}|_F\cong \mathcal{O}(1)^{\oplus r}$.
We also observe that every fiber of $\varphi$ is 
irreducible since the finite morphism $g$ above
is unibranch.

Now we apply Fujita's argument in the proof of 
\cite[Lemma (2.12)]{fs} to show that 
every fiber of $\varphi$ is reduced.
Let $z$ be an arbitrary point of $N$,
and let $F(\varphi)$ be a total fiber of $\varphi$ over $z$;
we know that $F(\varphi)_{\textrm{red}}$ is a smooth hyperquadric
and that $L|_{F(\varphi)_{\textrm{red}}}\cong \mathcal{O}(1)$.
In this paragraph, we may assume that $L$ is ample
by replacing $L$ with tensor product
of $L$ with the pull back of sufficiently ample line bundle
on $N$.
Take a large integer $a$ such that $aL$ is very ample
and let $D_1, \dots, D_m$ be a general member of $\lvert aL\rvert$.
Set $M'=D_1\cap \dots \cap D_m$.
Then $M'$ is nonsingular and $M'\cap F(\varphi)_{\textrm{red}}$
is a nonsingular subscheme consisting of $2a^m$ points.
Take a small enough neighborhood $V$ of $z$ with respect to 
the metric topology such that 
any connected component $V_{\lambda}$ of $\varphi^{-1}(V)\cap M'$
meets $F(\varphi)_{\textrm{red}}$ at only one point.
Let $\varphi_{\lambda}$ be the restriction of $\varphi$ to 
$V_{\lambda}$. 
We may assume that $\varphi_{\lambda}$ is a finite
morphism of degree $m_{\lambda}$.
Since the number of $\lambda$'s are equal to 
$\sharp (M'\cap F(\varphi))=2a^m$,
we see that $\deg (\varphi|_{M'})=
\sum_{\lambda}m_{\lambda}\geq 2a^m$.
On the other hand, a general fiber $F_g$ of $\varphi$
is a smooth hyperquadric and $L|_{F_g}\cong \mathcal{O}(1)$.
Hence $\deg (\varphi|_{M'})=2a^m$.
Therefore $m_{\lambda}=1$ for all $\lambda$,
and thus $\varphi_{\lambda}:V_{\lambda}\to V$ is 
bimeromorphic. By the analytic version of Zariski Main Theorem,
we infer that $\varphi_{\lambda}$ is biholomorphic.
Hence $N$ is smooth. Furthermore we see 
$F(\varphi)\cap V_{\lambda}$ consists of only one point
scheme-theoretically since $\varphi_{\lambda}$
is biholomorphic.
This implies that $F(\varphi)$ is generically reduced.
Since $M$ and $N$ are smooth
and $\dim F(\varphi)=\dim M-\dim N$,
we infer that $F(\varphi)$ is Cohen-Macaulay. 
Therefore $F(\varphi)$ is reduced.

By the consideration above, we see that 
for any closed fiber $F$ of $\psi$
the induced morphism 
$\mathbb{P}(\mathcal{E}|_F)\to {\pi'}^{-1}(\psi(F))$
is nothing but the projection 
$\mathbb{P}(\mathcal{E}|_F)=F\times \mathbb{P}^{r-1}
\to \mathbb{P}^{r-1}$:
every closed fiber of $\pi'$ is isomorphic to $\mathbb{P}^{r-1}$.
In particular, $\pi'$ is flat.
Here we apply Theorem~\ref{Brauer-Severi}.
Let $p$ be a closed point of $S$.
For a suitable finite \'{e}tale covering $U\to U'$ over a small open
neighborhood $U'$ of $p$,
the fiber product $N\times_S U$, denoted by $N_U$,
is isomorphic to $\mathbb{P}^{r-1}\times U$.
Let $\mathcal{O}_{N_U}(1)$ be the tautological line bundle
associated to the projective space bundle $N_U\to U$.
Consider the following commutative diagram
obtained from the diagram~(\ref{quadric diagram})
by the base change $U\to U'\subset S$.
\[
\begin{CD}
P_U         @>{\varphi_U}>> N_U            \\
@VV{\pi_U}V               @VV{\pi '_U}V \\
M_U         @>{\psi_U}>>    U.
\end{CD}
\]
Now the line bundle $L\otimes \varphi_U^*\mathcal{O}_{N_U}(-1)$
is trivial on each closed fiber of $\pi_U$.
Therefore we obtain a line bundle $A$ on $M_U$ such that 
$L\otimes \varphi_U^*\mathcal{O}_{N_U}(-1)=\pi_U^*A$.
Since $A|_{\tilde{F}}\cong \mathcal{O}(1)$ for every closed fiber 
$\tilde{F}$ of $\psi_U$, this implies that
$\psi$ is a quadric fibration in the \'{e}tale topology.

Finally, for the case of the complex topology,
it is easy to see that
$\pi'$ is locally trivialized in the complex topology.
Hence the same argument as above shows that
$\psi$ is a quadric fibration in the complex topology.
This completes the proof.
\end{proof}

\section{The case $\dim F-1<\tau r\leq \dim F$
}\label{dim F-1<tau r leq dim F}
Let $M$, $\mathcal{E}$, and $\psi :M\ra S$
be as in \S~\ref{overview}.

If $\psi$ is of fiber type
and $\dim F-1<\tau r\leq \dim F$,
we can give the following proposition.
\begin{prop}\label{relative=n/r}
Suppose that $\tau \geq 1$,
that $\psi$ is of fiber type,
and that $\dim F-1<\tau r\leq \dim F$
for a general fiber $F$ of $\psi$.
Set $d=\ulcorner \tau r\urcorner=\dim F$.
Then $S$ has only rational Gorenstein singularities.
Let $U$ denote the largest open subset of $S$
such that $\psi^{-1}(U)\to U$ is smooth.
Then we have the following:
\begin{enumerate}
\item $\psi^{-1}(U)\to U$ is 
a $\mathbb{P}^d$-fibration in the \'{e}tale (or complex) topology.
Moreover we have $\codim (S\setminus U, S)\geq 3$.
In particular we see that 
  \begin{enumerate}
     \item if $d=n-1$, then 
     $\psi$ is the $\mathbb{P}^{n-1}$-bundle in the Zariski topology;
     \item if $d=n-2$, then
     every closed fiber of $\psi$ is isomorphic to $\mathbb{P}^{n-2}$;
     \item if $d=n-3$, then
     $\psi$ has at most finite number of singular fibers.
  \end{enumerate}
Concerning $\mathcal{E}$ and the value of 
$\tau r$, we have either of the following:
  \begin{enumerate}
     \item[\textup{(i)}] $\mathcal{E}|_F\cong T_{\mathbb{P}^d}$ 
     for every closed fiber $F$ of $\psi^{-1}(U)\to U$,
     and $\tau r=d$;
     \item[\textup{(ii)}]  $\mathcal{E}|_F\cong 
     \mathcal{O}(1)^{\oplus (r-1)}\oplus \mathcal{O}(2)$
     for every closed fiber $F$ of $\psi^{-1}(U)\to U$,
     and and $\tau r=(d+1)r/(r+1)>d-1$
     (and hence $d\geq r>(d-1)/2$);
  \end{enumerate}
\item  $\psi^{-1}(U)\to U$ is a
quadric fibration in the \'{e}tale 
(or complex) 
topology (see \S~\ref{quadric fibration}),
and for every closed fiber $F$ of $\psi^{-1}(U)\to U$
we have $\mathcal{E}|_F\cong 
\mathcal{O}(1)^{\oplus r}$,
and $\tau r=d$.
Moreover we have $\codim (\Sing S, S)\geq 3$ if $d\geq 3$.
Furthermore if $d=n-1$ then
there exists a line bundle $\mathcal{O}_M(1)$ on $M$
such that $\mathcal{O}_M(1)|_F\cong \mathcal{O}_{\mathbb{Q}}(1)$
for every fiber $F$ of $\psi$
and that $\mathcal{E}\cong \psi^*\mathcal{E}'\otimes
\mathcal{O}_M(1)$
for some vector bundle 
$\mathcal{E}'$
of rank $r$ on $S$.
\end{enumerate}
\end{prop}
\begin{proof}
It follows immediately from Proposition~\ref{smooth}
that $S$ has only rational Gorenstein singularities.

Let $F_0$ be a closed fiber of $\psi^{-1}(U)\ra U$.
Then $\dim F_0=d$ and $K_{F_0}+\tau \det \mathcal{E}|_{F_0}=0$;
since $\tau \geq 1$, it follows from 
the case $n-1<\tau r\leq n$
that $(F_0, \mathcal{E}|_{F_0})$ 
and the value of $\tau r$ are
one of the following:\\
1) $(\mathbb{P}^{d}, T_{\mathbb{P}^{d}})$, and $\tau r=d$;\\ 
2) $(\mathbb{P}^{d}, \mathcal{O}(1)^{\oplus (r-1)}\oplus \mathcal{O}(2))$,
and $\tau r=(d+1)r/(r+1)>d-1$ (and hence $r>(d-1)/2$);\\
3) $(\mathbb{Q}^{d}, \mathcal{O}_{\mathbb{Q}}(1)^{\oplus r})$,
and $\tau r=d$;\\
4) $(\mathbb{P}({\mathcal F}), 
H({\mathcal F})\otimes {\psi'}^{*}{\mathcal G})$,where ${\mathcal F}$ 
is a vector bundle of rank $d$ on a smooth proper curve $C$, 
$\psi' :\mathbb{P}({\mathcal F})\ra C$ is 
the projection, and ${\mathcal G}$ is a vector bundle 
of rank $r$ on $C$, and $\tau r=d$.

We will show that $F_0$ cannot be isomorphic to 
a scroll $\mathbb{P}(\mathcal{F})$ over a curve $C$.
Suppose, to the contrary, that $F_0=\mathbb{P}(\mathcal{F})$;
in particular we assume that $d\geq 2$.
Here we also assume that $F_0$ is not isomorphic to 
a quadric surface $\mathbb{Q}^2$.
Then we see first that 
$h^1(\mathcal{O}_C)=h^1(\mathcal{O}_{F_0})=0$ since $F_0$ is Fano.
Hence $C=\mathbb{P}^1$.
Thus $\mathcal{F}$ is a direct sum of line bundles:
$\mathcal{F}=\oplus_{i=1}^d\mathcal{O}(f_i)$.
Here we may assume that $f_1\leq \cdots \leq f_d$.
Since $F_0$ is Fano and $\omega_{F_0}^{-1}=
H(\mathcal{F})^{\otimes d}\otimes 
{\psi'}^*\mathcal{O}(2-\sum_{i=1}^df_i)$,
we have $0<df_1+2-\sum_{i=1}^df_i$.
Therefore $(f_1,\ldots,f_d)$ is either
$(f_1,\ldots,f_1)$ or $(f_1,\ldots,f_1,f_1+1)$.
Hence we may assume that $(f_1,\ldots,f_d)$
is either $(1,\ldots,1)$ or $(1,\ldots,1,2)$.
Since $d\geq 2$ and 
$\mathcal{E}|_{F_0}=H(\mathcal{F})\otimes {\psi'}^*\mathcal{G}$,
this implies that $\mathcal{G}$ is a direct sum
of line bundles of non-negative degrees.
Since $\tau r=d$, we observe that 
$0=K_{F_0}+\tau \det \mathcal{E}|_{F_0}
={\psi'}^{*}(K_{\mathbb{P}^1}+\det \mathcal{F}+\tau \det \mathcal{G})$.
This shows that $\deg \mathcal{F}-2\leq 0$,
and hence we have $\mathcal{F}\cong \mathcal{O}(1)^{\oplus 2}$.
Therefore $F_0=\mathbb{Q}^2$, a contradiction.

Let $F_1$ be any other fiber of $\psi^{-1}(U)\ra U$
over a closed point $p\in U$;
$F_1$ is isomorphic to either $\mathbb{P}^{d}$ or $\mathbb{Q}^{d}$.
Extending $\mathcal{O}_{F_1}(1)$
to a line bundle on a small neighborhood of $F_1$
in the complex topology and looking at
the degree and Fujita's delta genus with respect to
the extended line bundle,
we find that 
any fiber over a closed point
in a small neighborhood of $p$
is isomorphic to $F_1$.
Since $U$ is connected, this implies that
any smooth closed fiber $F_1$ of $\psi$ is isomorphic to $F_0$.
Thus $\psi^{-1}(U)\ra U$ is a $\mathbb{P}^{d}$-fibration
or a quadric fibration.
Moreover it follows from Theorem~\ref{Brauer-Severi}
and Lemma~\ref{quadricbundle} that 
$\psi^{-1}(U)\ra U$ is a $\mathbb{P}^{d}$-fibration
or a quadric fibration in the \'{e}tale (or complex) topology
(it is easy to see that 
$\psi^{-1}(U)\ra U$ is a $\mathbb{P}^{d}$-fibration
in the complex topology
in case every fiber of $\psi^{-1}(U)\ra U$ is
$\mathbb{P}^d$).

Let $S_2$ be the intersection
$D_1\cap\cdots \cap D_{n-d-2}$
of general very ample divisors
$D_1, \ldots, D_{n-d-2}$ of $S$.
We claim here that the restricted morphism
$\psi^{-1}(S_2)\ra S_2$
is the contraction morphism of an extremal ray
unless $d=2$ and 
$\psi^{-1}(U)\ra U$ is a quadric fibration.
Here we consider a conic fibration
as a $\mathbb{P}^{1}$-fibration;
hence if $\psi^{-1}(U)\ra U$ is a quadric fibration
we assume that $d\geq 2$.
Now suppose that 
$\psi^{-1}(S_2)\ra S_2$ is not elementary.
Then by applying Proposition~\ref{Proposition 1.3}
to an $(d+2)$-dimensional
manifold $\psi^{-1}(S_2)$, 
we have $d+1< (d+3)/2$ in case 
$\psi^{-1}(U)\ra U$ is a $\mathbb{P}^{d}$-fibration,
and $d< (d+3)/2$
in case $\psi^{-1}(U)\ra U$ is a quadric fibration.
Therefore $d=2$ and $\psi^{-1}(U)\ra U$ must be a quadric fibration.
This proves the claim.

Now suppose that the restricted morphism
$\psi^{-1}(S_2)\ra S_2$
is the contraction morphism of an extremal ray.
It follows from Proposition~\ref{smooth}
that $S_2$ is smooth.
Therefore $\codim (\Sing S, S)\geq 3$.

Let $S_1$ be a general very ample divisor on $S_2$.
Then the restricted morphism $\psi^{-1}(S_1)\ra S_1$
is a contraction of an extremal ray
unless $d=2$ and $\psi^{-1}(U)\ra U$ is a quadric fibration;
indeed, if this is not elementary, again
by applying Proposition~\ref{Proposition 1.3}
to an $d+1$-dimensional
manifold $\psi^{-1}(S_1)$, 
we have $d+1\leq (d+2)/2$
in case $\psi^{-1}(S_1)\ra S_1$
is a $\mathbb{P}^d$-fibration
and $d\leq (d+2)/2$
in case $\psi^{-1}(S_1)\ra S_1$
is a quadric fibration.
Therefore $d=2$ and $\psi^{-1}(U)\ra U$ is a quadric fibration.

Suppose that $\psi^{-1}(S_1)\ra S_1$ is elementary
in the following;
note here that if $d=n-1$ then $S_1=S$ and $\psi^{-1}(S_1)=M$
so that $\psi^{-1}(S_1)\ra S_1$ is elementary even if $d=2$
and $\psi^{-1}(U)\ra U$ is a quadric fibration.
Let $U_1$ denote the largest open subset of $S_1$ such that
$\psi^{-1}(U_1)\ra U_1$ is smooth.
Here we see by Theorem~\ref{Brauer-Severi} and 
Lemma~\ref{quadricbundle}
that $\psi^{-1}(U_1)\ra U_1$ is,
in the \'{e}tale topology, a $\mathbb{P}^{d}$-bundle
or a quadric fibration,
and on the space $\psi^{-1}(V)=V\times_{U_1}\psi^{-1}(U_1)$
over any small \'{e}tale open set $V$ of $U_1$
exists a line bundle $H_V$ such that 
the restriction of $H_V$ to any closed fiber of $\psi^{-1}(V)\ra V$
is isomorphic to $\mathcal{O}(1)$.
Since $\dim U_1=1$, Tsen's theorem
implies that $H^2({U_1}_{\textrm{et}},\mathbb{G}_m)=0$, 
where ${U_1}_{\textrm{et}}$ denotes $U_1$ with \'{e}tale topology.
(See, e.g., \cite[III, p.108]{milne}.)
Hence modifying a glueing, if necessary,
we can glue these $H_V$'s in the \'{e}tale topology.
Moreover it follows from \cite[III.4.9]{milne} that
$H^1(\psi^{-1}(U_1)_{\textrm{et}}, \mathbb{G}_m)
=H^1(\psi^{-1}(U_1)_{\textrm{Zar}}, \mathcal{O}^{\times})
=\Pic \psi^{-1}(U_1)$.
Hence there exists an algebraic line bundle $H$ on $\psi^{-1}(U_1)$
such that $H|_{F}=\mathcal{O}_{F}(1)$ for any closed fiber $F$
of $\psi^{-1}(U_1)\ra U_1$.
Since $\psi^{-1}(S_1)$ is smooth and $H$ is algebraic,
we can extend $H$ to a line bundle on $\psi^{-1}(S_1)$,
which we also denote by $H$ by abuse of notation.

Let $F'$ be an arbitrary closed fiber of $\psi^{-1}(S_1)\ra S_1$.
Then $F'$ is irreducible and reduced;
indeed, if $F'$ were decomposed as $F'_1+F'_2$, then
the Cartier divisor $F'_1$ of $\psi^{-1}(S_1)$
would satisfy
the condition $F'_1.C_0'=0$
for a curve $C_0'$ in a general fiber of $\psi^{-1}(S_1)\to S_1$,
so that $F'_1$ itself must be a fiber of 
$\psi^{-1}(S_1)\to S_1$,
a contradiction.
If $\psi^{-1}(U_1)\ra U_1$ is a $\mathbb{P}^{d}$-fibration,
then its closed fiber $(F,H|_F)$ has Fujita's delta genus
$\Delta (F,H|_F)=0$ and  degree $H|_F^{d}=1$,
and therefore 
$(F',H|_{F'})$ also has the same delta genus and degree,
so that $(F',H|_{F'})\cong (\mathbb{P}^{d}, \mathcal{O}(1))$.
Thus $\psi^{-1}(S_1)\ra S_1$ is a $\mathbb{P}^{d}$-bundle
in the Zariski topology.
If $\psi^{-1}(U_1)\ra U_1$ is a quadric fibration,
then the polarized variety $(F,H|_F)$ has Fujita's delta genus
$\Delta (F,H|_F)=0$ and  degree $H|_F^{d}=2$,
and therefore 
$(F',H|_{F'})\cong (\mathbb{Q}^{d}, \mathcal{O}(1))$,
which may be singular.
Since we have $K_{F'}+\tau\det \mathcal{E}|_{F'}=0$
by the adjunction formula and we have $\tau r=d$,
we infer that $\det \mathcal{E}|_{F'}=\mathcal{O}(r)$.
Therefore $\mathcal{E}|_{F'}$ is isomorphic to 
$\mathcal{O}(1)^{\oplus r}$
by \cite[Lemma 3.6.1]{w3}
(note that \cite[Lemma 3.6.1]{w3} also holds for
a singular hyperquadric $\mathbb{Q}$,
which can be shown easily by 
the induction on $\dim \Sing \mathbb{Q}$).
Hence if $d=n-1$,
$\mathcal{E}$ can be written
in the form $\psi^*\mathcal{E}'\otimes H$
for some vector bundle $\mathcal{E}'$ of rank $r$ on $S$.

Suppose that $\psi^{-1}(U)\ra U$ is a $\mathbb{P}^{d}$-fibration.
Since $\psi^{-1}(S_1)\ra S_1$ is a $\mathbb{P}^d$-bundle,
$\psi^{-1}(S_2)\ra S_2$ has at most 
a finite number of singular fibers.
Moreover $\psi^{-1}(S_2)\to S_2$ has no divisorial fibers
since $\psi^{-1}(S_2)\to S_2$ is elementary.
Hence $\psi^{-1}(S_2)\to S_2$ has equidimensional fibers.
Therefore we infer that 
every closed fiber $F'$ of $\psi^{-1}(S_2)\ra S_2$ is isomorphic to 
${\mathbb P}^{d-1}$ by the same argument as in 
\cite[\S 2.2, e-mail note of T. Fujita]{abw}.
Hence $S_2\subset U$ and $\codim (S\setminus U, S)\geq 3$.

Finally, since the intersection number 
$(H(\mathcal{E}(-H))|_{(\psi\circ\pi)^{-1}(t)})^{2d-1}$ is constant
for all $t\in U$,
we obtain the cases (i) and (ii) of the proposition.
\end{proof}

\section{The case $n-1<\tau r\leq n-1$ and $\dim S>0$.}\label{section=n-1/r}
Let $\psi :M\ra S$ be as in \S~\ref{overview}.
We will give a proof of Theorem~\ref{Main theorem}
in case $n-2<\tau r\leq n-1$ and $\dim S>0$.
By inequality~(\ref{bound of tau r}) in 
\S~\ref{overview},
we have
two cases:\\
1) $\psi$ is of fiber type and $1\leq \dim S\leq 2$;\\
2) $\psi$ is divisorial and contracts $E(R)$ to a point.\\
For the case 2), we have already seen 
in \S~\ref{prelimdiv}
that $S$ is smooth, that 
$\psi$ is the blowing up at a point,
that $\tau r=n-1$,
that 
$\mathcal{E}\cong \psi^*\mathcal{E}'\otimes \mathcal{O}(-E)$
for some ample vector bundle $\mathcal{E}'$ of rank $r$ on $S$,
and that $\tau (S,\det\mathcal{E}')r\leq n-1$.
Suppose that $\tau (S,\det\mathcal{E}')r= n-1$,
and set $M_1=S$ and $\mathcal{E}_1=\mathcal{E}'$.
We will determine the structure of $(M_1,\mathcal{E}_1)$
in \S~\ref{morebir},
to obtain the case (12) of Theorem~\ref{Main theorem}.

Suppose that we are in the case 1) and
that $\psi$ is of fiber type. 
If $\dim S=1$,
then we can apply Proposition~\ref{relative=n/r} to obtain
the cases (8), (9) and a part of (10) of 
Theorem~\ref{Main theorem}.
If $\dim S=2$,
it follows from Proposition~\ref{relative}
that every closed fiber $F$ of $\psi$ is isomorphic to ${\mathbb P}^{n-2}$ 
and that $\mathcal{E}|_{F}\cong \mathcal{O}(1)^{\oplus r}$.
This is the case (11) of Theorem~\ref{Main theorem}.
\begin{rmk}
Suppose that $\psi$ is of fiber type and that $\dim S=2$.
If $r$ and $n-1$ are coprime, there exist
integers $a$ and $b$ such that
$-a(n-1)+br=1$.
Hence $aK_M+b\det \mathcal{E}$ makes
$\psi$ in fact a ${\mathbb P}^{n-2}$-bundle
in the Zariski topology.
On the contrary, if $r$ and $n-1$ are not coprime,
does there always exist a pair of an ample vector bundle $\mathcal{E}$
and $\psi$ such that $\psi$ is a $\mathbb{P}^{n-2}$-bundle
not in the Zariski topology but in the \'{e}tale topology ?
If $r=n-2=2$, an example is shown in \cite{abw}
\end{rmk}

\subsection{More on divisorial case}\label{morebir}
Let $\psi:M\to M_1$, $E$, and 
$(M_1,\mathcal{E}_1)$ be as in \S~\ref{prelimdiv}.

Note that there exists no rational curve $l_1$ passing through
$\psi (E)$ such that $\mathcal{E}_1|_{l_1}$ has
an $\mathcal{O}_{l_1}(1)$ as a quotient line bundle
since $\mathcal{E}\cong \psi^*\mathcal{E}_1\otimes
\mathcal{O}(-E)$ and $\mathcal{E}$ is ample.

Suppose that $\tau (M_1,\det \mathcal{E}_1)r=n-1$
and that $(M_1,\mathcal{E}_1)$ 
has a divisorial elementary contraction.
The note above then implies that 
the exceptional divisor $E_1$ of $M_1$ does not contain $\psi(E)$. 

Suppose that $\tau (M_1,\det \mathcal{E}_1)r=n-1$
and that every elementary contraction of 
$(M_1,\mathcal{E}_1)$ is of fiber type.
Then for every point $p$ of $M_1$
there exists an extremal rational curve passing through $p$;
in particular, $M_1$ contains an extremal rational curve $l_1$
passing though $\psi (E)$.
Since we must have $\det\mathcal{E}_1.l_1\geq 2r$
by the note above,
we have 
$2(n-1)=2\tau (M_1,\det \mathcal{E}_1)r
\leq \tau (M_1,\det \mathcal{E}_1)\det\mathcal{E}_1.l_1
=-K_{M_1}.l_1\leq n+1$.
Hence $n\leq 3$.
Now we use the assumption that $\tau\geq 1$;
since $\tau\geq 1$, we have $n-1=\tau r\geq r$,
and thus $r\leq 2$.
Suppose that $n=3$.
The note above and the classification
given in the case 1) of \S~\ref{section=n-1/r}
implies that $\rho (M_1)$ cannot be bigger than one
(note that in the case (10) of Theorem~\ref{Main theorem}
we have $r=n-1$ if $\tau r=n-1$).
Hence $\rho (M_1)=1$.
Suppose moreover that $r=2$. Then $\tau =1$.
The note above and the classification
given in \S~\ref{LC1=2shin} and \ref{Slr1=n-1}
then implies that 
$(M_1,\mathcal{E}_1)\cong (\mathbb{P}^3,\mathcal{O}(2)^{\oplus 2})$.
Suppose moreover that $r=1$. 
Then $K_{M_1}+2\mathcal{E}_1=0$.
Hence $(M_1,\mathcal{E}_1)$ is a Del Pezzo $3$-fold 
of degree $\mathcal{E}_1^3=\mathcal{E}^3+1\geq 2$.
If $2\leq \mathcal{E}_1^3\leq 7$, 
we can deduce from the classification \cite{fb}, \cite{isk}
of Del Pezzo $3$-folds that
for each point $x\in M_1$ there is a line passing through $x$.
This contradicts the note above. Hence $\mathcal{E}_1^3=8$
and $(M_1,\mathcal{E}_1)\cong (\mathbb{P}^3,\mathcal{O}(2))$.
If $n=2$, 
since every elementary contraction of $M_1$ is of fiber type
by assumption,
$M_1$ is isomorphic to $\mathbb{P}^2$
or a scroll over a curve; 
if $M_1\cong \mathbb{P}^2$ then 
$\mathcal{E}_1\cong \mathcal{O}(3)$
and if $M_1$ is a scroll then
$\mathcal{E}_1|_F\cong \mathcal{O}_F(2)$
for every fiber $F$ of the scroll.

Summing up, we can repeat this process (divisorial case)
until strict inequality 
$\tau (M_k,\det \mathcal{E}_k)r<n-1$ holds,
unless one of the following holds:\\
1) $M_1\cong \mathbb{P}^3$ and $\mathcal{E}_1\cong \mathcal{O}(2)$
or $\mathcal{O}(2)^{\oplus 2}$;\\
2) $(M_k, \mathcal{E}_k)\cong (\mathbb{P}^2,\mathcal{O}(3))$ 
($1\leq k\leq 8$);\\
3) $M_k$ is a $\mathbb{P}^1$-bundle over a curve
and $\mathcal{E}_k|_F\cong \mathcal{O}_F(2)$
for every fiber $F$ of the projection.
\qed
\begin{rmk}
Suppose that $(M,\mathcal{E})$ has one of the 
structures of type (8), (9), (10), (11), (12) (d)
of Theorem~\ref{Main theorem};
here we only assume the structure of $(M,\mathcal{E})$
with $r\leq n-1$
and not the nef value to satisfy $\tau r=n-1$.
Then we have $\tau r=n-1$
unless $(M,\mathcal{E})$ is one of the following:
$(\mathbb{P}^1\times \mathbb{P}^1,
\mathcal{O}(1)\boxtimes \mathcal{O}(2))$;
$(\mathbb{P}^2\times \mathbb{P}^1,
\mathcal{O}(1)\boxtimes \mathcal{O}(1))$;
or
$(\mathbb{P}^2\times \mathbb{P}^1,
(\mathcal{O}(1)\boxtimes \mathcal{O}(1))^{\oplus 2})$.
Namely if we suppose moreover that $\tau r>n-1$ 
then $(M,\mathcal{E})$
is one of the pairs listed above.
\end{rmk}

\section{The case $\dim F-1\geq \tau r>\dim F-2$}\label{tau r=dimF-1}
Let $M$, $\mathcal{E}$, 
and $\psi :M\ra S$
be as in \S~\ref{overview}.

The following proposition gives a rough
classification of $(M,\mathcal{E})$
with $\psi$ of fiber type
in case $\dim F-1\geq \tau r>\dim F-2$ for a general fiber $F$ of $\psi$.
\begin{prop}\label{relative=n-1/r}
Suppose that $\tau\geq 1$,
that $\psi$ is of fiber type,
and that $\dim F-1\geq \tau r>\dim F-2$
for a general fiber $F$ of $\psi$.
Set $d=\dim F=\ulcorner \tau r\urcorner+1$.
Then we have one of the following:
\begin{enumerate}
\item $d$ is odd, $r=(d-1)/2$, and $(F, \mathcal{E}|_F)\cong 
(\mathbb{P}^{d}, \mathcal{O}(1)^{\oplus (r-1)}\oplus \mathcal{O}(2))$;
\item $F$ is a Del Pezzo manifold, 
and $\mathcal{E}|_F\cong A^{\oplus r}$ for the ample
line bundle $A$ on $F$ such that $K_F+(d-1)A=0$
(the Picard number $\rho (F)$ of $F$ is not necessarily one
in this case);
\item $\tau =1$, the Picard number $\rho (F)$ of $F$ is one, 
and $K_F+\det \mathcal{E}|_F=0$;
\item $r=2$, $\tau =1$, $F\cong \mathbb{P}^2\times \mathbb{P}^1$,
and $\mathcal{E}|_F=
(\mathcal{O}_{\mathbb{P}^2}(1)
\oplus 
\mathcal{O}_{\mathbb{P}^2}(2))\boxtimes \mathcal{O}_{\mathbb{P}^1}(1)$
or $T_{\mathbb{P}^2}\boxtimes \mathcal{O}_{\mathbb{P}^1}(1)$;
\item $d-1>\tau r>d-2$, the Picard number $\rho (F)$ of $F$ is one, 
and $K_F+\tau \det \mathcal{E}|_F=0$.
\end{enumerate}
\end{prop}
\begin{proof}
Since $d-1\geq \tau r>d-2$ , $d=\dim F$, and 
$K_F+\tau \det \mathcal{E}|_F=0$,
we will apply the classification of the case $n-1\geq \tau r>n-2$;
note that we have already proved Theorem~\ref{Main theorem}
under an additional assumption that $\tau r>n-2$.

Suppose that the Picard number $\rho (F)$ of $F$ is one.
If $d-1>\tau r>d-2$, we are in the case (5)
of the proposition.
(The precise classification of $(F,\mathcal{E}|_F)$
with $d-1>\tau r>d-2$
corresponds to the cases (4), (13), (14), 
(15) where $n\neq 2$, and (16)
of Theorem~\ref{Main theorem}.)
Suppose that $d-1=\tau r$.
If $\tau =1$,
we are in the case (3) of the proposition.
If $\tau>1$, we see, 
by the classification of the case $\tau r=n-1$
and by the assumption $\rho (F)=1$,
that we are in the cases (1) or (2) of the proposition.

Suppose that the Picard number $\rho (F)$ of $F$ is bigger than one.
Then $F$ has at least two extremal rays.
Moreover we see, by the classification of the case 
$n-1\geq \tau r>n-2$, that 
we may assume that $d-1=\tau r$.

Suppose that some of the extremal rays of $F$
corresponds to a birational contraction morphism.
Then $(F,\mathcal{E}|_F)$ lies in 
the case (12)
of Theorem~\ref{Main theorem};
we see, in particular, that 
the birational contraction of the ray is divisorial.
Denote by $E$ the exceptional divisor.
Since $K_F+\tau \det \mathcal{E}|_F$ is trivial,
$(F,\mathcal{E}|_F)$ does not lie in the sub-case (a) 
of the case (12) of Theorem~\ref{Main theorem}.
Suppose that $(F,\mathcal{E}|_F)$ is in 
the sub-case (b) of 
the case (12).
Then $F$ is a $\mathbb{P}^1$-bundle 
$\mathbb{P}(
\mathcal{O}_{\mathbb{P}^2}\oplus \mathcal{O}_{\mathbb{P}^2}(-1)
)$ over $\mathbb{P}^2$,
and $\mathcal{E}|_F\cong
H(
\mathcal{O}_{\mathbb{P}^2}(2)\oplus \mathcal{O}_{\mathbb{P}^2}(1)
)^{\oplus r}$.
Hence $(F,\mathcal{E}|_F)$ lies in the case (2)
of the proposition.
If $(F,\mathcal{E}|_F)$ is in 
the sub-case (c) or (d) of the case (12) of Theorem~\ref{Main theorem},
then $r=1$, $\tau =1$, and $F$ is a Del Pezzo surface
with $-K_F=\mathcal{E}|_F$.
Hence $(F,\mathcal{E}|_F)$ lies in the case (2)
of the proposition.

Suppose that each extremal ray of $F$ corresponds
to a contraction morphism of fiber type.
Since $\rho(F)\geq 2$ by assumption,
the image of any contraction morphism is positive dimensional.
Thus $(F,\mathcal{E}|_F)$ lies in either of 
the cases (11), (8), (9), or (10) of Theorem~\ref{Main theorem}.

Suppose that some extremal ray 
corresponds to a contraction of 
the type (11) of Theorem~\ref{Main theorem};
let $f:F\to S'$ be the contraction onto a smooth surface $S'$.
Then $f$ is a $\mathbb{P}^{d-2}$-fibration
and $d\geq 3$.

Suppose that some other extremal ray also
corresponds to a contraction of 
the case (11) of Theorem~\ref{Main theorem};
let $g:F\to S''$ be a $\mathbb{P}^{d-2}$-fibration
onto a smooth surface $S''$.
Since the restriction of $f$
to a fiber of $g$ is finite by Lemma~\ref{finite of yz},
we see that $d-2\leq 2$.
If $d-2=2$, then the restriction of $f$
to a fiber $\mathbb{P}^{2}$ of $g$ is finite surjective 
onto a smooth projective surface $S'$.
Therefore we infer that 
$S'\cong \mathbb{P}^2$
since a non-constant morphism from $\mathbb{P}^2$
is finite
so that $S'$ must be a minimal rational surface
which does not admit a surjective morphism onto a curve.
Now, by Lemma~\ref{isom},
we know that the restriction of $f$
to any fiber of $g$ is an isomorphism.
Therefore $F\cong \mathbb{P}^2\times \mathbb{P}^2$.
Moreover we have $\mathcal{E}|_F\cong 
(\mathcal{O}(1)\boxtimes \mathcal{O}(1))^{\oplus r}$,
where $1\leq r\leq \tau r=3$.
Setting $A=\mathcal{O}(1)\boxtimes \mathcal{O}(1)$,
we have $K_F+(d-1)A=0$
and $\mathcal{E}\cong A^{\oplus r}$.
Hence this is in the case (2) of the proposition.
If $d-2=1$, then
$r\leq \tau r=2$. 
If $r=2$, since $K_F+\det \mathcal{E}|_F=0$
and $F$ has two $\mathbb{P}^1$-fibrations over surfaces, 
we see, by \cite{w2},
that $F$ is a Del Pezzo manifold
and that $\mathcal{E}|_F\cong A^{\oplus r}$
where $A$ is a line bundle satisfying $K_F+(d-1)A=0$.
This description is also valid in case $r=1$.
Hence these are also in the case (2) of the proposition.

Suppose that some other extremal ray 
corresponds to a contraction of 
the case (8) of Theorem~\ref{Main theorem};
let $g:F\to C$ be the contraction
onto a smooth curve $C$.
Then a general fiber of $g$ is a smooth hyperquadric
of dimension $d-1$.
Since the restriction of $f$
to a fiber of $g$ is finite by Lemma~\ref{finite of yz},
we see that $d-1=2$.
Hence the restriction of $f$
to a general fiber $\mathbb{Q}^{2}$ of $g$ is finite surjective 
onto a smooth projective surface $S'$.
Therefore we infer that 
$S'\cong \mathbb{P}^2$ or $\mathbb{Q}^2$.
If $S'\cong \mathbb{P}^2$,
then we see,
by Lemma~\ref{isom},
that a general fiber 
$\mathbb{Q}^2$ of $g$ 
is isomorphic to $S'$
via the restriction of $f$.
This is a contradiction.
Therefore $S'\cong \mathbb{Q}^2$.
Now if $r=2$ it follows from \cite{w2}
that $F$ is a Del Pezzo manifold
and that $\mathcal{E}|_F\cong A^{\oplus 2}$
where $A$ is a line bundle satisfying $K_F+(d-1)A=0$.
This description is also valid in case $r=1$.
Hence these are in the case (2) of the proposition.

Suppose that some other extremal ray 
corresponds to a contraction of 
the case (9) or (10) of Theorem~\ref{Main theorem};
let $g: F\to C$ be the contraction onto a curve $C$.
Then a fiber of $g$ is isomorphic to $\mathbb{P}^{d-1}$.
Since the restriction of $f$
to a fiber of $g$ is finite by Lemma~\ref{finite of yz}
and thus surjective 
onto a smooth projective surface $S'$,
we see that $d-1=2$
and that 
$S'\cong \mathbb{P}^2$.
Now it follows from Lemma~\ref{isom}
that the restriction of $f$
to any fiber of $g$ is an isomorphism.
Therefore $F\cong \mathbb{P}^2\times \mathbb{P}^1$.
Note here that $\tau r=d-1=2$
and that 
$\tau \det \mathcal{E}|_F=-K_F=
\mathcal{O}_{\mathbb{P}^2}(3)\boxtimes \mathcal{O}_{\mathbb{P}^1}(2)$.
Hence $\tau =1$, $r=2$, and $\mathcal{E}|_F=
(\mathcal{O}_{\mathbb{P}^2}(1)
\oplus 
\mathcal{O}_{\mathbb{P}^2}(2))\boxtimes \mathcal{O}_{\mathbb{P}^1}(1)$
or $T_{\mathbb{P}^2}\boxtimes \mathcal{O}_{\mathbb{P}^1}(1)$.
This is the case (4) of the proposition.

Suppose that two of the extremal rays of $F$
correspond to contraction morphisms
of the cases (8), (9), or (10) of Theorem~\ref{Main theorem}.
Then we see that $d=2$
by Lemma~\ref{finite of yz}.
Since the case where one of extremal rays of $F$
corresponds to a birational morphism
is already done in the above argument,
we may assume that $F$ is minimal.
Therefore we infer that 
$F\cong \mathbb{P}^1\times \mathbb{P}^1$
and that $\mathcal{E}|_F=
\mathcal{O}(2)\boxtimes \mathcal{O}(2)$.
Thus we have $K_F+(d-1)A=0$
and $\mathcal{E}|_F=A^{\oplus r}$,
and we are in the case (2) of the proposition.
This completes the proof of the proposition.
\end{proof}

\section{The case $\tau r=n-2$
and $\psi$ of fiber type with $\dim S>0$}\label{section=n-2/r}
Let $\psi :M\ra S$
be as in \S~\ref{overview}.
We will give a proof of Theorem~\ref{Main theorem}
under the condition
that $\tau r=n-2$
and that $\psi$ is of fiber type with $\dim S>0$.
By inequality~(\ref{bound of tau r}) in \S~\ref{overview},
we infer that $1\leq \dim S\leq 3$.

Suppose that $\dim S=3$.
Then Proposition~\ref{relative} shows that 
$S$ has only finite number of Gorenstein rational singularities,
that, except for a finite number of singular fibers of $\psi$,
every smooth fiber $F$ of $\psi$ is isomorphic to $\mathbb{P}^{n-3}$,
and that $\mathcal{E}|_F\cong \mathcal{O}(1)^{\oplus r}$.
This is the case (24) of Theorem~\ref{Main theorem}.

Suppose that $\dim S=2$.
Then we can apply Proposition~\ref{relative=n/r};
note here that $S$ is smooth even if $d=2$ 
and $\psi^{-1}(U)\to U$ is a quadric fibration
as can be seen from its proof,
since $\psi\colon M=\psi^{-1}(S_2)\to S_2=S$ is, by assumption,
a contraction morphism of an extremal ray.
Therefore we obtain
the cases (21), (22) and (23) of Theorem~\ref{Main theorem}.

Suppose that $\dim S=1$.
Now we can apply Proposition~\ref{relative=n-1/r};
suppose that we are in the case (1) of Proposition~\ref{relative=n-1/r}.
Then we can show that $\psi$ is a projective space bundle.
Therefore we get 
a special case of the case (10) of Theorem~\ref{Main theorem}.
The case (2) of Proposition~\ref{relative=n-1/r}
corresponds to the case (19) of Theorem~\ref{Main theorem}.
and the cases (3) and (4) correspond to 
the case (20) of Theorem~\ref{Main theorem}.

\bibliographystyle{amsplain}
\bibliography{nefeprint.bbl}
\end{document}